\documentclass[10pt]{article}
\usepackage{amsmath,amssymb}
\usepackage{color}
\usepackage{mathrsfs}
\usepackage{mathptmx}
\usepackage{verbatim}
\usepackage{epsfig}
\usepackage{CJK}
\usepackage{bm}
\usepackage{amsfonts}
\usepackage{amsthm}
\input{epsf.sty}
\usepackage{epsfig}
\usepackage{setspace}
\def\<{\langle}
\def\>{\rangle}

\def\be{\begin{equation}}
\def\ee{\end{equation}}
\def\ba{\begin{array}}
\def\ea{\end{array}}

\newtheorem{thm}{Theorem}[section]
\newtheorem{cor}[thm]{Corollary}
\newtheorem{lem}[thm]{Lemma}

\newtheorem{defn}[thm]{Definition}
\newtheorem{rem}{Remark}

\numberwithin{equation}{section}

\def\be{\begin{equation}}
\def\ee{\end{equation}}
\def\br{\begin{eqnarray}}
\def\er{\end{eqnarray}}

\title{Stability of KAM tori for nonlinear Schr\"odinger equation
\footnote{Supported by  the NNSFC 11101059, NNSFC 11271076 and
NNSFC 11121101 .}}
\author{$\mbox{Hongzi \ Cong}^1$ \hspace{0.2 cm}
        $\mbox{Jianjun \ Liu}^2$ \hspace{0.2 cm}
        $\mbox{Xiaoping \ Yuan}^3$ \hspace{0.2 cm}
        \\
${}^{1}$ School of Mathematical Sciences, Dalian University of
Technology \\
${}^{2}$ School of Mathematical Sciences, Sichuan University\\
${}^{3}$ School of Mathematical Sciences, Fudan University}
\begin{document}
\maketitle
 {\bf Abstract}: {\it  It is proved that the KAM tori (thus quasi-periodic solutions)  are long time stable for nonlinear Schr\"odinger equation.}

  \maketitle
\section{Introduction}

Since Kuksin's work \cite{K3} in 1987, the infinite dimensional KAM theory  has seen enormous  progress with application to partial differential equations (PDEs). There are too many references to list all of them. Here we refer to two books \cite{K1} and
\cite{Bour6} and two survey papers \cite{Bour5} and \cite{Kuk3}.

As an object to which the infinite dimensional KAM theory applies, one of typical models is nonlinear Schr\"odinger equation (NLS)
\be\label{y1} \sqrt{-1}\, u_t-\triangle u+V(x,\xi)u +|u|^2u+h.o.t.=0,\ee
subject to Dirichlet boundary condition or periodic boundary condition $x\in\mathbb T^d$, where the integer $d\ge 1$ is the spatial dimension of NLS.

{\it Case 1. $d=1$.} In 1993, Kuksin \cite{K1} proved that (\ref{y1}) possesses many invariant tori around the origin $u=0$ (thus quasi-periodic solutions of small initial values) for ``most" parameter vector $\xi$ when the potential $V$ depends on $\xi$ in some non-degenerate way. This kind of invariant tori obtained by KAM theory are usually called KAM tori. In 1996, Kuksin-P$\ddot{\mbox{o}}$schel \cite{K-P-2} further proved that the potential $V=V(x,\xi)$ can be replaced by a fixed constant potential $V\equiv C$.
All those results are obtained by KAM theory which involves the so-called second Melnikov conditions. By advantage of them the linearized equation along the KAM tori can be reduced to a linear equation of constant coefficient, thus the obtained KAM tori (thus quasi-periodic solutions) are linearly stable.

{\it Case 2. $d>1$.} This case is is significantly more complicated, since the second Melnikov conditions are violated seriously by multiplicity of the eigenvalues of $-\triangle$.
In his series of papers \cite{Bour1}-\cite{Bour6}, Bourgain developed a profound approach which was originally proposed by Craig-Wayne \cite{C-W}. It was proved by Bourgain \cite{B3}(1998) and \cite{Bour6}(2005) that (\ref{y1}) has many KAM tori for most $\xi$ when $V(x,\xi)$ is replaced by
   Fourier multiplier $M_\xi$. This approach which is today called C-W-B method does not involve the second Melnikov conditions. So is successfully overcome the difficulty that the second Melnikov conditions are seriously violated. Just because of lack of the second Melnikov conditions, the linear stability of the obtained KAM tori can not be derived incidentally. Recently Eliasson-Kuksin \cite{E-K} developed the classic KAM theory which involves the second Melnikov conditions and proved that the obtained KAM are linearly stable.

We call a KAM tori to be linearly stable, if the equilibrium point of the linearized equation along the KAM tori is Lyapunov stable and this kind stability can preserve undergoing a linear perturbation. Usually, a perturbation of some equation or system is {\it nonlinear}. So the nonlinear stability in some sense should be more interesting. In 2007, Eliasson \cite{Elia2} proposed to study the  nonlinear stability of the obtained KAM tori for (\ref{y1}).

In the present paper, we will prove that the obtained KAM tori of (\ref{y1}) with $d=1$ are of long time stability (a kind of nonlinear stability). More exactly, we have the following theorem:
\begin{thm}\label{T4}Consider the nonlinear
Schr$\ddot{\mbox{o}}$dinger equation
\begin{equation}\label{26}
\sqrt{-1}\, u_{t}=u_{xx}-M_{\xi}u+\varepsilon|u|^2u,
\end{equation}
subject to Dirichlet boundary conditions $u(t,0)=u(t,\pi)=0$, where
$M_{\xi}$ is a real
Fourier multiplier,
\begin{equation*}\label{091913}
M_{\xi} \sin jx=\xi_j\sin jx,\qquad \xi=(\xi_j)_{j\geq1}
\end{equation*}
and
\begin{equation*}
\xi\in\Pi:=\left\{\tilde\xi=(\tilde\xi_j)_{j\geq1}|\ \tilde{\xi}_j\in [1,2]/j,\ j\geq 1\right\}\subset\mathbb{R}^{\mathbb{N}}.
\end{equation*}For a sufficiently small
$\varepsilon>0$ depending on $n$ and $p$, there exists
a subset $\tilde\Pi\subset\Pi$, such that
for each $\xi\in\tilde\Pi$ equation (\ref{26}) possesses an $n$-dimensional KAM
torus $\mathcal{T}_{\xi}$ in Sobolev space $H^p_0([0,\pi])$.
Moreover, the KAM torus $\mathcal{T}_{\xi}$ of equation
(\ref{26}) is stable in long time, that is, for arbitrarily given $\mathcal M$ with
$0\leq \mathcal{M}\leq C(\varepsilon)$ (where $C(\varepsilon)$ is a constant depending on $\varepsilon$ and $C(\varepsilon)\rightarrow\infty$ as $\varepsilon\rightarrow0$) and any $p\geq 8(\mathcal{M}+7)^{4}+1$, there exists a
small positive $\delta_0$ depending on $n,p$ and
$\mathcal{M}$, such that for any
$0<\delta<\delta_0$ and any solution $u(t,x)$ of equation (\ref{26})
with the initial datum satisfying \color{black}
$${d}_{H^p_0[0,\pi]}(u(0,x),\mathcal{T}_{\xi}):=\inf_{w\in\mathcal{T}_\xi}||u(0,x)-w||_{H^p_0[0,\pi]}\leq \delta,$$
then
\begin{equation*}{d}_{H^p_0[0,\pi]}(u(t,x),\mathcal{T}_{\xi}):=\inf_{w\in\mathcal{T}_\xi}||u(t,x)-w||_{H^p_0[0,\pi]}\leq
2\delta,\qquad \mbox{for all} \ |t|\leq
{\delta}^{-\mathcal{M}}.
\end{equation*}
\end{thm}

Indeed, recently it has been proved by  Colliander-Keel-Staffilani-Takaota-Tao \cite{CKSTT} that  there exists Arnold diffusion in NLS equation. Thus, one should not expect to extend the stability $|t|\leq
{\delta}^{-\mathcal{M}}$ to $t\in (-\infty,+\infty)$ for all initial values.

There has been a lot of works on the long time stability of the equilibrium point $u=0$ and some approximate invariant tori for partial differential equations. See \cite{Bam1}-\cite{BN}, \cite{B96}, \cite{B6}, \cite{DS}, for example. Especially, we incorporate some important ideas on tame property from Bambusi-Gr$\acute{\mbox{e}}$bert's work \cite{BG}. See \S 3 for more remarks.

\color{black}
The rest of the present paper is organized as follows.\\
 In \S 2, we give some basic notations and an abstract theorem on the long time stability. \\
 In \S 3, we
will give some discussions and ideas of the proof of the abstract results. \\
In \S 4, we will discuss some properties of $p$-tame
property.\\
 In \S 5, we will discuss $p$-tame property of the solution of homological equation.\\
 In \S 6, we will construct a
partial normal form of order $\mathcal{M}+2$ in the
$\delta$-neighborhood of the KAM tori and show that the KAM tori are
stable in a long time (Theorem \ref{T3}).\\
 In \S 7, it is shown the
existence and long time stability of KAM tori (i.e.
quasi-periodic solutions) for the nonlinear
Schr$\ddot{\mbox{o}}$dinger equation (\ref{26}) according to the
above theorems and corollaries (Theorem \ref{T4}).\\
 In \S 8, it will
be given some technical lemmas.
In \S 9, the details of proof of Theorem {\ref{021102}} and Corollary \ref{081914} will be given, since it is a standard KAM proof, based on the results in Section 4 and Section 5.
\section*{Acknowledgements}{\it

In the Autumn of 2007,
Professor H. Eliasson gave a series of lectures on KAM theory for
Hamiltonian PDEs in Fudan University. In those lectures, he proposed
to study the normal form in the neighborhood of the invariant tori
and the nonlinear stability of the invariant tori.  The
authors are heartily grateful to Professor Eliasson.
}

\section{Some notations and the abstract results}
\subsection{Some notations}
To finish the proof of Theorem \ref{T4}, several abstract theorems will be given. To this end, we will introduce some notations
firstly. Given positive integers $n$ and $p$, by
$\mathbb{T}^n=\mathbb{C}^n/2\pi\mathbb{Z}^n$ denote the usual
$n$-dimensional torus and let  $\ell_{p}^2$ be the Hilbert space of
all
 complex sequences $w=(w_1,w_2,\dots)$ with
\begin{equation*}
||w||_{p}^2=\sum_{j\geq1}|w_j|^2j^{2p}<\infty.
\end{equation*}
Introduce an infinite dimensional symplectic phase space
\begin{equation*}
(x,y,q,\bar{q})\in\mathcal{P}^p:=\mathbb{T}^n\times\mathbb{C}^n\times
\ell^2_{p}\times \ell_{p}^2
\end{equation*} with the usual symplectic structure
\[dy\wedge dx+\sqrt{-1}dq\wedge d\bar q.\]
Given a subset $\Pi\subset\mathbb R^{\mathbb N}$ with positive
 measure in some sense (for example, in the sense of Gauss or Kolmogorov),
 here $\Pi$ will be regarded a parameter set. Let $N(y,q,\bar q;\xi)$ be an integrable Hamiltonian which depends on parameter $\xi\in\Pi$ and is of the form
 $$N(y,q,\bar
q;\xi)=\sum_{j=1}^{n}\omega_j(\xi)y_j+\sum_{j\geq1}\Omega_j(\xi)q_j\bar
q_j,$$where $\omega(\xi)=(\omega_1(\xi),\dots,\omega_n(\xi))$ is
called tangent frequency  and $
\Omega(\xi)=(\Omega_1(\xi),\Omega_2(\xi),\dots)$ is called normal
frequency. With the symplectic structure mentioned above, the motion
equation of $N(y,q,\bar q;\xi)$ is
\begin{equation}\label{0004}
\dot x=\omega(\xi),\quad \dot y=0,\quad \dot
q_j=\sqrt{-1}\Omega_j(\xi)q_j,\quad \dot{\bar
q}_j=-\sqrt{-1}\Omega_j(\xi)\bar q_j,\qquad j\geq 1.
\end{equation}
It is clear that, for each $\xi\in \Pi$, $(x(t),y(t),q(t),\bar
q(t))=(\omega(\xi)t,0,0,0)$ is a quasi-periodic solution to equation
(\ref{0004}) with rotational frequency $\omega(\xi)$. Moreover, let
$\hat{\mathbb{T}}^n=\mathbb{R}^n/2\pi\mathbb{Z}^n$, and
$$\mathcal{T}_0=\hat{\mathbb{T}}^n\times\{y=0\}\times\{q=0\}\times\{\bar
q=0\}.$$ Then $\mathcal{T}_0$ is an $n$-dimensional invariant torus with frequency
${\omega}(\xi)$ for equation (\ref{0004}).

Now consider a perturbation of the integrable Hamiltonian
$N(y,q,\bar q;\xi)$:
$$H(x,y,q,\bar q;\xi)=N(y,q,\bar q;\xi)+ R(x,y,q,\bar q;\xi),$$ where the perturbation $R(x,y,q,\bar q;\xi)$
 depends on the parameter $\xi\in\Pi$
  and the variable $(x,y,q,\bar q)\in \mathcal{P}^p$, and $R(x,y,q,\bar q;\xi)$
  is of small size in some sense ($p$-tame norm) which will be defined in the following steps.
\begin{defn}\label{090501} Let
\begin{equation*}
D(s)=\left\{x\in\mathbb{T}^n| \ ||\mbox{Im}\ x||\leq s\right\},
\end{equation*}where
$||\cdot||$ denotes the sup-norm for
complex vectors in $\mathbb C^n$ or $\mathbb{C}^{\mathbb{N}}$.
Consider a function $W(x;\xi):D(s)\times\Pi\rightarrow \mathbb{C}$
is analytic about the variable $x\in D(s)$ and $C^1$-smooth in the
Whitney's sense\footnote{In the whole of this paper, the derivatives
with respect to the parameter $\xi\in\Pi$ are understood in the
sense of Whitney.} about the parameter $\xi\in\Pi$ with the Fourier
series
$$W(x;\xi)=\sum_{k\in \mathbb Z^n}\widehat
W(k;\xi)e^{\sqrt{-1}\langle k,x\rangle},$$where $$\widehat
W(k;\xi):=\frac1{(2\pi)^n}\int_{\hat{\mathbb{T}}^n}W(x;\xi)e^{-\sqrt{-1}\langle
k,x\rangle}dx$$ is the $k$-th Fourier coefficient of the function
$W(x;\xi)$, and $\langle\cdot,\cdot\rangle$ denotes the usual inner
product, i.e.
\begin{equation*}
\langle k,x\rangle=\sum_{j=1}^nk_jx_j.
\end{equation*}
Then define the norm
\begin{equation}\label{092704}
||W||_{D(s)\times\Pi}=\sup_{\xi\in\Pi,j\geq1}{\sum_{k\in\mathbb{Z}^n}
\left(|\widehat{W}(k;\xi)|+|\partial_{\xi_j}
\widehat{W}(k;\xi)|\right)e^{|k|s}}.
\end{equation}
\end{defn}
\begin{defn}\label{060602}Let
\begin{equation*}
D(s,r)=\{(x,y)\in\mathbb{T}^n\times\mathbb{C}^n|\ ||\mbox{Im}\
x||\leq s,\ ||y||\leq r^2\}.
\end{equation*}
Consider a function $W(x,y;\xi):D(s,r)\times\Pi\rightarrow
\mathbb{C}$ is analytic about the variable $(x,y)\in D(s,r)$ and
$C^1$-smooth about the parameter $\xi\in\Pi$ with the following form
$$W(x,y;\xi)=\sum_{\alpha\in\mathbb{N}^n}W^{\alpha}(x;\xi)y^{\alpha}.$$
Define the norm
\begin{equation}\label{092703}
||W||_{D(s,r)\times\Pi}=\sum_{\alpha\in\mathbb{N}^n}
||W^{\alpha}(x;\xi)||_{D(s)\times\Pi}\;r^{2|\alpha|},
\end{equation}
where $|\alpha|=\sum_{j=1}^n|\alpha_j|$. In
this paper, always by $|\cdot|$ denotes $1$-norm for complex vectors in
$\mathbb{C}^n$ or $\mathbb{C}^{\mathbb{N}}$.
\end{defn}

 \begin{defn}\label{021002}
Introduce the complex $\mathcal T_0$-neighborhoods
\begin{equation*}
D(s,r_1,r_2)=\left\{(x,y,q,\bar q)\in \mathcal
P^p\left|\right.||\mbox{Im}\ x||\leq s,||y||\leq r_1^2,||q||_p+||\bar
q||_p\leq r_2\right\}.
\end{equation*} Let $r_1=r_2=r$. Consider a function
$W(x,y,q,\bar q;\xi):D(s,r,r)\times\Pi\rightarrow\mathbb{C}$ is
analytic about the variable $(x,y,q,\bar q)\in D(s,r,r)$ and
$C^1$-smooth about the parameter $\xi\in\Pi$ with the following form
$$W(x,y,q,\bar
q;\xi)=\sum_{\alpha\in\mathbb{N}^n,\beta,\gamma\in\mathbb{N}^{\mathbb{N}}}
W^{\alpha\beta\gamma}(x;\xi)y^{\alpha}q^{\beta}\bar q^{\gamma}.$$
Define the modulus $\lfloor W\rceil_{D(s,r)\times\Pi}(q,\bar q)$ of
$W(x,y,q,\bar q;\xi)$ by
\begin{equation}\label{081602}\lfloor W\rceil_{D(s,r)\times\Pi}(q,\bar q):=
\sum_{\beta,\gamma\in\mathbb{N}^{\mathbb{N}}}||W^{\beta\gamma}(x,y;\xi)||_{D(s,r)\times\Pi}\;
q^{\beta}\bar q^{\gamma},\end{equation}where
$$W^{\beta\gamma}(x,y;\xi)=\sum_{\alpha\in\mathbb{N}^n}W^{\alpha\beta\gamma}(x;\xi)y^{\alpha}.$$
\end{defn}Denote $z=(q,\bar{q}),$ i.e.$$
z=(z_j)_{j\in\bar{\mathbb{Z}}}\in\ell_p^2\times\ell_p^2,\qquad
\bar{\mathbb{Z}}=\mathbb{Z}-\{0\},$$ where $$z_{-j}=q_j,\qquad
z_j=\bar{q}_j,\qquad \qquad j\geq1.$$\footnote{In the whole of this
paper, we will use the notation $z=(q,\bar q)$.}Here we consider $q$
and $\bar{q}$ as two independent variables and define
$||z||_p=||q||_p+||\bar{q}||_p.$ For a homogeneous polynomial $W(z)$
of degree $h>0$, it is naturally associated with a symmetric
$h$-linear form $\widetilde W(z^{(1)},\dots,z^{(h)})$ such that $
\widetilde W(z,\dots,z)=W(z).$ More precisely, given a monomial
\begin{equation*}
W(z)=W^{\beta}z^{\beta}=W^{\beta}z_{j_1}\cdots z_{j_h},
\end{equation*}where $\beta=(\dots,\beta_{-2},\beta_{-1},\beta_1,\beta_2,\dots)\in\mathbb{N}^{\bar{\mathbb{Z}}}$ and
$|\beta|=\sum_{|j|\geq1}\beta_j=h,$ the symmetric $h$-linear form
$\widetilde W(z^{(1)},\dots,z^{(h)})$ is defined by
\begin{equation}\label{080801}
\widetilde
W(z^{(1)},\dots,z^{(h)})=\widetilde{W^{\beta}z^{\beta}}=\frac1{h!}\sum_{\tau_h}W^{\beta}z_{j_1}^{(\tau_h(1))}\cdots
 z_{j_h}^{(\tau_h(h))},
\end{equation}
where $\tau_h$ is an $h$-permutation. Now assume $W(z)$ is a
homogeneous polynomial of degree $h$ and of the  form $
W(z)=\sum_{|\beta|=h}W^{\beta}z^{\beta}, $  then define
$\widetilde W(z^{(1)},\dots,z^{(h)})$ by
\begin{equation}\label{080802}\widetilde
W(z^{(1)},\dots,z^{(h)})=\sum_{|\beta|=h}\widetilde{W^{\beta}z^{\beta}}.
\end{equation}

Basing on the above notations,  we will define $p$-tame norm of a
Hamiltonian vector field. Firstly, consider a Hamiltonian
\begin{equation}
\label{011601}W_h(x,y,z;\xi)=\sum_{\alpha\in\mathbb{N}^n,\beta\in\mathbb{N}^{\bar{\mathbb{Z}}},|\beta|=h}
W_h^{\alpha\beta}(x;\xi)y^{\alpha}z^{\beta},
\end{equation} where the
modulus of $W_h(x,y,z;\xi)$ is a homogeneous polynomial about $z$ of degree $h$ and
$W_h(x,y,z;\xi)$ itself is analytic about the variable $(x,y,z)\in D(s,r,r)$ and
$C^1$-smooth about the parameter $\xi\in\Pi$.  To simplify the
notation, we rewrite $W_h(x,y,z;\xi)$ as $W(x,y,z;\xi)$. By a little abuse of notation, let $W_z=(\sqrt{-1}{W_{\bar q}},-\sqrt{-1}W_{{q}})$
here and later. Notice the Hamiltonian vector field $X_W$ of
$W(x,y,z;\xi)$ is $(W_y,-W_x,W_z)$.  Denote
\begin{eqnarray}\label{091401}||(z^{h})||_{p,1}:=
\frac{1}{h}\sum_{j=1}^{h}||z^{(1)}||_1\cdots||z^{(j-1)}||_1||
z^{(j)}||_p||z^{(j+1)}||_1\cdots||z^{(h)}||_1.
\end{eqnarray}
\begin{defn}\label{071803}
In the normal direction of the Hamiltonian vector field $X_W$,
define the $p$-tame operator norm of $W_z$ by,
\begin{eqnarray}\label{091402}
|||{W_z}|||_{p,D(s,r)\times\Pi}^{T}: =\sup_{0\neq z^{(j)}\in \ell
^2_p\times\ell^2_p,1\leq j\leq h-1}
\frac{||\widetilde{\lfloor{W_z}\rceil}_{D(s,r)\times\Pi}(z^{(1)},\dots,z^{(h-1)})||_
{p}} {||(z^{h-1})||_{p,1}},\end{eqnarray} and
define the $p$-tame norm of $W_z$ by
\begin{eqnarray}\label{091403}
|||{W_z}|||_{p,D(s,r,r)\times\Pi}^T=\max
\left\{|||{W_z}|||_{p,D(s,r)\times\Pi}^T,|||{W_z}|||_{1,D(s,r)\times\Pi}^T\right\}r^{h-1}.
\end{eqnarray}\end{defn}
\begin{defn}\label{071802}
In the tangent direction of the Hamiltonian vector field $X_W$,
define the operator norm of $W_u$ ($u=x\ \mbox{or}\ y$) by,
\begin{eqnarray}\label{091404} |||W_{u}|||_{D(s,r)\times\Pi}:=\sup_{0\neq
z^{(j)}\in \ell ^2_1\times\ell^2_1,1\leq j\leq h}
\frac{||\widetilde{\lfloor{W_{u}}\rceil}_{D(s,r)\times\Pi}(z^{(1)},\dots,z^{(h)})||}
{||(z^{h})||_{1,1}},
\end{eqnarray}
and define the norm of $W_u$ ($u=x\ \mbox{or}\ y$) by
\begin{eqnarray}\label{091405}
|||{W_u}|||_{D(s,r,r)\times\Pi}:=|||{W_u}|||_{D(s,r)\times\Pi}r^{h}.
\end{eqnarray}
\begin{rem}Note that the dimensional of the tangent direction is
finite, so there is no so-called $p$-tame property. But $||W_u||_{D(s,r)\times\Pi}$ is required as a bounded map from $\ell^2_1\times\ell^2_1$ to $\mathbb{C}^n$ to guarantee the persistence of $p$-tame property under Poisson bracket.
\end{rem}
\end{defn}
\begin{defn}\label{051703} Define the $p$-tame norm of the Hamiltonian vector field $X_W$ as
follows,
\begin{eqnarray*}
|||X_W|||_{p,D(s,r,r)\times\Pi}^T =|||{W_y}|||_{D(s,r,r)\times\Pi}
+\frac1{r^2}|||{W_x}|||_{D(s,r,r)\times\Pi}
+\frac1r|||W_z|||_{p,D(s,r,r)\times \Pi}^T.
\end{eqnarray*}\end{defn}
\begin{defn}\label{080204} Consider a Hamiltonian $W(x,y,z;\xi)=\sum_{h\geq
0}W_{h}(x,y,z;\xi)$ is analytic about the variable $(x,y,z)\in
D(s,r,r)$ and $C^1$-smooth about the parameter $\xi\in\Pi$, where
$$W_{h}(x,y,z;\xi)=\sum_{\alpha\in\mathbb{N}^n,\beta\in\mathbb{N}^{\bar{\mathbb{Z}}},|\beta|=h}W_{h}^{\alpha\beta}(x;\xi)
y^{\alpha}z^{\beta}.$$ Then define $p$-tame norm of the Hamiltonian
vector field $X_W$ by
\begin{equation}\label{102201}
|||X_W|||_{p,D(s,r,r)\times\Pi}^T:=\sum_{h\geq 0}
|||X_{W_{h}}|||_{p,D(s,r,r)\times\Pi}^T. \end{equation} We say that
a Hamiltonian vector field $X_W$ (or a Hamiltonian $W(x,y,z;\xi)$)
has ${p}$-tame property on the domain $D(s,r,r)\times\Pi$, if and
only if $ |||X_W|||_{p,D(s,r,r)\times\Pi}^{T}<\infty. $
\end{defn}


\subsection{The abstract results}

 Now we have the following theorems:
\begin{thm}\label{T1}(Normal form of order 2)
Consider a perturbation of the integrable Hamiltonian
\begin{equation}\label{100803} H(x,y,q,\bar q;\xi)=N(y,q,\bar q;\xi)+R(x,y,q,\bar
q;\xi)
\end{equation}defined on the domain $D(s_0,r_0,r_0)\times\Pi$ with
$s_0,r_0\in(0,1]$, where
$$N(y,q,\bar q;\xi)=\sum_{j=1}^{n}\omega_j(\xi)y_j+\sum_{j\geq1}\Omega_j(\xi)q_j\bar
q_j$$ is a family of parameter dependent integrable Hamiltonian and
$$R(x,y,q,\bar q;\xi)=\sum_{\alpha\in\mathbb{N}^n,\beta,\gamma\in\mathbb{N}^{\mathbb{N}}}R^{\alpha\beta\gamma}(x;\xi)y^{\alpha}q^{\beta}\bar q^{\gamma}$$ is the perturbation. Suppose the tangent
frequency and normal frequency satisfy the following assumptions:

$Assumption\ A:Frequency\ Asymptotics.$  There exists absolute
constants $c_1, c_2>0$ such that
\begin{equation}\label{006}
|\Omega_i(\xi)-\Omega_j(\xi)|\geq c_1|i-j|(i+j),
\end{equation}
and \begin{equation}\label{091212} |\Omega_j(\xi)|\leq c_2j^2,
\end{equation} for all integers $i,j\ge 1$
uniformly on $\xi\in\Pi$;

$Assumption\ B: Twist \ conditions.$
\begin{eqnarray}\label{081902}
\partial_{\xi_j}\omega_i(\xi)=\delta_{ji}, \quad
\partial_{\xi_j}\Omega_{j'}(\xi)=\delta_{j(n+j')},\qquad 1\leq i\leq n,\ j,j'\ge
1.
\end{eqnarray}
The perturbation $R(x,y,q,\bar q;\xi)$ has $p$-tame property on the
domain $D(s_0,r_0,r_0)\times\Pi$ and satisfies the small assumption:
\begin{equation*}
 \varepsilon:=|||X_{R}|||^T_{p,D(s_0,r_0,r_0)\times\Pi}\leq
 \eta^{12}\epsilon, \qquad \mbox{ for some }\ \eta\in(0,1),
\end{equation*}where $\epsilon$ is a positive constant
depending on $s_0,r_0$ and $n$. Then there exists a subset
$\Pi_{\eta}\subset\Pi$ with the estimate
\begin{equation*}
\mbox{Meas}\ \Pi_{\eta}\geq(\mbox{Meas}\ \Pi)(1-O(\eta)).
\end{equation*}For each $\xi\in\Pi_{\eta}$, there is a symplectic map
$$\Psi: D(s_0/2,r_0/2,r_0/2)\rightarrow D(s_0,r_0,r_0),$$ such that
\begin{equation}\label{081601}\breve H(x,y,q,\bar q;\xi):=H\circ\Psi=
\breve N(y,q,\bar q;\xi)+ \breve R(x,y,q,\bar q;\xi),
\end{equation}where
\begin{equation}\label{082903}
\breve N(y,q,\bar
q;\xi)=\sum_{j=1}^{n}\breve{\omega}_j(\xi)y_j+\sum_{j\geq1}\breve{\Omega}_j(\xi)q_j\bar
q_j
\end{equation}
and
\begin{equation}\label{082902}
\breve R(x,y,q,\bar
q;\xi)=\sum_{\alpha\in\mathbb{N}^n,\beta,\gamma\in\mathbb{N}^{\mathbb{N}},2|\alpha|+|\beta|+|\gamma|\geq
3}\breve
R^{\alpha\beta\gamma}(x;\xi)y^{\alpha}q^{\beta}\bar{q}^{\gamma}.
\end{equation}
Moreover, the following estimates hold:\\
(1) for each $\xi\in\Pi_{\eta},$ the symplectic map
$\Psi:D(s_0/2,r_0/2,r_0/2)\rightarrow D(s_0,r_0,r_0)$
\color{black}satisfies
\begin{equation}\label{080101}
||\Psi-id||_{p,D(s_0/2,r_0/2,r_0/2)}\leq c\eta^6\epsilon,
\end{equation}
where
\begin{equation}\label{081102}
||\Psi-id||_{p,D(s_0/2,r_0/2,r_0/2)}=\sup_{w\in
D(s_0/2,r_0/2,r_0/2)}||(\Psi-id)w||_{\mathcal{P}^p,D(s_0,r_0,r_0)},
\end{equation}
and
\begin{equation}\label{032601}
||\tilde w||_{\mathcal{P}^p,D(s_0,r_0,r_0)}=||\tilde w_x||+\frac1{r_0^2}||\tilde w_y||+\frac1{r_0}||\tilde w_q||_p+\frac1{r_0}||\tilde w_{\bar{q}}||_p
 \end{equation}for each $\tilde w=(\tilde w_x,\tilde w_y,\tilde w_q,\tilde w_{\bar q})\in D(s_0,r_0,r_0)$;
 moreover,
\begin{equation}\label{080102}
|||D\Psi-Id|||_{p,D(s_0/2,r_0/2,r_0/2)}\leq c\eta^6\epsilon,
\end{equation}
where on the left-hand side hand we use the operator
norm\footnote{where $id$ denotes the identity map from $\mathcal
P^p\to\mathcal P^p$ and $Id$ denotes its tangent map. }
\begin{equation*}
|||D\Psi-Id|||_{p,D(s_0/2,r_0/2,r_0/2)}=\sup_{0\neq w\in
D(s_0/2,r_0/2,r_0/2) }\frac{||(D\Psi-Id) w||_{\mathcal{P}^p,D(s_0,r_0,r_0) }}{||
w||_{\mathcal{P}^p,D(s_0/2,r_0/2,r_0/2)}};
\end{equation*}
\\
\color{black} (2) the frequencies $\breve \omega(\xi)$ and $\breve
\Omega(\xi )$ satisfy
\begin{equation}\label{080103}
||\breve\omega(\xi)-\omega(\xi)||+\sup_{j\geq
1}||\partial_{\xi_j}(\breve\omega(\xi)-\omega(\xi))||\leq
{c\eta^8\epsilon},\end{equation}and
\begin{equation}\label{080104}
||\breve\Omega(\xi)-\Omega(\xi)||+\sup_{j\geq
1}||\partial_{\xi_j}(\breve\Omega(\xi)-\Omega(\xi))||\leq
{c{\eta^8}\epsilon};
\end{equation}
\\
(3) the Hamiltonian vector field $X_{\breve R}$ of the new perturbed
Hamiltonian $\breve R(x,y,q,\bar q;\xi)$ satisfies
\begin{equation}\label{080105}
 |||X_{\breve R}|||^T_{p,D(s_0/2,r_0/2,r_0/2)\times\Pi_{\eta}}\leq
 \varepsilon(1+{c\eta^6\epsilon}),
\end{equation}where $c>0$ is a constant depending on $s_0,r_0$ and $n$.
\end{thm}
\begin{cor}\label{L1}(The existence and the stability of KAM tori)
Consider the Hamiltonian$${\breve H}(x,y,q,\bar q;\xi)=\breve
N(y,q,\bar q;\xi)+\breve R(x,y,q,\bar q;\xi)$$ obtained in Theorem
\ref{T1}. For each $\xi\in\Pi_{\eta}$, there is an analytic
embedding invariant torus
$\mathcal{T}_{0}=\hat{\mathbb{T}}^n\times\{y=0\}\times\{q=0\}\times\{\bar
q=0\}$ with frequency $\breve{\omega}(\xi)$ for the Hamiltonian
$\breve{H}(x,y,q,\bar q;\xi)$, and $
{\mathcal{T}}:=\Psi^{-1}\mathcal{T}_{0}$ is an analytic embedding
invariant torus (i.e. so-called KAM torus) for the original
Hamiltonian $H(x,y,q,\bar q;\xi)$.

Moreover, given any small positive $\delta<r_0/10$, if $w(t)$ is a
solution of Hamiltonian vector field $X_{{H}}$ with the initial
datum $w(0)=(w_x(0),w_y(0), w_q(0), w_{\bar q}(0))$ satisfying
\begin{equation*}
d_{p}(w(0),\mathcal{T})\leq\delta,
\end{equation*}
then
\begin{equation*}
d_{p}(w(t),\mathcal{T})\leq 2\delta, \qquad\mbox{for all}\ |t|\leq
\delta^ {-1},
\end{equation*}
where the distance $d_p(w,v)$ between any two points
$$w=(w_x,w_y,w_q,w_{\bar q}),v=(v_x,v_y,v_q,v_{\bar q})\in D(s_0/4,4\delta,4\delta)$$
is defined by
\begin{equation}\label{111904}
d_p(w,v)=4\delta||w-v||_{\mathcal{P}^p,D(s_0/4,4\delta,4\delta)}
\end{equation}
and
\begin{equation}\label{111905}d_{p}(w,\mathcal{T}):=\inf_{v\in
\mathcal{T}}d_p(w,v).
\end{equation}
\end{cor}

\begin{thm}\label{T3}(The long time stability of KAM tori) Given any $0\leq\mathcal{M}\leq
(2c\eta^8\epsilon)^{-1}$ (the same $c,\eta$ and $\epsilon$ as stated
in Theorem \ref{T1}), there exist a small positive $\delta_0$
depending on $s_0,r_0,n$ and $\mathcal{M}$, and a subset
$\Pi_{\tilde\eta}\subset\Pi_{\eta}$ satisfying
\begin{equation}Meas\ \Pi_{\tilde\eta} \geq(Meas\
\Pi_{\eta})(1-O(\tilde\eta)),\end{equation} where $\tilde \eta$ is some constant in $(0,1)$. For any $p\geq
8(\mathcal{M}+7)^{4}+1$, $0<\delta<\delta_0$ and for each
$\xi\in\Pi_{\tilde\eta}$, the KAM tori ${\mathcal T}$is stable in
long time, i.e. if $w(t)$ is a solution of Hamiltonian vector field
$X_{{H}}$ with the initial datum $w(0)=( w_x(0), w_y(0),
w_q(0),w_{\bar q}(0))$ satisfying
\begin{equation*}
d_{p}( w(0),{\mathcal{T}})\leq \delta,
\end{equation*}
then
\begin{equation*}
d_{p}(w(t),{\mathcal{T}})\leq 2\delta,
\qquad\mbox{for all}\ |t|\leq \delta^ {-\mathcal{M}}.
\end{equation*}
\end{thm}


\begin{rem}\label{R1}\label{R1}Theorem \ref{T1} is essentially due to Kuksin \cite{K3,K1}. However, in \cite{K3,K1}, the symplectic map
$$\Psi:D(s_0/2,0,0)\rightarrow D(s_0,r_0,r_0),$$
so the normal form $H\circ \Psi$ is degenerate. One can extend the definition domain $D(s_0/2,0,0)$ of  $\Psi$ to $ D(s_0/2,r_0/2,r_0/2)$ (even to the whole space) in view of a remak by P\"oschel in \cite{P2,P1} and an observation that  $\Psi$ is linear in $y$ and quadratic in $(q,\bar q)$. In many recent KAM theorems by, for example, Eliasson-Kuksin\cite{E-K}, Gr$\acute{\mbox{e}}$bert-Thomann\cite{GT} and Berti-Biasco \cite{BB1}, the extension is done in this line. In particular, the detail is given out in \cite{GT}. Unfortunately, up to now we do not know how to fulfill the tame property of the perturbed vector filed $X_{\breve{R}}$ in the extended domain in this line. The tame property of $X_{\breve R}$ is key ingredient in our present paper.  On the other hand, in the earlier work by Wayne\cite{W} there is another KAM iteration procedure which is a bit different from Kuksin's\cite{K3,K1}. In wayne's procedure, the definition domain of $\Psi$ is just $D(s_0/2,r_0/2,r_0/2)$, not necessary to extend it to a larger domain. In the present paper, we adopt Wayne's iteration procedure directly so that the tame property of $X_{\breve R}$ can be verified explicitly. The proof of Theorem \ref{T1} is well known. The aim of  providing the proof in  \S 5 is to verify the  tame property of $X_{\breve R}$. If the reader acknowledge the fact of the tame property of $X_{\breve R}$, the section \S 5 can be skipped.
\end{rem}


\begin{rem}{ Since the parameter set
 $\Pi\subset\mathbb{R}^{\mathbb{N}}$ is of infinite dimension,
 the measure in the above theorems should be in the sense of Kolmogorov.
 Actually it is enough to assume the parameter set is of finite dimension. Write
\begin{equation*}
\xi=(\xi^{n},\xi^{\mathcal{N}},\xi^{\mathbb{N}})\in
\mathbb{R}^n\times \mathbb{R}^{\mathcal{N}}\times \ell^2
\end{equation*}with a large
$\mathcal{N}$ will be given in Theorem {\ref{T3}}. In the proof of
constructing normal form of order 2 in Theorem \ref{T1}, it is
enough to regard $\xi^n$ as parameters. In order to get partial
normal form of order $\mathcal{M}+2$ around the KAM torus, it
suffices to take the $(\xi^n,\xi^{\mathcal{N}})$ as parameters in
Theorems \ref{T3} and \ref{T4}. Therefore, the measure can be understood  as
Lebesgue measure in this paper. }
\end{rem}
\begin{rem}Instead of equation (\ref{26}), we also can prove the
existence and long time stability of KAM tori for general nonlinear
Schr$\ddot{\mbox{o}}$dinger equations, such as
\begin{equation}\label{082003}
{\bf i}u_{t}=u_{xx}-V(x)u-\frac{\partial g(x,u,\bar u)}{\partial
\bar{u}},\qquad x\in\mathbb{T},\ t\in\mathbb{R},
\end{equation}
where $V$ is a smooth and $2\pi$ periodic potential, and
$g(x,u_1,u_2)$ is a smooth function on the domain
$\mathbb{T}\times\mathcal{U}$, $\mathcal{U}$ being a neighborhood of
the origin in $\mathbb{C}\times\mathbb{C}$, $g$ has a zero of order
three at $(u_1,u_2)=(0,0)$ and that $g(x,u,\bar u)\in\mathbb{R}$.
 Equations (\ref{082003}) were discussed in \cite{BG} and shown that
 the origin is stable in long time by the infinite dimensional
 Birkhoff normal form theorem.
\end{rem}


\section{Some discussions and ideas of the proof}
We begin by discussing some basic observations in
Bambusi-Gr$\acute{\mbox{e}}$bert {\cite{BG}}.  Consider an
infinite dimensional Hamiltonian system
\begin{equation}\label{072501}
H(q,\bar q)=H_0(q,\bar q)+P(q,\bar q),\qquad q,\bar q\in\ell^2_{p},
\end{equation}
with symplectic structure $\sqrt{-1}dq\wedge d\bar q$, where
$$H_0(q,\bar q)=\sum_{j\geq 1}\Omega_jq_j\bar q_j$$ is the quadratic
part and $$P(q,\bar
q)=\sum_{\beta,\gamma\in\mathbb{N}^{\mathbb{N}},|\beta|+|\gamma|\geq
3}P^{\beta\gamma}q^{\beta}\bar q^{\gamma}$$ is a smooth function
having a zero of order at least three at the origin. Note that  the
Hamiltonian (\ref{072501}) is a normal form of order 2 around the
origin. As a dynamical consequence, a solution starting in the
$\delta$-neighborhood of the origin stays in the
$\delta$-neighborhood along the time $|t|\leq \delta^{-1}$. To show
the origin is stable in a longer time such as $|t|\leq
\delta^{-\mathcal{M}}$ for any $\mathcal{M}\geq 0$, a natural way is
to construct a normal form of order $\mathcal{M}+1$ around the
origin. To this end, split $P(q,\bar q)$ into two parts, which is
\begin{equation*}
P(q,\bar q)=P_1(q,\bar q)+P_2(q,\bar q),
\end{equation*}
where
\begin{equation*}
P_1(q,\bar
q)=\sum_{\beta,\gamma\in\mathbb{N}^{\mathbb{N}},|\beta|+|\gamma|\leq
\mathcal{M}+1}P^{\beta\gamma}q^{\beta}\bar{q}^{\gamma}
\end{equation*}
is the part of low order,
and
\begin{equation*}
P_2(q,\bar
q)=\sum_{\beta,\gamma\in\mathbb{N}^{\mathbb{N}},|\beta|+|\gamma|\geq
\mathcal{M}+2}P^{\beta\gamma}q^{\beta}\bar{q}^{\gamma}
\end{equation*}is the part of high order. In order to remove all non-normalized terms
\begin{equation*}
\sum_{\beta,\gamma\in\mathbb{N}^{\mathbb{N}},|\beta|+|\gamma|\leq
\mathcal{M}+1,|\beta-\gamma|\neq0}P^{\beta\gamma}q^{\beta}\bar{q}^{\gamma}
\end{equation*}in $P_1(q,\bar
q)$, the following non-resonant
conditions:
\begin{equation}\label{081201}
|\langle \beta-\gamma,\Omega\rangle|\geq C(\mathcal{M})
\end{equation}  are needed, where
 $\beta,\gamma\in\mathbb{N}^{\mathbb{N}}$ satisfying $$
|\beta|+|\gamma|\leq\mathcal{M}+1,\qquad |\beta-\gamma|\neq 0,$$
$\Omega=(\Omega_1,\Omega_2,\dots)$ and $C(\mathcal{M})>0$ is a
constant depending on $\mathcal{M}$. However, the conditions
(\ref{081201}) are hardly to hold for infinite dimensional
Hamiltonian systems, since there are too many inequalities in
(\ref{081201}). A key idea in {\cite{BG}} is that a large part of
the nonlinearity is `not relevant' according to tame property and
all the remaining nonlinearity can be eliminated using a suitable
non-resonant condition. More precisely, spilt the variable
$q=(q_1,q_2,\dots)$ into two parts with a given large $\mathcal{N}$,
i.e. let $q=(\tilde q,\hat q)$, where $\tilde
q=(q_1,\dots,q_{\mathcal{N}})$ is called the low frequency variable
and $\hat q=(q_{\mathcal{N}+1},q_{\mathcal{N}+2},\dots)$ is called
the high frequency variable. Rewrite $P_1(q,\bar q)$ as
\begin{equation*}
P_1(q,\bar q)=P_{11}(q,\bar q)+P_{12}(q,\bar q),
\end{equation*}
where
\begin{equation*}
P_{11}(q,\bar
q)=\sum_{\tilde{\beta},\tilde\gamma\in\mathbb{N}^{\mathcal{N}},
\hat{\beta},\hat{\gamma}\in\mathbb{N}^{\mathbb{N}}}
P^{\beta\gamma}\tilde
q^{\tilde\beta}\hat{q}^{\hat\beta}\bar{\tilde q}^{\tilde\gamma}\bar{\hat
q}^{\hat\gamma},\qquad |\tilde\beta|+|\hat\beta|+|\tilde\gamma|+|\hat\gamma|\leq\mathcal{M}+1,\quad
|\hat\beta|+|\hat\gamma|\leq 2,
\end{equation*}
and
\begin{equation*}
P_{12}(q,\bar
q)=\sum_{\tilde{\beta},\tilde\gamma\in\mathbb{N}^{\mathcal{N}},
\hat{\beta},\hat{\gamma}\in\mathbb{N}^{\mathbb{N}}}
P^{\beta\gamma}\tilde
q^{\tilde\beta}\hat{q}^{\hat\beta}\bar{\tilde q}^{\tilde\gamma}\bar{\hat
q}^{\hat\gamma},\qquad |\tilde\beta|+|\hat\beta|+|\tilde\gamma|+|\hat\gamma|\leq\mathcal{M}+1,\quad
|\hat\beta|+|\hat\gamma|\geq 3,
\end{equation*}
with $\beta=(\tilde\beta,\hat\beta)$, $\gamma=(\tilde\gamma,\hat\gamma)$, $\tilde\beta=(\beta_1,\dots,\beta_{\mathcal{N}})$, $\hat\beta=(\beta_{\mathcal{N}+1},\beta_{\mathcal{N}+2},\dots)$,
$\tilde\gamma=(\gamma_1,\dots,\gamma_{\mathcal{N}})$ and $\hat\gamma=(\gamma_{\mathcal{N}+1},\gamma_{\mathcal{N}+2},\dots)$.
As stated in \cite{BG}, $P_{12}(q,\bar q)$ is the `non-relevant'
part. Moreover, the non-normalized terms in $P_{11}(q,\bar q)$ can be
removed by the non-resonant conditions
\begin{equation}\label{081202}
|\langle \tilde\beta-\tilde\beta,\tilde\Omega\rangle+\langle\hat\beta-\hat\gamma,\hat\Omega
\rangle|\geq C(\mathcal{M},\mathcal{N})
\end{equation}with
\begin{equation*}
|\tilde\beta|+|\hat\beta|+|\tilde\gamma|+|\hat\gamma|\leq\mathcal{M}+1, |\hat\beta|+|\hat\gamma|\leq2,
|\tilde\beta-\tilde\gamma|+|\hat\beta-\hat\gamma|\neq 0,
\end{equation*}
where $\tilde\Omega=(\Omega_1,\dots,\Omega_{\mathcal{N}})$,
$\hat{\Omega}=(\Omega_{\mathcal{N}+1},\Omega_{\mathcal{N}+2},\dots)$
and $C(\mathcal{M},\mathcal{N})$ is a positive constant depending on
$\mathcal{M}$ and $\mathcal{N}$. Note that there are less
inequalities in condition (\ref{081202}) than in conditions
(\ref{081201}), because of $|\hat\beta|+|\hat\gamma|\le 2$. More importantly, the non-resonant conditions
(\ref{081202}) are satisfied for many infinite dimensional
Hamiltonian systems. As a result in \cite{BG}, it is shown that if the
frequency $\Omega=(\tilde\Omega,\hat\Omega)$ satisfies the non-resonant conditions
(\ref{081202}), then there exists a symplectic transform $\Phi$ such
that
\begin{equation}\label{080107}
H\circ \Phi=H_0+Z+Q_1+Q_2,
\end{equation}
where $Z$ depends on the actions $I_j=|q_j|^2$,
$Q_1=O(||q||_p^{\mathcal{M}+2})$ and $Q_2=O(||\hat q||_p^{3})$ for
any $\mathcal{M}\geq 0$. Based on the partial normal form
(\ref{080107}) and tame property, the following dynamical result in
\cite{BG} is obtained: any solution with data in the
$\delta$-neighborhood of origin still stays in the
$\delta$-neighborhood of origin for time $|t|\leq
\delta^{-\mathcal{M}+1}$. Furthermore, the method in \cite{BG} can
be used to construct almost-global existence solutions for many
PDEs. For example, see Bambusi-Delort-Gr$\acute{\mbox{e}}$bert-Szeftel \cite{BDGS}.

In the present paper, we will prove long time stability of KAM
tori for infinite dimensional Hamiltonian systems. Note that the standard KAM
techniques constructs a non-degenerate normal form of order 2 (see Corollary  \ref{L1}). As a
consequence, it is easy to show the existence, linear stability and time $\delta^{-1}$
stability of KAM tori. In order to obtain a longer time stability such as
$\delta^{-\mathcal{M}}$ for any $\mathcal{M}\geq 0$, it is natural
to construct a normal form of order $\mathcal{M}+1$ around the KAM
tori. To this end, we split $\breve {R}(x,y,q,\bar q;\xi)$ (see (\ref{082902})) into two
parts, which is
\begin{equation*}
\breve {R}(x,y,q,\bar q;\xi)=\breve {R}_1(x,y,q,\bar q;\xi)+\breve
{R}_2(x,y,q,\bar q;\xi),
\end{equation*}
where
\begin{equation*}
\breve {R}_1(x,y,q,\bar
q;\xi)=\sum_{\alpha\in\mathbb{N}^n,\beta,\gamma\in\mathbb{N}^{\mathbb{N}},\atop 2|\alpha|+|\beta|+|\gamma|\leq
\mathcal{M}+1}\breve{R}^{\alpha\beta\gamma}(x;\xi)y^{\alpha}q^{\beta}\bar
q^{\gamma}
\end{equation*}
is the part of low order, and
\begin{equation*}
\breve {R}_2(x,y,q,\bar
q;\xi)=\sum_{\alpha\in\mathbb{N}^n,\beta,\gamma\in\mathbb{N}^{\mathbb{N}},\atop
2|\alpha|+|\beta|+|\gamma|\geq
\mathcal{M}+2}\breve{R}^{\alpha\beta\gamma}(x;\xi)y^{\alpha}q^{\beta}\bar
q^{\gamma}
\end{equation*}is the part of high order.
Expand $\breve{R}^{\alpha\beta\gamma}(x;\xi)$ into Fourier series $$
\breve{R}^{\alpha\beta\gamma}(x;\xi)=\sum_{k\in\mathbb{Z}^n}\widehat{\breve{R}^{\alpha\beta\gamma}}(k;\xi)e^{\sqrt{-1}\langle
k,x\rangle} ,$$ where $\widehat{\breve{R}^{\alpha\beta\gamma}}(k;\xi)$ is
the $k$-th Fourier coefficient of $\breve{R}^{\alpha\beta\gamma}(x;\xi)$.
Then we must remove all non-normalized terms in $\breve
{R}_1(x,y,q,\bar q;\xi)$, which are
\begin{equation*}
\sum_{\alpha\in\mathbb{N}^n,\beta,\gamma\in\mathbb{N}^{\mathbb{N}},\atop2|\alpha|+|\beta|+|\gamma|\leq \mathcal{M}+1, |k|+|\beta-\gamma|\neq
0}
\sum_{k\in\mathbb{Z}^n}\widehat{\breve{R}^{\alpha\beta\gamma}}(k;\xi)e^{\sqrt{-1}\langle
k,x\rangle}y^{\alpha}q^{\beta}\bar q^{\gamma},
\end{equation*}
where the following non-resonant conditions are needed
\begin{equation}\label{072503}
|\langle k,\omega\rangle+\langle \beta-\gamma,\Omega\rangle|\geq
C(k,\mathcal{M})
\end{equation}for any
$k\in\mathbb{Z}^n$ and $\beta,\gamma\in\mathbb{N}^{\mathbb{N}}$
satisfying
$$|\beta|+|\gamma|\leq \mathcal{M}+1,\qquad
|k|+|\beta-\gamma|\neq 0,$$where $C(k,\mathcal{M})$ is a positive
constant depending on $k$ and $\mathcal{M}$. Note that there are
more inequalities in conditions (\ref{072503}) than in conditions
(\ref{081201}). Therefore, the non-resonant conditions
(\ref{072503}) are more hardly to hold.

Following the idea in \cite{BG}, we would like to construct a
partial normal form around the KAM tori instead of a normal form.
But we have to face the following problems: (1) are the
`weakened' non-resonant conditions satisfied when constructing a
partial normal form of high order in the neighborhood of the KAM
tori? (2) how should we define tame property in the case that the
tangent direction exists ($n> 0$)? (3) does the tame property
preserve under the KAM iteration (infinite number of symplectic
transformations) and normal form iteration (finite number of symplectic
transformations)? Since  the number of the transformations is infinite and the transformations involve the action-angle variable  in the KAM iteration, the problem (3) is very hard to solve.

We will solve the above problems as follows. Firstly, following \cite{BG} we split
the normal
 variable $q$ into two parts with a given large $\mathcal{N}$, i.e. let $q=(\tilde q,\hat
 q)$, where $\tilde q=(q_1,\dots,q_{\mathcal{N}})$ is the low frequency
 normal variable and
 $\hat{q}=(q_{\mathcal{N}+1},q_{\mathcal{N}+2},\dots)$ is the high
 frequency normal variable. Rewrite $\breve R_1(x,y,q,\bar q;\xi)$ as
 \begin{equation*}
\breve R_1(x,y,q,\bar q;\xi)=\breve R_{11}(x,y,q,\bar q;\xi)+\breve
R_{12}(x,y,q,\bar q;\xi),
 \end{equation*}
 where\begin{equation*}
 \breve R_{11}(x,y,q,\bar q;\xi)=\sum_{\alpha\in\mathbb{N}^n,\tilde\beta,
 \tilde\gamma\in\mathbb{N}^{\mathcal{N}},\hat\beta,\hat\gamma\in\mathbb{N}^{\mathbb{N}},
 \atop 2|\alpha|+|\tilde\beta|+|\hat\beta|+|\tilde\gamma|+|\hat\gamma|\leq\mathcal{M}+1,
 |\hat\beta|+|\hat\gamma|\leq 2}{\breve R^{\alpha\beta\gamma}}(x;\xi)y^{\alpha}\tilde
 q^{\tilde\beta}{\hat q}^{\hat\beta}\bar{\tilde{q}}^{\tilde\gamma}\bar{\hat
 q}^{\hat\gamma}
 \end{equation*}
 and
 \begin{equation*}
 \breve R_{12}(x,y,q,\bar q;\xi)=\sum_{\alpha\in\mathbb{N}^n,\tilde\beta,\tilde\gamma\in\mathbb{N}^{\mathcal{N}}
 ,\hat\beta,\hat\gamma\in\mathbb{N}^{\mathbb{N}},
 \atop 2|\alpha|+|\tilde\beta|+|\hat\beta|+|\tilde\gamma|+|\hat\gamma|\leq\mathcal{M}+1,
 |\hat\beta|+|\hat\gamma|\geq3}{\breve R^{\alpha\beta\gamma}}(x;\xi)y^{\alpha}\tilde
 q^{\tilde\beta}{\hat q}^{\hat\beta}\bar{\tilde{q}}^{\tilde\gamma}\bar{\hat
 q}^{\hat\gamma}.
 \end{equation*}
Then we can obtain a partial normal form
 of order $\mathcal{M}+1$ by removing the non-normalized terms in $\breve R_{11}(x,y,q,\bar
 q;\xi)$ which are
 \begin{equation*}
 \sum_{\alpha\in\mathbb{N}^n,\tilde\beta,
 \tilde\gamma\in\mathbb{N}^{\mathcal{N}},\hat\beta,\hat\gamma\in\mathbb{N}^{\mathbb{N}},
 \atop 2|\alpha|+|\tilde\beta|+|\hat\beta|+|\tilde\gamma|+|\hat\gamma|\leq\mathcal{M}+1,
 |\hat\beta|+|\hat\gamma|\leq 2}\sum_{k\in\mathbb{Z}^n}\widehat{\breve R^{\alpha\beta\gamma}}(k;\xi)e^{\sqrt{-1}\langle k,x \rangle}y^{\alpha}\tilde
 q^{\tilde\beta}{\hat q}^{\hat\beta}\bar{\tilde{q}}^{\tilde\gamma}\bar{\hat
 q}^{\hat\gamma}
 \end{equation*}with
 \begin{equation*}|k|+|\tilde\beta-\tilde\gamma|+|\hat\beta-\hat\gamma|\neq 0
 \end{equation*}under the following non-resonant conditions
 \begin{equation}\label{072505}
 |\langle k,\omega\rangle+\langle
\tilde\beta-\tilde\gamma,\tilde\Omega\rangle+\langle\hat\beta-\hat\gamma,\hat\Omega\rangle|\geq
\frac{\tilde\eta}{C(\mathcal{M},\mathcal{N})(|k|+1)^{\tau}},
 \end{equation}
 where $k\in\mathbb{Z}^n$ and
 $$|\tilde\beta|+|\tilde\gamma|+|\hat\beta|+|\hat\gamma|\leq \mathcal{M}+1,|\hat\beta|+|\hat\gamma|\leq
 2,|k|+|\tilde\beta-\tilde\gamma|+|\hat\gamma-\hat\gamma|\neq 0,$$
 and $C(\mathcal{M},\mathcal{N})>0$ is a constant depending on
 $\mathcal{M}$ and $\mathcal{N}$ will be given in (\ref{090302}).
 Note the conditions (\ref{072505}) are reduced to the non-resonant conditions (\ref{081201}) 
 when $n=0$,  and they are similar to the standard non-resonant conditions
 \begin{equation*}
 |\langle k,\omega\rangle+\langle  l,\Omega\rangle|\geq
\frac{\tilde\eta}{(|k|+1)^{\tau}},\qquad k\in\mathbb{Z}^n,l\in
\mathbb{Z}^{\mathbb{N}},|l|\leq2,|k|+|l|\neq 0
 \end{equation*}in KAM
 technique while without splitting the normal variable $q$.
In Section \ref{091406}, we show the non-resonant conditions
(\ref{072505}) are satisfied by removing
 the parameters of a small measure.

 Secondly, note an important fact that when the dimensional of the tori
is 0 ($n=0$), the tori can be considered as a point, which is just the case discussed in \cite{BG}. We define
$p$-tame norm ($p$-tame property)
 of Hamiltonian vector field $X_W$ by the following steps:

 {\it Step 1,} Write $W(x,y,z;\xi)=\sum_{\beta}W^{\beta}(x,y;\xi)z^{\beta}$ and take the norm $||\cdot||_{D(s,r)\times\Pi}$ of $W^{\beta}(x,y;\xi)$:
   \[\lfloor W\rceil_{D(s,r)\times\Pi}(z):=\sum_{\beta}||W^{\beta}(x,y;\xi)||_{D(s,r)\times\Pi} \, z^{\beta}.\]
   (See (\ref{081602}) in Definition \ref{021002} for the
 detail).

 {\it Step 2,} note that
 the modulus $\lfloor W\rceil_{D(s,r)\times\Pi}(z)$ depends on normal variables $z=(q,\bar q)$ only (independent of tangent variable $(x,y)$ and parameter $\xi$). Then we follow the method in \cite{BG} to define the $p$-tame norm of $W_z$ (the normal direction of Hamiltonian vector field $X_W$). More precisely, first consider a Hamiltonian $W(x,y,z;\xi)$, the modulus of which is a homogeneous polynomial about $z$ of degree $h$ (see (\ref{011601})). Then define $p$-tame operator
 norm $$|||\cdot|||_{p,D(s,r)\times\Pi}^T$$of $W_z$ (see (\ref{091402})
 in Definition \ref{071803}).

 {\it Step 3,} for a general Hamiltonian
 $$W(x,y,z;\xi)=\sum_{h\geq0}W_h(x,y,z;xi),$$
 where $\lfloor W_h\rceil_{D(s,r)\times\Pi}(z)$ is a homogeneous polynomial about $z$ of degree $h$. It is natural to define $p$-tame
 norm \begin{equation}\label{011701}
 |||W_z|||_{p,D(s,r,r)\times\Pi}^T=\sum_{h\geq1}|||(W_h)_z|||_{p,D(s,r)\times\Pi}^Tr^{h-1}.
 \end{equation}However, $p$-tame
 norm defined by (\ref{011701}) is not enough to show the persistence of $p$-tame norm under Poisson bracket. To this end, $\lfloor W_h\rceil_{D(s,r)\times\Pi}(z)$ is required as a bounded map from
$\ell^2_1\times\ell^2_1$ into
$\ell^2_1\times\ell^2_1$, which is also used in \cite{BG}. Then we can define the $p$-tame norm of $(W_h)_z$ by
\begin{equation*}
|||{(W_h)_z}|||_{p,D(s,r,r)\times\Pi}^T=\max
\left\{|||{(W_h)_z}|||_{p,D(s,r)\times\Pi}^T,|||{(W_h)_z}|||_{1,D(s,r)\times\Pi}^T\right\}r^{h-1}
\end{equation*}(see (\ref{091403})
 in Definition \ref{071803}).

{\it Step 4,} we deal with the tangent direction of Hamiltonian vector field $X_W$. Note that $W_x$ and $W_y$ are finite dimensional, so there is no so-called $p$-tame property. But to guarantee the persistence of
$p$-tame property under
Poisson bracket, define the operator norm
\begin{equation*}
|||\cdot|||_{D(s,r)\times\Pi}
\end{equation*}and
the norm
\begin{equation*}
|||\cdot|||_{D(s,r,r)\times\Pi}
\end{equation*}of $W_x$ and $W_y$ by (\ref{091404}) and (\ref{091405}) in
Definition \ref{071802} respectively, where $||W_x||_{D(s,r)\times\Pi}$ and $||W_y||_{D(s,r)\times\Pi}$ are required as bounded maps
from
$\ell^2_1\times\ell^2_1$ into
$\mathbb{C}^n$.

{\it Step 5,} give the
 definition of $p$-tame norm ($p$-tame property) of Hamiltonian vector field $X_W$ for a general Hamiltonian $W(x,y,z;\xi)$ in Definition {\ref{080204}}.

Thirdly, we should prove $p$-tame property
 survives under KAM iteration and normal form iteration, which is the key part
 in this paper (see details in Section \ref{092701}). The essential difference between this paper and \cite{BG}
 is that before constructing a partial normal form of $\mathcal{M}+1$, we
 need to use infinite symplectic transformations (KAM iteration) to
 obtain a normal form of order 2 (see Theorem \ref{T1}), while in \cite{BG} a normal form of order 2 is already there (see (\ref{072501})).
 Thus, we should prove $p$-tame property is preserved under infinite
 symplectic transformations. To this end, we need estimate the $p$-tame norm of Poisson bracket of two Hamiltonian (Theorem \ref{021102}), the composition of a Hamiltonian with a Hamiltonian flow (Theorem \ref{081914}) and the solution of homological equation (Theorem \ref{0005}).

  Also, we have to face
 frequency shift in KAM iteration and more complicated small divisors than in \cite{BG} because
 of the existence of tangent direction of Hamiltonian vector field.
Here, we point out that we can only obtain time $|t|\leq
\delta^{-\mathcal{M}+1}$ stability of the KAM tori for
$\mathcal{M}\leq (2c\eta^8\epsilon)^{-1}$, comparing to any
$\mathcal{M}\geq 0$ in \cite{BG} because of the problem of the
frequency shift.

Finally, basing on the partial normal form of high order and $p$-tame property, we get
 the dynamical
 consequence that solutions starting in the $\delta$-neighborhood of the KAM
torus still remain in the $\delta$-neighborhood of the KAM torus for
time $|t|\leq \delta^{-\mathcal{M}+1}$, i.e. most of KAM tori are
long time stable.

At the end of this section, we will give some simple estimate in following remarks.
\begin{rem}\label{011604}In view of (\ref{092704}) in Definition \ref{090501}, it is easy to verify that
\begin{equation}\label{011603}
||\cdot||_{D(s-\sigma)\times\Pi}\leq||\cdot||_{D(s)\times\Pi},
\end{equation}for $0<\sigma<s$. Moreover, in view of (\ref{092703}) in Definition \ref{060602}, the following inequalities hold
\begin{equation}\label{011702}
||\cdot||_{D(s-\sigma,r)\times\Pi}\leq||\cdot||_{D(s,r)\times\Pi},
\end{equation}
and
\begin{equation}\label{011703}
||\cdot||_{D(s,r-\sigma')\times\Pi}\leq||\cdot||_{D(s,r)\times\Pi},
\end{equation}
where $0<\sigma'<r$. Furthermore, (\ref{011702}) and (\ref{011703}) implies
\begin{equation}\label{011709}
||\cdot||_{D(s-\sigma,r-\sigma')\times\Pi}\leq||\cdot||_{D(s,r)\times\Pi}.
\end{equation}
\end{rem}
\begin{rem}Given a Hamiltonian $W(x,y,z;\xi)$ has $p$-tame property on the
domain $D(s,r,r)\times\Pi$. It is easy to verify that
\begin{equation}
|||X_W|||_{p,D(s-\sigma,r,r)\times\Pi}^{T}\leq |||X_W|||_{p,D(s,r,r)\times\Pi}^{T},
\end{equation}where $0<\sigma<r$. But the following inequality usually is false
\begin{equation*}
|||X_W|||_{p,D(s,r-\sigma',r-\sigma')\times\Pi}^{T}\leq |||X_W|||_{p,D(s,r,r)\times\Pi}^{T},
\end{equation*}where $0<\sigma'<r$.
However, if let $0<\sigma'<r/2$, we have the following estimate
\begin{equation}\label{011708}
|||X_W|||_{p,D(s,r-\sigma',r-\sigma')\times\Pi}^{T}\leq 4|||X_W|||_{p,D(s,r,r)\times\Pi}^{T},
\end{equation}
since $0<\sigma'<r/2$ implies $r/2<r-\sigma'<r$.
\end{rem}
\begin{rem}
Based on (\ref{091402})
 in Definition \ref{071803}, for each
$(x,y,z)\in\mathcal{P}^p$ and $\xi\in\Pi$, the following estimate holds
\begin{equation}\label{011602}
||(W_h)_z(x,y,z;\xi)||_p\leq
|||{(W_h)_z}|||_{p,D(s,r)\times\Pi}^{T}||z||_p||z||_1^{\max\{h-2,0\}}.
\end{equation}
\end{rem}
\begin{rem}Based on (\ref{091404}) in
Definition \ref{071802}, for each $(x,y,z)\in\mathcal{P}^p$ and $\xi\in\Pi$, the following estimates hold
\begin{equation} ||(W_h)_x(x,y,z;\xi)||\leq
|||{(W_h)_x}|||_{D(s,r)\times\Pi}||z||_1^{h},
\end{equation}
and
\begin{equation} ||(W_h)_y(x,y,z;\xi)||\leq
|||{(W_h)_y}|||_{D(s,r)\times\Pi}||z||_1^{h}.
\end{equation}
\end{rem}
\begin{rem}\label{080803}Note that
 $$||\cdot||_{D(s,r)\times\Pi}\geq0.$$ Then in view of
Definition \ref{071803} and Definition \ref{071802}, it is easy to
verify that
\begin{eqnarray*}
|||{W_z}|||_{p,D(s,r)\times\Pi}^{T}& =&\sup_{0\neq z^{(j)}\in \ell
^2_p\times\ell^2_p,1\leq j\leq h-1,z^{(j)}\geq 0}
\frac{||\widetilde{\lfloor{W_z}\rceil}_{D(s,r)\times\Pi}(z^{(1)},\dots,z^{(h-1)})||_
{p}} {||(z^{h-1})||_{p,1}},\end{eqnarray*}and
\begin{eqnarray*}|||W_{u}|||_{D(s,r)\times\Pi}& =&\sup_{0\neq z^{(j)}\in \ell
^2_1\times\ell^2_1,1\leq j\leq h,z^{(j)}\geq0}
\frac{||\widetilde{\lfloor{W_{u}}\rceil}_{D(s,r)\times\Pi}(z^{(1)},\dots,z^{(h)})||}
{||(z^{h})||_{1,1}},
\end{eqnarray*}where $u=x$ or $y$, $z^{(j)}\geq 0$
means $z_i^{(j)}\geq0$ for $i\in\bar{\mathbb{Z}}$ and
$z^{(j)}=(z_i^{(j)})_{i\in\bar{\mathbb{Z}}}$. Without loss of
generality, we always assume that each entry of $z^{(j)}$ is
non-negative, when estimating $p$-tame norm of Hamiltonian vector
field in the rest of this paper. Also see the same discussion in
Remark 4.5 in \cite{BG}.
\end{rem}

\section{Properties of the Hamiltonian with $p$-tame
property}\label{092701}
The following theorem will show that $p$-tame property persists under Poisson bracket.
\begin{thm}\label{021102}
Suppose that both Hamiltonian functions
$$U(x,y,z;\xi)=\sum_{\beta\in\mathbb{N}^{\bar{\mathbb{Z}}}}U^{\beta}(x,y;\xi)z^{\beta}$$
and
$$V(x,y,z;\xi)=\sum_{\beta\in\mathbb{N}^{\bar{\mathbb{Z}}}}V^{\beta}(x,y;\xi)z^{\beta},$$
satisfy $p$-tame property on the domain $D(s,r,r)\times\Pi$,
where
\begin{equation*}
U^{\beta}(x,y;\xi)=\sum_{\alpha\in\mathbb{N}^{{n}}}U^{\alpha\beta}(x;\xi)y^{\alpha}
,\end{equation*} and
\begin{equation*}
V^{\beta}(x,y;\xi)=\sum_{\alpha\in\mathbb{N}^{{n}}}V^{\alpha\beta}(x;\xi)y^{\alpha}
.\end{equation*} Then the Poisson bracket $\{U,V\}(x,y,z;\xi)$ of
$U(x,y,z;\xi)$ and $V(x,y,z;\xi)$ with respect to the symplectic structure $dy\wedge dx+\sqrt{-1}\sum_{j\geq 1}dz_{-j}\wedge dz_j$ has $p$-tame property on the
domain $D(s-\sigma,r-\sigma',r-\sigma')\times\Pi$ for
$0<\sigma<s,0<\sigma'<r/2$. Moreover, the following inequality holds
\begin{eqnarray}&&|||X_{\{U,V\}}|||_{p,D(s-\sigma,r-\sigma',r-\sigma')\times\Pi}^T
\nonumber\\
\label{091120}&\leq&C\max\left\{\frac 1{\sigma},\frac
{r}{\sigma'}\right\}|||X_{U}|||_{p,D(s,r,r)\times\Pi}^T
|||X_{V}|||_{p,D(s,r,r)\times\Pi}^T,
\end{eqnarray}
where $C>0$ is an absolute constant.
\end{thm}
\begin{proof}
By a direct calculation,
\begin{equation*}
X_{\{U,V\}}=(\{U,V\}_y,-\{U,V\}_x,\sqrt{-1}\{U,V\}_{\bar
q},-\sqrt{-1}\{U,V\}_{ q})=DX_U\cdot X_V-DX_V\cdot X_U,
\end{equation*}where
\begin{eqnarray*}\label{010}DX_U\cdot X_V&=&\left(\begin{array}{cccc}
U_{yx}&U_{yy}&U_{yq}&U_{y\bar{q}}\\-U_{xx}&-U_{xy}&-U_{xq}&-U_{x\bar
q}\\
\sqrt{-1}U_{\bar qx}&\sqrt{-1}U_{\bar qy}&\sqrt{-1}U_{\bar qq}&\sqrt{-1}U_{\bar q\bar{q}}\\
-\sqrt{-1}U_{ qx}&-\sqrt{-1}U_{ qy}&-\sqrt{-1}U_{ qq}&-\sqrt{-1}U_{ q\bar{q}}
\end{array}\right)\left(\begin{array}{cccc}
V_{y}\\-V_{x}\\ \sqrt{-1}V_{\bar q}\\
-\sqrt{-1}V_{ q}
\end{array}\right)
\end{eqnarray*}
and \begin{eqnarray*} DX_V\cdot X_U&=&\left(\begin{array}{cccc}
V_{yx}&V_{yy}&V_{yq}&V_{y\bar{q}}\\-V_{xx}&-V_{xy}&-V_{xq}&-V_{x\bar
q}\\\sqrt{-1}V_{\bar qx}&\sqrt{-1}V_{\bar qy}&\sqrt{-1}V_{\bar qq}&\sqrt{-1}V_{\bar q\bar{q}}\\
-\sqrt{-1}V_{ qx}&-\sqrt{-1}V_{qy}&-\sqrt{-1}V_{ qq}&-\sqrt{-1}V_{ q\bar{q}}
\end{array}\right)\left(\begin{array}{cccc}
U_{y}\\-U_{x}\\
\sqrt{-1}U_{\bar q}\\
-\sqrt{-1}U_{q}
\end{array}\right).
\end{eqnarray*}
Thus, there are 32 terms in $X_{\{U,V\}}$, and we classify the 32
terms into 4 cases (for simplicity, we will omit the coefficients
$\pm\sqrt{-1}$ or $-1$ sometimes, which do not affect the estimate
of
$p$-tame norm below):\\
$\bf{Case}\ 1: finite -  finite$ . $$\sum_{j=1
}^nU_{xx_j}V_{y_j},\quad \sum_{j=1
}^nU_{xy_j}V_{x_j},\quad
\sum_{j=1
}^nU_{yx_j}V_{y_j},\quad \sum_{j=1
}^nU_{yy_j}V_{x_j},$$and
$$\sum_{j=1
}^nV_{xx_j}U_{y_j},\quad \sum_{j=1
}^nV_{xy_j}U_{x_j},
\quad \sum_{j=1
}^nV_{yx_j}U_{y_j},\quad \sum_{j=1
}^nV_{yy_j}U_{x_j};$$ $\bf{Case}\ 2: finite-infinite$. $$\sum_{j\geq
1}U_{xq_j} V_{\bar q_j},\quad \sum_{j\geq
1}U_{x\bar{q}_j}V_{q_j},\quad \sum_{j\geq 1}U_{yq_j}V_{\bar
q_j},\quad \sum_{j\geq 1}U_{y\bar q_j}V_{q_j},$$and
$$\sum_{j\geq 1}V_{xq_j}U_{\bar q_j},\quad \sum_{j\geq 1}V_{x\bar q_j}U_{q_j},\quad \sum_{j\geq 1}
V_{yq_j}U_{\bar q_j},\quad \sum_{j\geq 1}V_{y\bar q_j}U_{q_j};$$
$\bf{Case}\ 3: infinite-finite$.
$$\sum_{j=1
}^nU_{qx_j}V_{y_j},\quad \sum_{j=1
}^nU_{qy_j}V_{x_j},
\quad \sum_{j=1
}^nU_{\bar{q}x_j}V_{y_j},\quad \sum_{j=1
}^nU_{\bar{q}y_j}V_{x_j},$$and $$\sum_{j=1
}^nV_{qx_j}U_{y_j},\quad \sum_{j=1
}^nV_{qy_j}U_{x_j},\quad
\sum_{j=1
}^nV_{\bar{q}x_j}U_{y_j},\quad \sum_{j=1
}^nV_{\bar{q}y_j}U_{x_j};$$ $\bf{Case}\ 4:infinite-infinite$.
$$\sum_{j\geq 1}U_{qq_j}V_{\bar q_j},\quad \sum_{j\geq 1}U_{q\bar q_j}V_{q_j},\quad
\sum_{j\geq 1}U_{\bar{q}q_j}V_{\bar q_j},\quad \sum_{j\geq
1}U_{\bar{q}\bar{q}_j}V_{q_j},$$and $$\sum_{j\geq 1}V_{qq_j}U_{\bar
q_j},\quad \sum_{j\geq 1}V_{q\bar q_j}U_{q_j},\quad \sum_{j\geq
1}V_{\bar qq_j}U_{\bar q_j},\quad \sum_{j\geq 1}V_{\bar q\bar
q_j}U_{q_j}.$$We will drop the index $\Pi$ for shorten notations and
regard $x$ and $y$ as scalar for simplicity, when estimating
$p$-tame norm of Hamiltonian vector field $X_{\{U,V\}}$ below.

Suppose $U(x,y,z;\xi)$ and $V(x,y,z;\xi)$ are homogeneous polynomials about $z$, that is
\begin{equation}\label{102301}U(x,y,z;\xi):=U_h(x,y,z;\xi)=\sum_{\beta\in\mathbb{N}^{\bar{\mathbb{Z}}},
|\beta|=h}U^{\beta}(x,y;\xi)z^{\beta} \end{equation} and
\begin{equation}\label{102302}V(x,y,z;\xi):=V_l(x,y,z;\xi)=\sum_{\beta\in\mathbb{N}^{\bar{\mathbb{Z}}},|\beta|=l}
V^{\beta}(x,y;\xi)z^{\beta},\end{equation}
for some $h,l\in\mathbb{N}.$

\textbf{Step 1. Estimate
$|||\{U,V\}_x|||_{D(s-\sigma,r-\sigma',r-\sigma')}$ and $|||\{U,V\}_y|||_{D(s-\sigma,r-\sigma',r-\sigma')}$.}

In this step, we will give the following estimates
\begin{equation}\label{042806}
\frac1{(r-\sigma')^2}|||\{U,V\}_x|||_{D(s-\sigma,r-\sigma',r-\sigma')}\leq
\frac {32}{e\sigma}|||
X_U|||_{p,D(s,r,r)}^T|||X_V|||_{p,D(s,r,r)}^{T},
\end{equation}
and
\begin{equation}\label{042807}
|||\{U,V\}_y|||_{D(s-\sigma,r-\sigma',r-\sigma')}\leq \frac
{8r}{\sigma'}||| X_U|||_{p,D(s,r,r)}^{T}|||X_V|||_{p,D(s,r,r)}^{T}.
\end{equation}

Note that
\begin{eqnarray*}
\{U,V\}_x&=&U_{xx}V_y-U_{xy}V_x-V_{xx}
U_y+V_{xy}U_x\quad (\mbox{case 1})\\
&&+\sqrt{-1}\sum_{j\geq 1}(U_{xq_j}V_{\bar{q}_j}-U_{x\bar
q_j}V_{q_j}-V_{xq_j}U_{\bar{q}_j}+V_{x\bar q_j}U_{q_j})\quad
(\mbox{case 2}),
\end{eqnarray*}
and
\begin{eqnarray*}
\{U,V\}_y&=&U_{yx}V_y-U_{yy}V_x-V_{yx}
U_y+V_{yy}U_x\quad (\mbox{case 1})\\
&&+\sqrt{-1}\sum_{j\geq 1}(U_{yq_j}V_{\bar{q}_j}-U_{y\bar
q_j}V_{q_j}-V_{yq_j}U_{\bar{q}_j}+V_{y\bar q_j}U_{q_j})\quad
(\mbox{case 2}).
\end{eqnarray*}
Without loss of generality, we just consider the term $U_{xx}V_y$, which is in case $1$, and the term $\sum_{j\geq 1}(U_{yq_j}V_{\bar{q}_j}-U_{y\bar
q_j}V_{q_j})$, which is in case 2, and it is sufficient to give the following estimates $$\frac{1}{(r-\sigma')^2}|||U_{xx}V_y|||_{D(s-\sigma,r-\sigma',r-\sigma')}\leq
\frac {4}{e\sigma}|||
X_U|||_{p,D(s,r,r)}^T|||X_V|||_{p,D(s,r,r)}^{T},$$
and
$$|||\sum_{j\geq 1}(U_{yq_j}V_{\bar{q}_j}-U_{y\bar
q_j}V_{q_j})|||_{D(s-\sigma,r-\sigma',r-\sigma')}\leq
\frac {r}{\sigma'}|||
X_U|||_{p,D(s,r,r)}^T|||X_V|||_{p,D(s,r,r)}^{T}$$
respectively.

Let $m=h+l$ and $\tau_m$ be an $m$-permutation. Since $U(x,y,z;\xi)$ and $V(x,y,z;\xi)$ have $p$-tame property on
the domain $D(s,r,r)\times\Pi$, then by (\ref{091404}) in
Definition \ref{071802}, we get
\begin{equation}\label{080302'}
\widetilde{\lfloor
U_{x}\rceil}_{D(s,r)}(z^{(\tau_{m}(1))},\dots,z^{(\tau_{m}(h))})\leq|||
U_{x}|||_{D(s,r)}||z^{(\tau_{m}(1))}||_1\cdots||z^{(\tau_{m}(h))}||_{1},
\end{equation}
and
\begin{equation}\label{080303}
\widetilde{\lfloor V_{y}\rceil}_{D(s,r)}(z^{(\tau_m(h+1))},\dots,
z^{(\tau_{m}(m))})\leq|||
V_y|||_{D(s,r)}||z^{(\tau_{m}(h+1))}||_1\cdots||z^{(\tau_{m}(m))}||_{1}.
\end{equation}
By the generalized Cauchy estimate in Lemma \ref{050402},
\begin{eqnarray*}\label{090502}
||U_{xx}||_{D(s-\sigma,r)}\leq \frac
{1}{e\sigma}||U_x||_{D(s,r)}.
\end{eqnarray*}
Hence, in view of (\ref{081602}) in Definition 3, the
definition of symmetric linear form (see (\ref{080802})) and Remark \ref{080803}, we have
\begin{eqnarray}\label{080301'}
\widetilde{\lfloor
U_{xx}\rceil}_{D(s-\sigma,r)}(z^{(\tau_m(1))},\dots,z^{(\tau_m(h))})\leq
\frac {1}{e\sigma}\widetilde{\lfloor
U_x\rceil}_{D(s,r)}(z^{(\tau_m(1))},\dots,z^{(\tau_m(h))}).
\end{eqnarray}Based on the inequalities (\ref{080302'})-(\ref{080301'}),
\begin{eqnarray}
&&\nonumber\widetilde{\lfloor
U_{xx}\rceil}_{D(s-\sigma,r)}(z^{(\tau_{m}(1))},\dots,z^{(\tau_{m}(h))})
\widetilde{\lfloor
V_{y}\rceil}_{D(s-\sigma,r)}(z^{(\tau_m(h+1))},\dots,
z^{(\tau_{m}(m))})\\
&\leq&\frac {1}{e\sigma}||| U_{x}|||_{D(s,r)}|||
V_y|||_{D(s,r)}||z^{(\tau_{m}(1))}||_1\cdots||z^{(\tau_{m}(m))}||_{1}\label{080306}.
\end{eqnarray}Then we obtain
\begin{eqnarray}\nonumber&&\widetilde{\lfloor
U_{xx}V_y\rceil}_{D(s-\sigma,r-\sigma')}(z^{(1)}, \dots,z^{(m)})
\\
\nonumber&\leq&\widetilde{\lfloor
U_{xx}V_y\rceil}_{D(s-\sigma,r)}(z^{(1)}, \dots,z^{(m)})\qquad \mbox{(based on Remark \ref{011604} and Remark \ref{080803})}
\\
&\leq&\nonumber\label{080608}\widetilde{\lfloor
U_{xx}\rceil_{D(s-\sigma,r)}\lfloor
V_y\rceil_{D(s-\sigma,r)}}(z^{(1)},\dots,z^{(m)})\qquad \mbox{(based on Lemma \ref{051701})}\\
 &=&\frac1{m!}\sum_{\tau_{m}}\widetilde{\lfloor
U_{xx}\rceil}_{D(s-\sigma,r)}(z^{(\tau_{m}(1))},\dots,z^{(\tau_{m}(h))})
\widetilde{\lfloor
V_{y}\rceil}_{D(s-\sigma,r)}(z^{(\tau_m(h+1))},\dots,
z^{(\tau_{m}(m))}),\nonumber\\
&&\mbox{(based on Lemma \ref{042201})}\label{080305}\\
&\leq&\frac {1}{e\sigma}||| U_{x}|||_{D(s,r)}|||
V_y|||_{D(s,r)}||(z^{m})||_{1,1},\label{080309}
\end{eqnarray}where the last inequality is based on (\ref{080306}) and
the fact that
\begin{eqnarray}\label{080307}
\frac{1}{m!}\sum_{\tau_m}||z^{(\tau_{m}(1))}||_1\cdots||z^{(\tau_{m}(m))}||_{1}
=||z^{(1)}||_{1}\cdots ||z^{(m)}||_{1}=||(z^{m})||_{1,1}.
\end{eqnarray}
By (\ref{091404}) in Definition \ref{071802} and the inequality
(\ref{080309}), it is easy to see that
\begin{equation}\label{081605}
|||{ U_{xx} V_y}|||_{D(s-\sigma,r-\sigma')}\leq \frac {1}{e\sigma}|||
U_{x}|||_{D(s,r)}||| V_y|||_{D(s,r)}.
\end{equation}
Finally, we obtain
\begin{eqnarray}
&&\nonumber\frac1{(r-\sigma')^2}|||{
U_{xx}V_y}|||_{D(s-\sigma,r-\sigma',r-\sigma')}\\
 \nonumber &=&\frac1{(r-\sigma')^2}|||{
U_{xx}V_y}|||_{D(s-\sigma,r-\sigma')}(r-\sigma')^{h+l}\nonumber\qquad
\mbox{(by (\ref{091405}) in Definition \ref{071802})}\\
\nonumber&\leq&
 \frac {1}{e\sigma}\cdot\frac {r^2}{(r-\sigma')^2}\cdot\frac{1}{r^2}||| U_{x}|||_{D(s,r)}|||
V_y|||_{D(s,r)}(r-\sigma')^{h+l}\\&&\nonumber \mbox{(based on the
inequality (\ref{081605}))}\\
\nonumber&\leq&
 \frac {1}{e\sigma}\cdot\frac {r^2}{(r-\sigma')^2}\cdot\frac{1}{r^2}||| U_{x}|||_{D(s,r)}|||
V_y|||_{D(s,r)}r^{h+l}\\
&\leq&\frac {4}{e\sigma}\left(\frac1{r^2}||| U_{x}|||_{D(s,r)}r^h\right)\left(|||
V_y|||_{D(s,r)}r^{l}\right)\nonumber\\
 &&\mbox{(based on $0<\sigma'<r/2$ implying $r/2<r-\sigma'<r$)}\nonumber\\
&=&\frac {4}{e\sigma}\left(\frac1{r^2}|||
U_{x}|||_{D(s,r,r)}\right)|||V_y|||_{D(s,r,r)}\qquad \mbox{(by
(\ref{091405}) in Definition \ref{071802})}\nonumber\\
&\leq&\frac {4}{e\sigma}\left(|||
X_U|||^T_{p,D(s,r,r)}\right)|||X_V|||^T_{p,D(s,r,r)}.
\label{042803}
\end{eqnarray}
Denote by
$$U_{yz}\cdot
V_{z}=\sqrt{-1}\sum_{j\geq 1}(U_{yq_j}V_{\bar q_j}-U_{y\bar
q_j}V_{q_j}).$$
Let $j=m-2=h+l-2$ (here we assume $j\geq0$, otherwise $U_z=0$ or $V_z=0$). Then we obtain
\begin{eqnarray}
&&\nonumber\widetilde{ \lfloor U_{yz}\cdot
V_{z}}\rceil_{D(s-\sigma,r-\sigma')}(z^{(1)},\dots,z^{(j)})\\
&\leq&\nonumber\frac1{j!}\sum_{\tau_{j}}\widetilde{ \lfloor
U_{yz}\rceil}_{D(s-\sigma,r-\sigma')}(z^{\tau_{j}(1)},\dots,z^{\tau_{j}(h-1)})\cdot
\widetilde{ \lfloor
V_{z}\rceil}_{D(s-\sigma,r-\sigma')}(z^{\tau_{j}(h)},\dots,z^{\tau_{j}(j)})\\
&&\nonumber \mbox{(following the proof of the inequality
(\ref{080305}) and $\tau_j$ is a $j$-permutation)}
\\&\leq&\nonumber\frac1{j!}\sum_{\tau_{j}}||\widetilde{ \lfloor
U_{yz}\rceil}_{D(s-\sigma,r-\sigma')}(z^{\tau_{j}(1)},\dots,z^{\tau_{j}(h-1)})||_1
||\widetilde{ \lfloor
V_{z}\rceil}_{D(s-\sigma,r-\sigma')}(z^{\tau_{j}(h)},\dots,z^{\tau_{j}(j)})||_1\\
&&\nonumber\mbox{(based on the inequality $ |z\cdot \tilde
z|=|\sum_{j\in\bar{\mathbb{Z}}}z_j\tilde z_{j}|\leq||z||_0||\tilde
z||_0\leq ||z||_1||\tilde z||_1
$})\\
&\leq&
\nonumber\frac1{j!}\sum_{\tau_{j}}\frac{1}{r\sigma'}||\widetilde{
\lfloor
U_{z}\rceil}_{D(s-\sigma,r)}(z^{\tau_{j}(1)},\dots,z^{\tau_{j}(h-1)})||_1
||\widetilde{ \lfloor
V_{z}\rceil}_{D(s,r)}(z^{\tau_{j}(h)},\dots,z^{\tau_{j}(j)})||_1\\
&&\nonumber \mbox{(based on the generalized Cauchy estimate (\ref{071902}) in Lemma \ref{050403})}\\
&\leq&\label{080308}\frac{1}{r\sigma'}|||
U_{z}|||_{1,D(s,r)}^T||| V_{z}|||_{1,D(s,r)}^T||(z^{j})||_{1,1},
\end{eqnarray}
where the last inequality is based on the inequality (\ref{011702}) in Remark \ref{011604}, the formula (\ref{080307}) for $m=j$ and the inequalities
\begin{equation*} ||\widetilde{\lfloor
U_{z}\rceil}_{D(s,r)}(z^{(\tau_{j}(1))},\dots,z^{(\tau_{j}(h-1))})||_1\leq|||
U_{z}|||_{1,D(s,r)}^T||z^{(\tau_{j}(1))}||_1\cdots||z^{(\tau_{j}(h-1))}||_{1},
\end{equation*}
and
\begin{equation*}
||\widetilde{\lfloor V_{z}\rceil}_{D(s,r)}(z^{(\tau_j(h))},\dots,
z^{(\tau_{j}(j))})||_1\leq|||
V_z|||_{1,D(s,r)}^T||z^{(\tau_{j}(h))}||_1\cdots||z^{(\tau_{j}(j))}||_{1},
\end{equation*}since $U(x,y,z;\xi)$ and
$V(x,y,z;\xi)$ having $p$-tame property on the domain
$D(s,r,r)\times\Pi$ (see (\ref{091402}) in Definition \ref{071803}
for $p=1$).

According to the estimate (\ref{080308})
and in view of (\ref{091404}) in Definition \ref{071802}, we obtain
\begin{equation}\label{081606}
|||{ U_{yz}\cdot V_{z}}|||_{D(s-\sigma,r-\sigma')}\leq \frac
{1}{r\sigma'}||| U_{z}|||_{1,D(s,r)}^T||| V_{z}|||_{1,D(s,r)}^T.
\end{equation}Furthermore,
\begin{eqnarray}
\nonumber&&|||{ U_{yz}V_z}|||_{D(s-\sigma,r-\sigma',r-\sigma')}\\
\nonumber &=&|||{
U_{yz}V_z}|||_{D(s-\sigma,r-\sigma')}(r-\sigma')^{h+l-2}\nonumber\qquad
\mbox{(by (\ref{091405}) in Definition \ref{071802})}
\\ \nonumber&\leq&
 \frac {1}{r\sigma'}||| U_{z}|||^T_{1,D(s,r)}|||
V_z|||^T_{1,D(s,r)}r^{h+l-2}\\&& \nonumber \mbox{(based on
the inequality (\ref{081606}), $r-\sigma'<r$ and the assumption $j=h+l-2\geq0$)}\\
&=&\frac {r}{\sigma'}\left(\frac1{r}||| U_{z}|||^T_{1,D(s,r)}r^{h-1}\right)\left(\frac1r|||
V_z|||^T_{1,D(s,r)}r^{l-1}\right)\nonumber\\
&\leq&\frac {r}{\sigma'}\left(\frac1{r}|||
U_{z}|||^T_{p,D(s,r,r)}\right)\left(\frac1r|||V_y|||^T_{p,D(s,r,r)}\right)\qquad \mbox{(by
(\ref{091403}) in Definition \ref{071803})}\nonumber\\
&\leq&\frac {r}{\sigma'}|||
X_{U}|||^T_{p,D(s,r,r)}|||X_V|||^T_{p,D(s,r,r)}\nonumber\label{091803}.
\end{eqnarray}

\textbf{Step 2. Estimate $p$-tame norm of the terms in case 3, which
are
\begin{eqnarray*}
&&U_{zx}V_y=(U_{qx}V_y,U_{\bar qx}V_y),\qquad V_{zx}U_y=(V_{qx}
U_y,V_{\bar qx}U_y) \\
&&U_{zy}V_x=(U_{qy}V_x,U_{\bar qy}V_x),\qquad V_{zy}U_x=(V_{qy}
U_x,V_{\bar qy}U_x).
\end{eqnarray*}}
Firstly, we will estimate
$|||U_{zx}V_y|||^T_{p,D(s-\sigma,r-\sigma',r-\sigma')}$. Let
$\tilde m=h+l-1$. Following the proof of (\ref{080309}), we obtain
\begin{eqnarray}
||\widetilde{\lfloor U_{zx}
V_{y}\rceil}_{D(s-\sigma,r-\sigma')}(z^{(1)},\dots,z^{(\tilde m)})||_p
\leq\frac{1}{e\sigma}||| U_{z}|||^T_{p,D(s,r)}||| V_{y
}|||_{D(s,r)}f(z),\label{080609}
\end{eqnarray}where
\begin{equation*}
f(z)=\frac1{\tilde m!}\frac1{h-1}\sum_{\tau_{\tilde
m}}\sum_{j=1}^{h-1}||z^{\tau_{\tilde
m}(1)}||_1\cdots||z^{\tau_{\tilde m}(j-1)}||_1 ||z^{\tau_{\tilde
m}(j)}||_p||z^{\tau_{\tilde m}(j+1)}||_1\cdots||z^{\tau_{\tilde{m}}
({\tilde{m}})}||_1 ,\end{equation*} and using the inequalities
\begin{eqnarray*} &&||\widetilde{\lfloor
U_{z}\rceil}_{D(s,r)}(z^{(\tau_{\tilde
m}(1))},\dots,z^{(\tau_{\tilde m}(h-1))})||_p\\&\leq& |||
U_{z}|||^T_{p,D(s,r)}\left(\frac1{h-1}\sum_{j=1}^{h-1}||z^{\tau_{\tilde
m}(1)}||_1\cdots||z^{\tau_{\tilde m}(j-1)}||_1 ||z^{\tau_{\tilde
m}(j)}||_p||z^{\tau_{\tilde m}(j+1)}||_1\cdots||z^{\tau_{\tilde{m}}
({h-1})}||_1\right)
\end{eqnarray*}
and
\begin{eqnarray*}
\widetilde{\lfloor V_{y}\rceil}_{D(s,r)}(z^{(\tau_{\tilde
m}(h))},\dots,z^{(\tau_{\tilde m}(\tilde m)})\leq||| V_{y
}|||_{D(s,r)}||z^{\tau_{\tilde m}(h)}||_1\cdots||z^{\tau_{\tilde{m}}
({\tilde{m}})}||_1,
\end{eqnarray*}since $U(x,y,z;\xi)$ and $V(x,y,z;\xi)$ have $p$-tame property on the domain
$D(s,r,r)\times\Pi$. If
\begin{equation}f(z)=||(z^{\tilde m})||_{p,1}
 \end{equation}(which will be prove in Lemma \ref{012101}), then
in view of (\ref{091402}) in Definition \ref{071803} and the
inequality (\ref{080609}), we obtain
\begin{equation}\label{080701}
|||U_{zx}V_y|||_{p,D(s-\sigma,r-\sigma')}^{T}\leq \frac{1}{e\sigma}|||
U_{z}|||^T_{p,D(s,r)}||| V_{y }|||_{D(s,r)}.
\end{equation}
In particular, when $p=1$, (\ref{080701}) reads
\begin{equation}\label{080702}
|||U_{zx}V_y|||_{1,D(s-\sigma,r-\sigma')}^{T}\leq \frac{1}{e\sigma}|||
U_{z}|||^T_{1,D(s,r)}||| V_{y }|||_{D(s,r)}.
\end{equation}
Hence,
\begin{eqnarray}
&&\nonumber\frac1{r-\sigma'}|||U_{zx}V_y|||_{p,D(s-\sigma,r-\sigma',r-\sigma')}^{T}\\
&=&\nonumber\frac1{r-\sigma'}\max\left\{|||U_{zx}V_y|||_{p,D(s-\sigma,r-\sigma')}^{T},
|||U_{zx}V_y|||_{1,D(s-\sigma,r-\sigma')}^{T}\right\}(r-\sigma')^{h+l-1}\\&&\nonumber\mbox{(in
view of (\ref{091403}) in Definition \ref{071803})}\\
\nonumber&\leq&
\frac{1}{e\sigma}\max\left\{|||U_{z}|||_{p,D(s,r)}^{T}||| V_{y
}|||_{D(s,r)},|||U_{z}|||_{1,D(s,r)}^{T}||| V_{y
}|||_{D(s,r)}\right\}(r-\sigma')^{h+l-2}\\
&&\mbox{(based on the inequalities (\ref{080701}) and (\ref{080702}))}
\nonumber\\
\nonumber&\leq&
\frac{1}{e\sigma}\cdot \frac{r}{(r-\sigma')}\cdot\frac 1r\max\left\{|||U_{z}|||_{p,D(s,r)}^{T},|||U_{z}|||_{1,D(s,r)}^{T}\right\}||| V_{y
}|||_{D(s,r)}r^{h+l-1}\\
&&\mbox{(based on $r-\sigma'<r$, and the assumption $h+l-1\geq0$ otherwise $U_z=0$)}
\nonumber\\&\leq&\nonumber\frac{2}{e\sigma}\left(\frac1r\max\left
\{|||U_{z}|||_{p,D(s,r)}^{T}, |||U_{z}|||_{1,D(s,r)}^{T}
\right\}r^{h-1}\right)|||V_{y
}|||_{D(s,r)}r^{l}\\
&&\mbox{(based on $0<\sigma'<r/2$ implying $r/2<r-\sigma'<r$)}\nonumber\\&=&\frac{2}{e\sigma}\left(\frac1r|||
U_{z}|||^T_{p,D(s,r,r)}\right)||| V_{y }|||_{D(s,r,r)}\qquad
\mbox{(by (\ref{091403}) in Definition \ref{071803})}\label{042805}.
\end{eqnarray}
Following
the proof of the inequality (\ref{042805}), we obtain
 \begin{eqnarray}
 \frac1{r-\sigma'}|||V_{zx}U_y|||_{p,D(s-\sigma,r-\sigma',r-\sigma')}^{T}&\leq&
\frac{2}{e\sigma}\left(\frac1r|||
V_{z}|||^T_{p,D(s,r,r)}\right)||| U_{y }|||_{D(s,r,r)},\\
 \frac1{r-\sigma'}|||U_{zy}V_x|||_{p,D(s-\sigma,r-\sigma',r-\sigma')}^{T}&\leq&
\frac{2r}{\sigma'}||| U_{z}|||^T_{p,D(s,r,r)}\left(\frac1{r^2}|||
V_{x
}|||_{D(s,r,r)}\right),\label{091806}\\
 \frac1{r-\sigma'}|||V_{zy}U_x|||_{p,D(s-\sigma,r-\sigma',r-\sigma')}^{T}&\leq&
\frac{2r}{\sigma'}||| V_{z}|||^T_{p,D(s,r,r)}\left(\frac1{r^2}|||
U_{x }|||_{D(s,r,r)}\right),\label{091807}
\end{eqnarray}where, to prove the
inequalities (\ref{091806}) and (\ref{091807}), we use the
generalized Cauchy estimate (\ref{071902})
instead of (\ref{071901}) in Lemma \ref{050403}.

\textbf{Step 3. Estimate
$p$-tame norm of the terms in case 4, which are
$$U_{zz}\cdot V_{z}=\sqrt{-1}(\sum_{j\geq1}(U_{\bar{q}q_j}V_{\bar
q_j}-U_{\bar{q}\bar{q}_j}V_{q_j}),\sum_{j\geq1}(U_{qq_j}V_{\bar
q_j}-U_{q\bar q_j}V_{q_j}))$$ and
$$V_{zz}\cdot U_{z}=\sqrt{-1}(\sum_{j\geq1}(V_{\bar{q}q_j}U_{\bar
q_j}-V_{\bar{q}\bar{q}_j}U_{q_j}),\sum_{j\geq1}(V_{qq_j}U_{\bar
q_j}-V_{q\bar q_j}U_{q_j})).$$}

Firstly, note an important fact that $\lfloor U\rceil_{D(s,r)}(z)$
and $\lfloor V\rceil_{D(s,r)}(z)$ are two Hamiltonian depending only
on the normal variable $z$ (independent of the tangent variables $(x,y)$ and parameter $\xi$). Moreover,
\begin{equation*}
{\left\{\lfloor U\rceil_{D(s,r)}(z),\lfloor
V\rceil_{D(s,r)}(z)\right\}}_z=\lfloor
U_{zz}\rceil_{D(s,r)}\cdot\lfloor V_{z}\rceil_{D(s,r)}-\lfloor
V_{zz}\rceil_{D(s,r)}\cdot \lfloor U_{z}\rceil_{D(s,r)}.
\end{equation*}
Following the proof of Lemma 4.12 in \cite{BG}, we obtain
\begin{eqnarray}
&&\nonumber|||{\lfloor U_{zz}\rceil_{D(s,r)}\cdot \lfloor
V_{z}\rceil}_{D(s,r)}|||^T_{p,D(s,r)}+|||{\lfloor
V_{zz}\rceil_{D(s,r)}\cdot \lfloor
U_{z}\rceil}_{D(s,r)}|||^T_{p,D(s,r)}\\&\leq&\nonumber
(h+l-2)|||\lfloor U_z\rceil_{D(s,r)}|||^T_{p,D(s,r)}|||\lfloor
V_z\rceil_{D(s,r)}|||^T_{p,D(s,r)}\\&=&\label{082204} (h+l-2)|||
U_z|||^T_{p,D(s,r)}|||V_z|||^T_{p,D(s,r)}.
\end{eqnarray}
In particular, when $p=1$, the inequality (\ref{082204}) reads
\begin{eqnarray}
&&\nonumber|||{\lfloor U_{zz}\rceil_{D(s,r)}\cdot \lfloor
V_{z}\rceil}_{D(s,r)}|||^T_{1,D(s,r)}+|||{\lfloor
V_{zz}\rceil_{D(s,r)}\cdot \lfloor
U_{z}\rceil}_{D(s,r)}|||^T_{1,D(s,r)}\\&\leq&\label{092910}
(h+l-2)||| U_z|||^T_{1,D(s,r)}|||V_z|||^T_{1,D(s,r)}.
\end{eqnarray}By Lemma \ref{051701}, we have the following
inequalities
\begin{equation}\label{082205}
||\widetilde{\lfloor U_{zz}\cdot
V_{z}\rceil}_{D(s,r)}(z^{(1)},\dots,z^{(h+l-3)})||_p\leq
||\widetilde{\lfloor U_{zz}\rceil_{D(s,r)}\cdot \lfloor
V_{z}\rceil}_{D(s,r)}(z^{(1)},\dots,z^{(h+l-3)})||_p,
\end{equation}
and
\begin{equation}\label{082206}
||\widetilde{\lfloor V_{zz}\cdot
U_{z}\rceil}_{D(s,r)}(z^{(1)},\dots,z^{(h+l-3)})||_p\leq
||\widetilde{\lfloor V_{zz}\rceil_{D(s,r)}\cdot \lfloor
U_{z}\rceil}_{D(s,r)}(z^{(1)},\dots,z^{(h+l-3)})||_p.
\end{equation}Then
in view of (\ref{091402}) in Definition \ref{071803} and the
inequalities (\ref{082204})-(\ref{082206}), we obtain
\begin{equation}\label{082208}
|||U_{zz}\cdot V_{z}|||_{p,D(s,r)}^T+|||V_{zz}\cdot
U_{z}|||_{p,D(s,r)}^T\leq (h+l-2)|||
U_z|||^T_{p,D(s,r)}|||V_z|||^T_{p,D(s,r)}
\end{equation}
and
\begin{equation}\label{082208'}
|||U_{zz}\cdot V_{z}|||_{1,D(s,r)}^T+|||V_{zz}\cdot
U_{z}|||_{1,D(s,r)}^T\leq (h+l-2)|||
U_z|||^T_{1,D(s,r)}|||V_z|||^T_{1,D(s,r)}
\end{equation}respectively.
 Hence, we have
\begin{eqnarray}\label{091808}
\nonumber&&\frac1{r-\sigma'}\left(|||U_{zz}\cdot
V_{z}|||_{p,D(s-\sigma,r-\sigma',r-\sigma')}^T+|||V_{zz}\cdot
U_{z}|||_{p,D(s-\sigma,r-\sigma',r-\sigma')}^T\right)\\
\nonumber&=&\frac1{r-\sigma'}\max\left\{|||U_{zz}\cdot
V_{z}|||_{p,D(s-\sigma,r-\sigma')}^T,|||U_{zz}\cdot
V_{z}|||_{1,D(s-\sigma,r-\sigma')}^T\right\}(r-\sigma')^{h+l-3}\\
\nonumber&&+\frac1{r-\sigma'}\max\left\{|||V_{zz}\cdot
U_{z}|||_{p,D(s-\sigma,r-\sigma')}^T,|||V_{zz}\cdot
U_{z}|||_{1,D(s-\sigma,r-\sigma')}^T\right\}(r-\sigma')^{h+l-3}\\
&&\mbox{(in view of (\ref{091403}) in Definition
\ref{071803})}\nonumber \\
&\leq&\max\left\{|||U_{zz}\cdot V_{z}|||_{p,D(s,r)}^T,|||U_{zz}\cdot
V_{z}|||_{1,D(s,r)}^T\right\}(r-\sigma')^{h+l-4}\nonumber\\
\nonumber&&+\max\left\{|||V_{zz}\cdot
U_{z}|||_{p,D(s,r)}^T,|||V_{zz}\cdot
U_{z}|||_{1,D(s,r)}^T\right\}(r-\sigma')^{h+l-4}\qquad \mbox{(in view of (\ref{011709}))}\\
\nonumber&\leq&2(h+l-2)\max\left\{|||U_{z}|||_{p,D(s,r)}^T|||
V_{z}|||_{p,D(s,r)}^T,|||U_{z}|||_{1,D(s,r)}^T|||
V_{z}|||_{1,D(s,r)}^T\right\}(r-\sigma')^{h+l-4}\\&&\nonumber\mbox{(based on the inequalities (\ref{082208}) and (\ref{082208'}))}\\
\nonumber&\leq&\frac2{(r-\sigma')\sigma'}\max\left\{|||U_{z}|||_{p,D(s,r)}^T|||
V_{z}|||_{p,D(s,r)}^T,|||U_{z}|||_{1,D(s,r)}^T|||
V_{z}|||_{1,D(s,r)}^T\right\}r^{h+l-2} \\ && \nonumber(\mbox{using
the inequality}\ k(r-\sigma')^{k-1}\leq\frac{r^{k}}{\sigma'})
\\&\leq& \frac2 {(r-\sigma')\sigma'}\nonumber
|||U_z|||^{T}_{p,D(s,r,r)}|||V_{z}|||^{T}_{p,D(s,r,r)}\qquad \mbox{(in
view of (\ref{091403}) in Definition \ref{071803})}\nonumber
\\&\leq& \frac {4r}{\sigma'}
\left(\frac1r|||U_z|||^{T}_{p,D(s,r,r)}\right)\left(\frac1r|||V_z|||^{T}_{p,D(s,r,r)}\right),
\label{021011}
\end{eqnarray}
where the last inequality is based on $0<\sigma'<r/2$ implies $r/2<r-\sigma'<r.$
\\
\textbf{Step 4. Estimate
$$|||\{U,V\}_z|||^T_{p,D(s-\sigma,r-\sigma',r-\sigma')}\quad \mbox{
and}\quad
|||X_{\{U,V\}}|||^T_{p,D(s-\sigma,r-\sigma',r-\sigma')}.$$} By
the inequalities (\ref{042805})-(\ref{091807}) and (\ref{021011}),
we obtain
 \begin{eqnarray}
\nonumber&&\frac1{r-\sigma'}|||\{U,V\}_z|||_{p,D(s-\sigma,r-\sigma',r-\sigma')}^{T}\\
&\leq& 12\max\left\{\frac {1}{\sigma},\frac
{r}{\sigma'}\right\}|||
X_U|||_{p,D(s,r,r)}^T|||X_V|||_{p,D(s,r,r)}^{T}\label{042808}.
\end{eqnarray}
By Definition \ref{051703} and in view of the inequalities
(\ref{042806}), (\ref{042807}) and (\ref{042808}), it is easy to see
that
\begin{equation}\label{080207}
|||X_{ \{U,V\}}|||_{p,D(s-\sigma,r-\sigma',r-\sigma')}^T\leq
C\max\left\{\frac 1{\sigma},\frac
{r}{\sigma'}\right\}|||X_U|||_{p,D(s,r,r)}^T|||X_V|||_{p,D(s,r,r)}^T,
\end{equation}
where $C>0$ is an absolute constant.

Finally, consider general Hamiltonian
$$U(x,y,z;\xi)=\sum_{h\geq0}U_h(x,y,z;\xi)$$ and
$$V(x,y,z;\xi)=\sum_{l\geq0}V_l(x,y,z;\xi),$$where
\begin{equation*}\label{102301}U_h(x,y,z;\xi)=\sum_{\beta\in\mathbb{N}^{\bar{\mathbb{Z}}},
|\beta|=h}U^{\beta}(x,y;\xi)z^{\beta} \end{equation*}
and
\begin{equation*}\label{102302}V_l(x,y,z;\xi)=\sum_{\beta\in\mathbb{N}^{\bar{\mathbb{Z}}},|\beta|=l}
V^{\alpha\beta}(x,y;\xi)z^{\beta}. \end{equation*}By a direct calculation,
we have
\begin{eqnarray*}&&|||X_{\{U,V\}}|||_{p,D(s-\sigma,r-\sigma',r-\sigma')\times\Pi}^T\\
&=&|||X_{\{\sum_{h\geq0}U_{h},\sum_{l\geq0}V_{l}\}}
|||_{p,D(s-\sigma,r-\sigma',r-\sigma')\times\Pi}^T\\
&=&|||\sum_{h,l\geq0}X_{\{U_{h},V_{l}\}}
|||_{p,D(s-\sigma,r-\sigma',r-\sigma')\times\Pi}^T\\
&\leq&\sum_{h,l\geq0}|||X_{\{U_{h},
V_{l}\}}|||_{p,D(s-\sigma,r-\sigma',r-\sigma')\times\Pi}^T\\
&\leq&\sum_{h,l\geq0}C\max\left\{\frac 1{\sigma^{n+1}},\frac
{r}{\sigma'}\right\} |||X_{U_{h}}|||_{p,D(s,r,r)\times\Pi}^T
|||X_{V_{l}}|||_{p,D(s,r,r)\times\Pi}^T\qquad \mbox{(by (\ref{080207}))}\\
&=&C\max\left\{\frac 1{\sigma},\frac
{r}{\sigma'}\right\}\left(\sum_{h\geq0}
|||X_{U_{h}}|||_{p,D(s,r,r)\times\Pi}^T\right)
\left(\sum_{l\geq0}|||X_{V_{l}}
|||_{p,D(s,r,r)\times\Pi}^T\right)\\
&=&C\max\left\{\frac 1{\sigma},\frac
{r}{\sigma'}\right\}|||X_{U}|||_{p,D(s,r,r)\times\Pi}^T
|||X_{V}|||_{p,D(s,r,r)\times\Pi}^T,
\end{eqnarray*}
where the last equality is based on Definition \ref{051703}.
\end{proof}
\begin{rem}
In view of the estimate (\ref{091120}), the coefficient \begin{equation}\label{000007}C\max\left\{\frac 1{\sigma},\frac
{r}{\sigma'}\right\}
\end{equation}
is necessary in normal form iteration (but not important in KAM iteration). The reason is that while constructing a partial normal form of higher order (see Section 6), $\sigma\approx s$, $r\approx \rho$ and $\sigma'\approx\rho/\mathcal{M}$, then (\ref{000007}) reads a constant independent of $\rho$, but depending on $n,s$ and $\mathcal{M}$, which make the first step of the normal form iteration work.
\end{rem}
\color{red}
The theorem below show that $p$-tame property persists under Hamiltonian phase flow.
\color{black}\begin{thm}\label{081914}
Consider two Hamiltonian $U(x,y,z;\xi)$ and $V(x,y,z;\xi)$
satisfying $p$-tame property on the domain $D(s,r,r)\times\Pi$ for
some $0<s,r\leq1$. Given $0<\sigma<s, 0<\sigma'<r/2$, suppose
\begin{equation}\label{081112}|||X_U|||_{p,D(s,r,r)\times\Pi}^T\leq \frac1{2A},
\end{equation}
where
\begin{equation}\label{090504}A=4Ce\max\left\{\frac{1}{\sigma},\frac{r}{\sigma'}\right\}
\end{equation}
and $C>0$ is the constant given in (\ref{091120}) in Theorem
\ref{021102}. Then for each $|t|\leq 1$, we have
\begin{equation*}
|||X_{V\circ
X_U^t}|||_{p,D(s-\sigma,r-\sigma',r-\sigma')\times\Pi}^T\leq
2|||X_V|||_{p,D(s,r,r)\times\Pi}^T.
\end{equation*}
\end{thm}
\begin{proof}
Let
\begin{equation*}
W^{(0)}(x,y,z;\xi)=V(x,y,z;\xi),
\end{equation*}
and \begin{equation*} W^{(j)}(x,y,z;\xi)=\{W^{(j-1)},U\}(x,y,z;\xi),
\qquad j\geq1.
\end{equation*}
Hence,
\begin{eqnarray*}
W^{(1)}(x,y,z;\xi)&=&\{W^{(0)},U\}(x,y,z;\xi)=\{V,U\}(x,y,z;\xi),\\
W^{(2)}(x,y,z;\xi)&=&\{W^{(1)},U\}(x,y,z;\xi)=\{\{V,U\},U\}(x,y,z;\xi),\\
&\dots&\\
W^{(j)}(x,y,z;\xi)&=&\{\dots\{\{V,U\},U\},U\},\dots\}, U\}(x,y,z;\xi)\qquad (\mbox{there
are $j$ times U}).
\end{eqnarray*}
For $j\geq1$, let $\sigma_j=\frac{\sigma}{j}$ and
$\sigma_j'=\frac{\sigma'}j$. Hence, we obtain
\begin{eqnarray}\nonumber
&&|||X_{W^{(j)}}|||_{p,D(s-\sigma,r-\sigma',r-\sigma')\times\Pi}^{T}\\
\nonumber&=&
|||X_{W^{(j)}}|||_{p,D(s-j\sigma_j,r-j\sigma'_j,r-j\sigma'_j)\times\Pi}^{T}\\&\leq&\nonumber
\left(4C\max\left\{\frac{1}{\sigma_j},\frac{r}{\sigma'_j}\right\}\right)^{j}
\left(|||X_{V}
|||_{p,D(s,r,r)\times\Pi}^{T}\right)\left(|||X_{U}|||_{p,D(s,r,r)\times\Pi}^{T}\right)^j\\
&&\mbox{(based on Theorem {\ref{021102}} and (\ref{011708}))}\nonumber\\
&=&\label{081110}j^j\left(4C\max\left\{\frac{1}{\sigma},\frac{r}{\sigma'}\right\}\right)^{j}
\left(|||X_{V}
|||_{p,D(s,r,r)\times\Pi}^{T}\right)\left(|||X_{U}|||_{p,D(s,r,r)\times\Pi}^{T}\right)^j.
\end{eqnarray}
Using the inequality
\begin{equation}
j^j<j!e^j
\end{equation}and in view of the inequality (\ref{081110}), we have
\begin{eqnarray}\nonumber
&&\frac1{j!}|||X_{W^{(j)}}|||_{p,D(s-\sigma,r-\sigma',r-\sigma')\times\Pi}^{T}\\&\leq&\nonumber
e^{j}\left(4C\max\left\{\frac{1}{\sigma},\frac{r}{\sigma'}\right\}\right)^{j}
\left(|||X_{V}
|||_{p,D(s,r,r)\times\Pi}^{T}\right)\left(|||X_{U}|||_{p,D(s,r,r)\times\Pi}^{T}\right)^j\\
&=&\nonumber\left(|||X_{V}
|||_{p,D(s,r,r)\times\Pi}^{T}\right)\left(A|||X_{U}|||_{p,D(s,r,r)\times\Pi}^{T}\right)^j\
\mbox{(in view of
(\ref{090504}))}\\
&=&\label{081111}{2^{-j}}|||X_{V}
|||_{p,D(s,r,r)\times\Pi}^{T}\qquad\qquad\qquad(\mbox{based on the
inequality (\ref{081112})}).\end{eqnarray} Expand the Hamiltonian
$V\circ X_U^t(x,y,z;\xi)$ into Taylor series about $t$ at $t=0$, and
we have
\begin{equation}\label{081108}
V\circ X_U^t(x,y,z;\xi)=\sum_{j\geq 0}
\frac{t^{j}}{j!}W^{(j)}(x,y,z;\xi).
\end{equation}
Then
\begin{eqnarray*}
&&|||X_{V\circ
X_U^t}|||_{p,D(s-\sigma,r-\sigma',r-\sigma')\times\Pi}^{T}\\
&=&|||X_{\sum_{j\geq 0} \frac{t^{j}}{j!}W^{(j)}}|||_{p,D(s-\sigma,r-\sigma',r-\sigma')\times\Pi}^{T}\qquad(\mbox{in view of (\ref{081108})})\\
&\leq&\sum_{j\geq 0}
\frac{1}{j!}|||X_{W^{(j)}}|||_{p,D(s-\sigma,r-\sigma',r-\sigma')\times\Pi}^{T}\qquad \mbox{(in view of $|t|\leq 1$)}\\
&\leq&\sum_{j\geq 0} {2^{-j}}|||X_{V}
|||_{p,D(s,r,r)\times\Pi}^{T}\qquad\qquad\qquad (\mbox{based on the inequality (\ref{081111})})\\
&=&2|||X_{V} |||_{p,D(s,r,r)\times\Pi}^{T}.
\end{eqnarray*}
\end{proof}
\color{red} Based on the theorem below, we will estimate the $p$-tame norm of the solution of homological equation.
\color{black}
\begin{thm}\label{0005}Consider two Hamiltonian
$$U(x,y,z;\xi)=\sum_{\alpha\in\mathbb{N}^n,\beta\in\mathbb{N}^{\bar{\mathbb{Z}}}}
U^{\alpha\beta}(x;\xi)y^{\alpha}z^{\beta}$$ and
$$V(x,y,z;\xi)=\sum_{\alpha\in\mathbb{N}^n,\beta\in\mathbb{N}^{\bar{\mathbb{Z}}}}
V^{\alpha\beta}(x;\xi)y^{\alpha}z^{\beta}.$$ Suppose $V(x,y,z;\xi)$
has $p$-tame property on the domain $D(s,r,r)\times \Pi$, i.e
$$|||X_V|||_{p,D(s,r,r)\times\Pi}^T<\infty.$$ For each
$\alpha\in\mathbb{N}^n,\beta\in\mathbb{N}^{\bar{\mathbb{Z}}},k\in
\mathbb{Z}^n, j\geq 1$ and some fixed constant $\tau>0$, assume the
following inequality holds
\begin{equation}\label{0007}
|\widehat{U^{\alpha\beta}}(k;\xi)|+|\partial_{\xi_j}\widehat{U^{\alpha\beta}}(k;\xi)|\leq
(|k|+1)^{\tau}(|\widehat{V^{\alpha\beta}}(k;\xi)|+|\partial_{\xi_j}\widehat{V^{\alpha\beta}}(k;\xi)|)
,\end{equation} where $\widehat{U^{\alpha\beta}}(k;\xi)$ and
$\widehat{V^{\alpha\beta}}(k;\xi)$ are the $k$-th Fourier
coefficients of $U^{\alpha\beta}(x;\xi)$ and
$V^{\alpha\beta}(x;\xi)$, respectively. Then, $U(x,y,z;\xi)$ has
$p$-tame property on the domain $D(s-\sigma,r,r)\times\Pi$ for
$0<\sigma<s$. Moreover, we have
\begin{equation}\label{0012} |||X_U|||_{p,D(s-\sigma,r,r)\times\Pi}^T\leq
\frac{c}{\sigma^{\tau}} |||X_V|||_{p,D(s,r,r)\times\Pi}^T,
\end{equation}
where $c>0$ is a constant depending on $s$ and $\tau$.
\end{thm}
\begin{proof}Without loss of generality, we
suppose $U(x,y,z;\xi)=U_h(x,y,z;\xi)$ and
$V(x,y,z;\xi)=V_h(x,y,z;\xi)$.

Firstly, we will estimate $|||U_x|||_{D(s-\sigma,r,r)\times\Pi}$.
For $1\leq i \leq n$, note
$$U_{x_i}(x,y,z;\xi)=\sum_{\alpha\in\mathbb{N}^n,\beta\in\mathbb{N}^{\bar{\mathbb{Z}}},|\beta|=h}
U_{x_i}^{\alpha\beta}(x;\xi)y^{\alpha}z^{\beta}$$ and
$$V_{x_i}(x,y,z;\xi)=\sum_{\alpha\in\mathbb{N}^n,\beta\in\mathbb{N}^{\bar{\mathbb{Z}}},|\beta|=h}
V_{x_i}^{\alpha\beta}(x;\xi)y^{\alpha}z^{\beta}.$$ Expand
$U_{x_i}^{\alpha\beta}(x;\xi)$ and $V_{x_i}^{\alpha\beta}(x;\xi)$
into Fourier series, which are
\begin{equation*}
U_{x_i}^{\alpha\beta}(x;\xi)=\sum_{k\in\mathbb{Z}^n}\widehat{U_{x_i}^{\alpha\beta}}(k,\xi)e^{\sqrt{-1}\langle
k,x\rangle}
\end{equation*}
and
\begin{equation*}
V_{x_i}^{\alpha\beta}(x;\xi)=\sum_{k\in\mathbb{Z}^n}\widehat{V_{x_i}^{\alpha\beta}}(k,\xi)e^{\sqrt{-1}\langle
k,x\rangle}.
\end{equation*}Note that
\begin{equation*}
\widehat{U_{x_i}^{\alpha\beta}}(k,\xi)=k_i\widehat{U^{\alpha\beta}}(k,\xi),
\end{equation*}and
\begin{equation*}
\widehat{V_{x_i}^{\alpha\beta}}(k,\xi)=k_i\widehat{V^{\alpha\beta}}(k,\xi),
\end{equation*}where $\widehat{U^{\alpha\beta}}(k,\xi)$ and
$\widehat{V^{\alpha\beta}}(k,\xi)$ are $k$-th Fourier coefficients
of $U^{\alpha\beta}(x;\xi)$ and $V^{\alpha\beta}(x;\xi)$
respectively. By the inequality (\ref{0007}), for $j\geq1$,
we have
\begin{equation}\label{082201}
|\widehat{U_{x_i}^{\alpha\beta}}(k;\xi)|+|\partial_{\xi_{j}}\widehat{U_{x_i}^{\alpha\beta}}(k;\xi)|\leq
(|k|+1)^{\tau}(|\widehat{V_{x_i}^{\alpha\beta}}(k;\xi)|+|\partial_{\xi_{j}}\widehat{V_{x_i}^{\alpha\beta}}(k;\xi)|)
.\end{equation}
By a simple calculation, we obtain
\begin{eqnarray}\nonumber&&\sum_{k\in\mathbb{Z}^n}
\left(|\widehat{U_{x_i}^{\alpha\beta}}(k;\xi)|+|\partial_{\xi_j}\widehat{U_{x_i}^{\alpha\beta}}(k;\xi)|\right)e^{|k|(s-\sigma)}
\\&\leq&\nonumber\sum_{k\in\mathbb{Z}^n}(|k|+1)^{\tau}
\left(|\widehat{V_{x_i}^{\alpha\beta}}(k;\xi)|+|\partial_{\xi_j}\widehat{V_{x_i}^{\alpha\beta}}(k;\xi)|\right)e^{|k|(s-\sigma)}\\
&&\nonumber \mbox{(according to the inequality (\ref{082201}))}
\\&=&\nonumber\sum_{k\in\mathbb{Z}^n}
\left(|\widehat{V_{x_i}^{\alpha\beta}}(k;\xi)|+|\partial_{\xi_j}
\widehat{V_{x_i}^{\alpha\beta}}(k;\xi)|\right)e^{|k|s}(|k|+1)^{\tau}e^{-|k|\sigma}\\
&\leq&\nonumber\left(\sum_{k\in\mathbb{Z}^n}
\left(|\widehat{V_{x_i}^{\alpha\beta}}(k;\xi)|+|\partial_{\xi_j}
\widehat{V_{x_i}^{\alpha\beta}}(k;\xi)|\right)e^{|k|s}\right)\left(\sup_{k\in\mathbb{Z}^n}(|k|+1)^{\tau}e^{-|k|\sigma}
\right)\\
&\leq&\label{092702}
\frac{c}{\sigma^{\tau}}\left(\sum_{k\in\mathbb{Z}^n}
\left(|\widehat{V_{x_i}^{\alpha\beta}}(k;\xi)|+|\partial_{\xi_j}
\widehat{V_{x_i}^{\alpha\beta}}(k;\xi)|\right)e^{|k|s}\right)
,\end{eqnarray} where $c>0$ is a constant depending on $s$ and $\tau$.
By the inequality (\ref{082201}) and the definition
of the norm $||\cdot||_{D(s)\times\Pi}$ (see (\ref{092704}) in
Definition \ref{090501}), we obtain
\begin{equation}\label{080810}
||U_{x_i}^{\alpha\beta}||_{D(s-\sigma)\times\Pi}\leq
\frac{c}{\sigma^{\tau}} ||V_{x_i}^{\alpha\beta}||_{D(s)\times\Pi}.
\end{equation}
Hence, in view of (\ref{092703}) in Definition \ref{060602},
we obtain
\begin{equation}\label{080811}
||U^{\beta}_{x_i}||_{D(s-\sigma,r)\times\Pi}\leq
\frac{c}{\sigma^{\tau}} ||V_{x_i}^{\beta}||_{D(s,r)\times\Pi},
\end{equation}
where
$$U^{\beta}_{x_i}(x,y;\xi)=\sum_{\alpha\in\mathbb{N}^n}U_{x_i}^{\alpha\beta}(x;\xi)y^{\alpha}$$
and
$$V^{\beta}_{x_i}(x,y;\xi)=\sum_{\alpha\in\mathbb{N}^n}V_{x_i}^{\alpha\beta}(x;\xi)y^{\alpha}.$$
Then, following the proof of (\ref{042803}), it is easy to verify that
\begin{eqnarray}|||U_{x}|||_{D(s-\sigma,r,r)\times\Pi}
\leq\label{080813}
\frac{c}{\sigma^{\tau}}|||V_{x}|||_{D(s,r,r)\times\Pi}.
\end{eqnarray}
Similar to the proof of (\ref{080813}), we get
\begin{eqnarray}|||U_{y}|||_{D(s-\sigma,r,r)\times\Pi}
\leq\label{080809}
\frac{c}{\sigma^{\tau}}|||V_{y}|||_{D(s,r,r)\times\Pi},
\end{eqnarray}
and
\begin{eqnarray}|||U_{z}|||_{D(s-\sigma,r,r)\times\Pi}
\leq\label{012001}
\frac{c}{\sigma^{\tau}}|||V_{z}|||_{D(s,r,r)\times\Pi}.
\end{eqnarray}
Finally, using the inequalities (\ref{080813})-(\ref{012001}), and in view of Definition \ref{051703}, we
obtain the estimate
\begin{equation*}|||X_U|||_{p,D(s-\sigma,r,r)\times\Pi}^T\leq
\frac{c}{\sigma^{\tau}}|||X_V|||_{p,D(s,r,r)\times\Pi}^T.
\end{equation*}
\end{proof}
As in \cite{P1}, define the weighted norm of Hamiltonian vector field $X_U$ on the domain
$D(s,r,r)\times\Pi$ by
\begin{equation}\label{081101}
|||X_U|||_{\mathcal{P}^p,D(s,r,r)\times\Pi}=\sup_{(x,y,z;\xi)\in
D(s,r,r)\times\Pi}||X_U||_{\mathcal{P}^p,D(s,r,r)}.
\end{equation}
\color{red}Then the following two theorems give the estimate of the norm of time-$1$ map.
\color{black}
\begin{thm}\label{012002}Give a Hamiltonian
\begin{equation*}\label{100802}U(x,y,z;\xi)=\sum_{\beta\in\mathbb{N}^{\bar{\mathbb{Z}}}}
U^{\beta}(x,y;\xi)z^{\beta} \end{equation*} satisfying $p$-tame
property on the domain $D(s,r,r)\times\Pi$ for some $0<s,r\leq 1$.
Then we have
\begin{equation}\label{081105}
|||X_U|||_{\mathcal{P}^p,D(s,r,r)\times\Pi}\leq
|||X_U|||_{p,D(s,r,r)\times\Pi}^T.
\end{equation}
\end{thm}
\begin{proof}
Without loss of generality, assume
\begin{equation*}
U(x,y,z;\xi):=U_h(x,y,z;\xi)=\sum_{\beta\in\mathbb{N}^{\bar{\mathbb{Z}}},
|\beta|=h}U^{\beta}(x,y;\xi)z^{\beta}.
\end{equation*}Denote
\begin{equation*}
U^{\beta-1_j}_{z_j}(x,y;\xi)=\beta_jU^{\beta}(x,y;\xi),
\end{equation*}
where\begin{equation*}
\beta-1_j=(\dots,\beta_{j-1},\beta_j-1,\beta_{j+1},\dots).
\end{equation*}Based on the definition
$||\cdot||_{D(s,r)\times\Pi}$ (see (\ref{092703}) in Definition \ref{060602}), for
each $(x,y;\xi)\in D(s,r)\times\Pi$, we have
\begin{equation}\label{081610}
|U^{\beta-1_j}_{z_j}(x,y;\xi)|\leq
||U^{\beta-1_j}_{z_j}||_{D(s,r)\times\Pi}.
\end{equation}Then
\begin{eqnarray}
&&\nonumber|U_{z_j}(x,y,z;\xi)|\\
\nonumber&=&
\left|\sum_{\beta\in\mathbb{N}^{\bar{\mathbb{Z}}},|\beta|=h}
U^{\beta-1_j}_{z_j}(x,y;\xi)z^{\beta-1_j}\right|\\
&\leq&\nonumber\sum_{\beta\in\mathbb{N}^{\bar{\mathbb{Z}}},|\beta|=h}
|U^{\beta-1_j}_{z_j}(x,y;\xi)z^{\beta-1_j}|\\
&\leq&\nonumber{\sum_{\beta\in\mathbb{N}^{\bar{\mathbb{Z}}},|\beta|=h}
||U^{\beta-1_j}_{z_j}||_{D(s,r)\times\Pi}z^{\beta-1_j}}\\&&
\nonumber \mbox{(based on the inequality (\ref{081610}) and each
entry of $z^{(j)}$ is non-negative)}\\
&=&\label{090503}\widetilde{\lfloor
U_{z_j}\rceil}_{D(s,r)\times\Pi}(z^{(1)},\dots,z^{(h-1)})|_{z^{(j)}=z,
1\leq j\leq h-1}.
\end{eqnarray}
Moreover,
\begin{eqnarray}
&&\nonumber||U_z(x,y,z;\xi)||_{p}\\&\leq&
\nonumber||\widetilde{\lfloor
U_{z}\rceil}_{D(s,r)\times\Pi}(z^{(1)},\dots,z^{(h-1)})|_{z^{(j)}=z,
1\leq j\leq h-1}||_p\qquad \mbox{(based on the inequality (\ref{090503}))}\\
&\leq&
\nonumber|||U_z|||_{p,D(s,r)\times\Pi}^T||z||_p||z||_1^{h-2}\qquad
(\mbox{based on $z^{(j)}=z, 1\leq j\leq h-1$})\\
&\leq& \nonumber|||U_z|||_{p,D(s,r)\times\Pi}^T||z||_p^{h-1}\\
&\leq& \nonumber|||U_z|||_{p,D(s,r)\times\Pi}^Tr^{h-1}\\
&=&\label{081106}|||U_z|||^{T}_{p,D(s,r,r)\times\Pi}.
\end{eqnarray}
For each $(x,y,z;\xi)\in D(s,r,r)\times\Pi$, we obtain
\begin{eqnarray}
||U_x(x,y,z;\xi)||\leq \label{081107}|||U_x|||_{D(s,r,r)\times\Pi},
\end{eqnarray}
and
\begin{eqnarray}
||U_y(x,y,z;\xi)||\leq \label{081107'}|||U_y|||_{D(s,r,r)\times\Pi},
\end{eqnarray}
by following the proof of (\ref{081106}).
Hence,
\begin{eqnarray*}
&&\nonumber|||X_U|||_{\mathcal{P}^p,D(s,r,r)\times\Pi}\\&=&\nonumber\sup_{(x,y,z;\xi)\in
D(s,r,r)\times\Pi}\left(||U_y||+\frac1{r^2}||U_x||+\frac1r||U_z||_p\right)\qquad \mbox{(in view of the formula (\ref{081101}))}\\
&\leq&\nonumber
|||U_y|||_{D(s,r,r)\times\Pi}+\frac1{r^2}|||U_x|||_{D(s,r,r)\times\Pi}+\frac1r|||U_z|||^T_{p,D(s,r,r)\times\Pi}\\&&\nonumber
\mbox{(in view of the inequalities (\ref{081106})-(\ref{081107'}))}\\
&=&|||X_U|||^{T}_{p,D(s,r,r)\times\Pi}.\label{091103}
\end{eqnarray*}
\end{proof}
\begin{thm}\label{081503}
Suppose the Hamiltonian
$$U(x,y,z;\xi)=\sum_{\beta\in\mathbb{N}^{\bar{\mathbb{Z}}}}U^{\beta}(x,y;\xi)z^{\beta}$$
has $p$-tame property on the domain $D(s,r,r)\times\Pi$ for some
$0<s,r\leq 1$. Let $X_U^{t}$ be the phase flow generalized by the
Hamiltonian vector field $X_U$. Given $0<\sigma<s$ and
$0<\sigma'<r/2$, assume
\begin{equation*}\label{091810}|||X_U|||_{p,D(s,r,r)\times\Pi}^T<\min\{\sigma,\sigma'\}.
\end{equation*}
Then, for each $\xi\in\Pi$ and each $|t|\leq 1$, one has
\begin{equation}\label{091101}
||X_U^t-id||_{p,D(s-\sigma,r-\sigma',r-\sigma')}\leq
|||X_U|||_{p,D(s,r,r)\times\Pi}^T.
\end{equation}
\end{thm}
\begin{proof}
The inequality (\ref{091101}) can be proven directly based on Theorem \ref{012002} and following the proof of Lemma A.4 in \cite{P1}.
\end{proof}
\section{Proof of Theorem \ref{T1} and Corollary \ref{L1}}
\subsection{The $p$-tame property of the solution of homological equation}
\subsubsection{The derivation of homological equation}
Recall the perturbation of the integrable Hamiltonian (see
(\ref{100803}))
$$H(x,y,q,\bar q;\xi)=N(y,q,\bar q;\xi)+R(x,y,q,\bar q;\xi).$$ Denote
$R(x,y,q,\bar q;\xi)=R^{low}(x,y,q,\bar q;\xi)+R^{high}(x,y,q,\bar
q;\xi),$ where
  \begin{eqnarray*}
R^{low}(x,y,q,\bar
q;\xi)=\sum_{\alpha\in\mathbb{N}^n,\beta,\gamma\in\mathbb{N}^{\mathbb{N}},
2|\alpha|+|\beta|+|\gamma|\leq
2}R^{\alpha\beta\gamma}(x;\xi)y^{\alpha}q^{\beta}\bar q^{\gamma},
\end{eqnarray*}
and \begin{eqnarray*} R^{high}(x,y,q,\bar
q;\xi)=\sum_{\alpha\in\mathbb{N}^n,\beta,\gamma\in\mathbb{N}^{\mathbb{N}},2|\alpha|+|\beta|+|\gamma|\geq
3}R^{\alpha\beta\gamma}(x;\xi)y^{\alpha}q^{\beta}\bar q^{\gamma}.
\end{eqnarray*}The symplectic
coordinate change will be produced by the time-$1$ map
$X^t_F|_{t=1}$ of the Hamiltonian vector field $X_F$, where
$F(x,y,q,\bar q;\xi)$ is in the form\begin{eqnarray*} &&F(x,y,q,\bar
q;\xi)\\&=&F^{low}(x,y,q,\bar
q;\xi)\\&=&\sum_{\alpha\in\mathbb{N}^n,\beta,\gamma\in\mathbb{N}^{\mathbb{N}},2|\alpha|+|\beta|+|\gamma|\leq
2}F^{\alpha\beta\gamma}(x;\xi)y^{\alpha}q^{\beta}\bar q^{\gamma}.
\end{eqnarray*}
Using Taylor formula,
\begin{eqnarray*}
H_{+}&:=&H\circ
X_F^t|_{t=1}\\&=&H+\{H,F\}+\int^1_0(1-t)\{\{H,F\},F\}\circ
X_F^tdt\\
&=&N+\{N,F\}+\int^1_0(1-t)\{\{N,F\},F\}\circ X_F^tdt\\
&&+R^{low}+\int_0^1\{R^{low},F\}\circ X_F^tdt\\
&&+R^{high}+\{R^{high},F\}+\int^1_0(1-t)\{\{R^{high},F\},F\}\circ
X_F^tdt.
\end{eqnarray*}
Then we obtain the modified homological equation
\begin{equation}\label{16}
\{N,F\}+R^{low}+\{R^{high},F\}^{low}=N_+-N,
\end{equation}
where $N_{+}$ will be given in (\ref{081501}) below. If the homological equation
(\ref{16}) is solved, then the new perturbation term $R_{+}$ can be written as
\begin{eqnarray}
R_{+}&=&R^{high}+\{R^{high},F\}^{high}\label{011501}\\
&&+\int^1_0(1-t)\{\{N+R^{high},F\},F\}\circ X_F^tdt\label{011502}\\
&&+\int_0^1\{R^{low},F\}\circ X_F^tdt\label{011503}.
\end{eqnarray}
Note that we do not need to eliminate the terms in (\ref{011501}) at next step of KAM procedure, so (\ref{011501}) is not necessary to be small. On the other hand,
(\ref{011502}) is quadratic in $F$ and (\ref{011503}) contains the terms $R^{low}$ and $F$, which guarantee $R_{+}^{low}$ small enough. Therefore, the domain $D(s,r,r)$ is not required to shrink too fast such that we can obtain a non-degenerate normal form of order 2 directly.
\subsubsection{The solvability of homological equation (\ref{16})}
To solve the homological equation (\ref{16}), we should know what is
the term $\{R^{high},F\}^{low}$ exactly. Following Kuksin's notation
in \cite{K1}, write
\begin{equation*}
R^{low}=R^{x}+R^{y}+R^{1}+R^{2},
\end{equation*}where
\begin{eqnarray}\label{102405}R^{x}&=&
R^{x}(x;\xi)=\sum_{\alpha\in\mathbb{N}^n,\beta,\gamma\in\mathbb{N}^{\mathbb{N}},2|\alpha|+|\beta|+|\gamma|=0}R^{\alpha\beta\gamma}(x;\xi),\\
\label{102406}R^y&=&\sum_{j=1}^nR^{y_j}(x;\xi)y_j=\sum_{\alpha\in\mathbb{N}^n,\beta,\gamma\in\mathbb{N}^{\mathbb{N}},|\alpha|=1,\
|\beta|+|\gamma|=0 }R^{\alpha\beta\gamma}(x;\xi)y^{\alpha},\\
 \nonumber R^1&=&\sum_{j\geq 1} \left(R^{q_j}(x;\xi)q_j+
R^{\bar{q}_j}(x;\xi)\bar{q}_j\right)\\&=&\sum_{\alpha\in\mathbb{N}^n,\beta,\gamma\in\mathbb{N}^{\mathbb{N}},|\alpha|=0,
|\beta|+|\gamma|=1 }R^{\alpha\beta\gamma}(x;\xi)q^{\beta}\bar
q^{\gamma},\\
\nonumber R^{2}&=&\sum_{i,j\geq1} \left(R^{q_iq_j}(x;\xi)q_iq_j+
R^{q_i\bar{q}_j}(x;\xi)q_i\bar{q}_j+
R^{\bar{q}_i\bar{q}_j}(x;\xi)\bar{q}_i\bar{q}_j\right)\\
&=&\sum_{\alpha\in\mathbb{N}^n,\beta,\gamma\in\mathbb{N}^{\mathbb{N}},|\alpha|=0,
|\beta|+|\gamma|=2 }R^{\alpha\beta\gamma}(x;\xi)q^{\beta}\bar
q^{\gamma}.\label{102407}
\end{eqnarray}Moreover, write
$ F=F^{x}+F^{y}+F^{1}+F^{2} $ and $R^{high}=\sum_{j=0}^4R^{(j)}$,
where
\begin{eqnarray*}
R^{(0)}&=&\sum_{\alpha\in\mathbb{N}^n,\beta,\gamma\in\mathbb{N}^{\mathbb{N}},|\alpha|=2,
|\beta|+|\gamma|=0 }R^{\alpha\beta\gamma}(x;\xi)y^{\alpha},\\
R^{(1)}&=&\sum_{\alpha\in\mathbb{N}^n,\beta,\gamma\in\mathbb{N}^{\mathbb{N}},|\alpha|=1,
|\beta|+|\gamma|=1 }R^{\alpha\beta\gamma}(x;\xi)y^{\alpha}q^{\beta}
\bar{q}^{\gamma},\\
R^{(2)}&=&\sum_{\alpha\in\mathbb{N}^n,\beta,\gamma\in\mathbb{N}^{\mathbb{N}},|\alpha|=1
,|\beta|+|\gamma|=2 }R^{\alpha\beta\gamma}(x;\xi)
y^{\alpha}q^{\beta}\bar{q}^{\gamma},\\
R^{(3)}&=&\sum_{\alpha\in\mathbb{N}^n,\beta,\gamma\in\mathbb{N}^{\mathbb{N}},|\alpha|=0,
|\beta|+|\gamma|=3 }R^{\alpha\beta\gamma}(x;\xi)
q^{\beta}\bar{q}^{\gamma},\\
R^{(4)}&=&\sum_{\alpha\in\mathbb{N}^n,\beta,\gamma\in\mathbb{N}^{\mathbb{N}},2|\alpha|+|\beta|+|\gamma|\geq5\
\mbox{or}\ |\beta|+|\gamma|\geq 4}R^{\alpha\beta\gamma}(x;\xi)
y^{\alpha}q^{\beta}\bar{q}^{\gamma}\nonumber .\end{eqnarray*} By a
direct calculation,
\begin{eqnarray*}
\{R^{high},F\}^{low}&=&\sum_{j=1}^n R_{y_j}^{(0)}
F^{x}_{x_j}+\sum_{j=1}^n R_{y_j}^{(1)}F^{x}_{x_j}+\sum_{j=1}^n
R^{(1)}_{y_j}F^1_{x_j}+\sum_{j=1}^n R^{(2)}_{y_j}
F^{x}_{x_j}\\&&+\sqrt{-1}\sum_{j\geq 1} \left(R_{\bar
q_j}^{(1)}F^1_{q_j}-R_{q_j}^{(1)}
F^1_{\bar{q}_j}+R_{\bar{q}_j}^{(3)} F^1_{q_j}- R_{q_j}^{(3)}
F^1_{\bar{q}_j}\right).
\end{eqnarray*}More precisely,
\begin{eqnarray*}
\{R^{high},F\}^{low}=\{R^{high},F\}^y+\{R^{high},F\}^1+\{R^{high},F\}^2,
\end{eqnarray*}where
\begin{eqnarray*}
\{R^{high},F\}^y&=&\sum_{j=1}^n R_{y_j}^{(0)}
F^{x}_{x_j}+\sqrt{-1}\sum_{j\geq 1} \left(R_{\bar
q_j}^{(1)}F^1_{q_j}-R_{q_j}^{(1)} F^1_{\bar{q}_j}\right),\\
\{R^{high},F\}^1&=&\sum_{j=1}^n R_{y_j}^{(1)}F^{x}_{x_j},\\
\{R^{high},F\}^2&=&\sum_{j=1}^n R^{(1)}_{y_j}F^1_{x_j}+\sum_{j=1}^n
R^{(2)}_{y_j} F^{x}_{x_j}+\sqrt{-1}\sum_{j\geq 1}
\left(R_{\bar{q}_j}^{(3)} F^1_{q_j}- R_{q_j}^{(3)}
F^1_{\bar{q}_j}\right).
\end{eqnarray*}
Let $\partial_{\omega}=\omega\cdot
\partial_x$ and $W=\{R^{high},F^x+
F^1\}.$ Then the homological equation (\ref{16}) decomposes into
\begin{eqnarray}
\label{2401}\partial_{\omega}F^{x}+\widehat{N^{x}}&=&R^{x},\\
(\partial_{\omega}+\sqrt{-1}\Omega_j)F^{q_j}&=&
\label{2403}R^{q_j}+\{R^{high},F^x\}^{q_j},\qquad j\geq1,\\
(\partial_{\omega}-\sqrt{-1}\Omega_j)F^{\bar{q}_j}
&=&\label{2402}R^{\bar{q}_j}+\{R^{high},F^x\}^{\bar{q}_j},\qquad j\geq1,\\
\partial_{\omega}F^{y_j}+\widehat{N^{y_j}}&=&
\label{2404}R^{y_j}+W^{y_j},\qquad 1\leq j\leq n,\\
(\partial_{\omega}+\sqrt{-1}\Omega_i+\sqrt{-1}\Omega_j)F^{q_iq_j}&=&
\label{2407}R^{q_iq_j}+W^{q_iq_j},\qquad i,j\geq1,\\
(\partial_{\omega}+\sqrt{-1}\Omega_i -\sqrt{-1}\Omega_j)F^{q_i\bar
q_j}+\delta_{ij}\widehat{N^{q_i\bar{q}_j}}
&=&\label{2405} R^{q_i\bar{q}_j}+W^{q_i\bar{q}_j},\qquad i,j\geq1,\\
(\partial_{\omega}-\sqrt{-1}\Omega_i -\sqrt{-1}\Omega_j)F^{\bar
q_i\bar q_j}
&=&\label{2406}R^{\bar{q}_i\bar{q}_j}+W^{\bar{q}_i\bar{q}_j},\qquad
i,j\geq1,\end{eqnarray} where
\begin{equation*}\label{070801}\widehat{N^{x}}(\xi)=\widehat{R^{x}}(0;\xi),\quad
\widehat{N^{y_j}}(\xi)=\widehat{R^{y_j}}(0;\xi)+\widehat{W^{y_j}}(0;\xi),\quad
\widehat{N^{q_i\bar{q}_j}}(\xi)=\widehat{R^{q_i\bar{q}_j}}(0;\xi)+
\widehat{W^{q_i\bar{q}_j}}(0;\xi). \end{equation*}\color{red} Hence, we can
solve the homological equations (\ref{2401})-(\ref{2406}) one by
one, i.e. the homological equation (\ref{16}) is solvable. \color{black}Moreover, set
\begin{equation}\label{081501}N_+=N+
\widehat{N^{x}}(\xi)+\sum_{j=1}^n\widehat{N^{y_j}}(\xi)y_j+
\sum_{j\geq 1}\widehat{N^{q_j\bar{q}_j}}(\xi)q_j\bar q_j,
\end{equation}
\begin{equation}
\omega_{+j}=\omega_j+\widehat{N^{y_j}}(\xi),\qquad 1\leq j\leq n,
\end{equation}
and \begin{equation}
\Omega_{+j}=\Omega_j+\widehat{N^{q_j\bar{q}_j}}(\xi),\qquad j\geq1.
\end{equation}
\subsubsection{The solution of homological equation (\ref{16})}
\begin{thm}\label{20}
Consider a perturbation of the integrable Hamiltonian
$$H(x,y,q,\bar q;\xi)=N(y,q,\bar q;\xi)+R(x,y,q,\bar q;\xi),$$ where
$$N(y,q,\bar q;\xi)=\sum_{j=1}^{n}\omega_j(\xi)y_j+\sum_{j\geq1}\Omega_j(\xi)q_j\bar
q_j$$ is a parameter dependent integrable Hamiltonian and
$$R(x,y,q,\bar q;\xi)=R^{low}(x,y,q,\bar q;\xi)+R^{high}(x,y,q,\bar q;\xi).$$  Suppose assumptions A and B are fulfilled for
$\omega(\xi)$ and $\Omega(\xi)$,
\begin{equation}\label{081910}
 |||X_{R^{low}}|||_{p,D(s,r,r)\times\Pi}^T\leq \varepsilon,
 \end{equation}
and
 \begin{equation}\label{081911}
 |||X_{R^{high}}|||_{p,D(s,r,r)\times\Pi}^T\leq 1,
\end{equation}for some $0<s,r\leq 1$. For some fixed constant $\tau>n+1$,
let\begin{equation}\label{102403}
\mathcal{R}_{kl}=\left\{\xi\in\Pi:\left|\langle
k,\omega(\xi)\rangle+\langle
l,\Omega(\xi)\rangle\right|\leq\frac{\eta}{(|
k|+1)^{\tau}}\right\},\qquad
k\in\mathbb{Z}^n,l\in\mathbb{Z}^{\mathbb{N}},
\end{equation}and
let
\begin{equation}\label{102404}
\tilde{\Pi}=\Pi\backslash\bigcup_{k\in\mathbb{Z}^n,|l|\leq2,
|k|+|l|\neq0}\mathcal{R}_{kl}.
\end{equation}
Then for each $\xi\in\tilde{\Pi}$, the homological equation
(\ref{16})
$$\{N,F\}+R^{low}+\{R^{high},F\}^{low}=N_+-N$$ has a solution $F(x,y,q,\bar q;\xi)$
with the estimates\begin{equation}\label{012003}
 |||X_F|||^T_{p,D(s-\sigma/2,r-\sigma/2,r-\sigma/2)\times\tilde\Pi}\lessdot \frac{\varepsilon}{\eta^{6}{\sigma^{6\tau+5}}},
 \end{equation}
 and
 \begin{equation}\label{012004}
|||X_{N_+-N}|||^T_{p,D(s-\sigma/2,r-\sigma/2,r-\sigma/2)\times\tilde\Pi}\lessdot
\frac{\varepsilon}{\eta^{4}{\sigma^{4\tau+4}}},
\end{equation}where $0<\sigma<\frac1{10}\min\{s,r\}$ and $a\lessdot b$ means there exists a constant $c>0$ depending on $n$ and $\tau$ such that $a\leq cb$. More precisely, we have
\begin{eqnarray}
|||X_{F^x}|||^T_{p,D(s-\tilde\sigma,r,r)\times\tilde\Pi}&\lessdot& \frac{\varepsilon}{\eta^2\sigma^{2\tau+1}}\label{012005},\\
|||X_{F^{1}}|||^T_{p,D(s-3\tilde\sigma,r-3\tilde\sigma,r-3\tilde\sigma)\times\tilde\Pi}&\lessdot&
\frac{\varepsilon}{\eta^4\sigma^{4\tau+3}}\label{012006},\\
|||X_{F^{y}}|||^T_{p,D(s-5\tilde\sigma,r-5\tilde\sigma,r-5\tilde\sigma)\times\tilde\Pi}&\lessdot&
\frac{\varepsilon}{\eta^6\sigma^{6\tau+5}}\label{012007},\\
|||X_{F^{2}}|||^T_{p,D(s-5\tilde\sigma,r-5\tilde\sigma,r-5\tilde\sigma)
\times\tilde\Pi}&\lessdot&\frac{\varepsilon}{\eta^6\sigma^{6\tau+5}},\label{012008}
\end{eqnarray}
where $\tilde\sigma=\sigma/10$.

Moreover, the new Hamiltonian $H_+(x,y,q,\bar q;\xi)$ has the
following form
\begin{equation*}H_{+}(x,y,q,\bar q;\xi):=H\circ X_F^t|_{t=1}=N_{+}(y,q,\bar q;\xi)
+R_{+}(x,y,q,\bar q;\xi),
\end{equation*}
where $$N_{+}(y,q,\bar q;\xi)=\sum_{j=1}^n \omega_{+j}(\xi)y_j
+\sum_{j\geq 1}\Omega_{+j}(\xi)q_j\bar q_j$$ and
\begin{eqnarray} &&\nonumber R_{+}(x,y,q,\bar q;\xi)\\&=&\nonumber R^{high}+\{R^{high},F\}^{high}+
\int^1_0(1-t)\{\{N,F\},F\}\circ X_F^tdt \\&&
+\int_0^1\{R^{low},F\}\circ
X_F^tdt+\int^1_0(1-t)\{\{R^{high},F\},F\}\circ X_F^tdt\label{082002},
\end{eqnarray}with the following estimates hold:\\
(1) for each $\xi\in\tilde\Pi$, the symplectic map
$\Phi=X_{F}^t|_{t=1}$ satisfies
\begin{equation*}
||\Phi-id||_{p,D(s-\sigma/2,r-\sigma/2,r-\sigma/2)}\lessdot
\frac{\varepsilon}{\eta^{6}{\sigma^{6\tau+5}}} ,\end{equation*}
and
\begin{equation*}
|||D\Phi-Id|||_{p,D(s-\sigma,r-\sigma,r-\sigma)}\lessdot
\frac{\varepsilon}{\eta^6\sigma^{6\tau+6}};
\end{equation*}
(2) the frequencies $\omega_{+}(\xi)$ and $\Omega_{+}(\xi)$ satisfy
\begin{equation*}
||\omega_{+}(\xi)-\omega(\xi)||+\sup_{j\geq
1}||\partial_{\xi_j}(\omega_{+}(\xi)-\omega(\xi))||\lessdot
\frac{\varepsilon}{\eta^{4}{\sigma^{4\tau+4}}},\end{equation*}and
\begin{equation*}
||\Omega_{+}(\xi)-\Omega(\xi)||+\sup_{j\geq
1}||\partial_{\xi_j}(\Omega_{+}(\xi)-\Omega(\xi))||\lessdot
\frac{\varepsilon}{\eta^{4}{\sigma^{4\tau+4}}};
\end{equation*}(3) the perturbation $R_{+}(x,y,q,\bar q;\xi)$ satisfies
\begin{equation*}
|||X_{R_{+}^{low}}|||^T_{p,D(s-\sigma,r-\sigma,r-\sigma)\times\tilde\Pi}\lessdot
\left(\frac{\varepsilon}{\eta^6\sigma^{6\tau+6}}\right)^2,
\end{equation*}and
\begin{eqnarray*}
|||X_{R_{+}^{high}}|||^T_{p,D(s-\sigma,r-\sigma,r-\sigma)\times\tilde\Pi}
\lessdot1+\frac{\varepsilon}{\eta^6\sigma^{6\tau+6}}+\left(\frac{\varepsilon}{\eta^6\sigma^{6\tau+6}}\right)^2;\end{eqnarray*}\\
(4) the measure of the subset $\tilde\Pi$ of $\Pi$ satisfies
\begin{equation}\label{091701}
Meas\ \tilde\Pi\geq (Meas\ \Pi)(1-O(\eta)).
\end{equation}
\end{thm}
\begin{proof}First of all, we will give two simple
estimates.\\
\textbf{No. 1.} In view of (\ref{102403}) and (\ref{102404}), for
each $\xi\in\tilde\Pi$ and each $k\in\mathbb{Z}^n,
l\in\mathbb{Z}^{\mathbb{N}}$ satisfying $|l|\leq 2,|k|+|l|\neq 0$,
we have
\begin{equation*}
|\langle k,\omega(\xi)\rangle+\langle l,\Omega(\xi)\rangle|\geq
\frac{\eta}{(|k|+1)^\tau}.
\end{equation*}
 Hence,
\begin{equation}\label{081901}
{|\langle k,\omega(\xi)\rangle+\langle
l,\Omega(\xi)\rangle|}^{-1}\leq \frac{(|k|+1)^\tau}{\eta}.
\end{equation}
 Moreover, for each $j\geq 1$,
\begin{eqnarray*}
\partial_{\xi_j}\left(\langle k,\omega(\xi)\rangle+\langle
l,\Omega(\xi)\rangle\right)^{-1} =\frac{\langle
k,\partial_{\xi_j}\omega(\xi)\rangle+\langle
l,\partial_{\xi_j}\Omega(\xi)\rangle}{(\langle
k,\omega(\xi)\rangle+\langle l,\Omega(\xi)\rangle)^2}.
\end{eqnarray*}
Then in view of twist conditions (\ref{081902}), the inequality
(\ref{081901}) and $|l|\leq 2$, we obtain
\begin{eqnarray}
&&\nonumber|\partial_{\xi_j}\left(\langle
k,\omega(\xi)\rangle+\langle
l,\Omega(\xi)\rangle\right)^{-1}|\\
&=&\nonumber\left|\frac{\langle
k,\partial_{\xi_j}\omega(\xi)\rangle+\langle
l,\partial_{\xi_j}\Omega(\xi)\rangle}{(\langle
k,\omega(\xi)\rangle+\langle
l,\Omega(\xi)\rangle)^2}\right|\\
\nonumber
&\leq&\frac{(|k|+1)^{2\tau}(|k|+2)}{\eta^2}\nonumber\\
&\leq&\frac{2(|k|+1)^{2\tau+1}}{\eta^2} \label{081903}.
\end{eqnarray}
\textbf{No. 2.} Recall the notations (\ref{102405})-(\ref{102407}),
and let $U^{0}=U^{x}+U^{y}$. Then in view of Definition
\ref{080204}, we have
\begin{equation}\label{082403}
|||X_{U^{i}}|||^{T}_{p,D(s,r,r)\times\Pi}\leq
|||X_{U^{low}}|||^{T}_{p,D(s,r,r)\times\Pi}\leq
|||X_{U}|||^{T}_{p,D(s,r,r)\times\Pi},\quad 0\leq i\leq 2,
\end{equation}
and
\begin{equation}\label{082901}
|||X_{U^{high}}|||^{T}_{p,D(s,r,r)\times\Pi}\leq
|||X_{U}|||^{T}_{p,D(s,r,r)\times\Pi}.
\end{equation} Moreover,
\begin{equation}\label{082401}
|||X_{U^{x}}|||^{T}_{p,D(s,r,r)\times\Pi}\leq
|||X_{U^{0}}|||^{T}_{p,D(s,r,r)\times\Pi}
\end{equation}
and
\begin{equation}\label{082402}
|||X_{U^{y}}|||^{T}_{p,D(s,r,r)\times\Pi}\leq
|||X_{U^{0}}|||^{T}_{p,D(s,r,r)\times\Pi}.
\end{equation}
The inequalities (\ref{082401}) and (\ref{082402}) can be proven by
following the proof of Theorem \ref{0005}.

Now we will prove the theorem by the following steps.

 \textbf{Step 1.
Estimate $|||X_F|||^T_{p,D(s-\sigma/2,r-\sigma/2,r-\sigma/2)\times\tilde\Pi}$ and $
|||X_{N_+-N}|||^T_{p,D(s-\sigma/2,r-\sigma/2,r-\sigma/2)\times\tilde\Pi}
$.}\\
\color{red}Based on (\ref{011703}) and in view of $F=F^x+F^y+F^1+F^2$, it is sufficient to prove the inequalities of (\ref{012005})-(\ref{012008}) hold. Since the proofs of (\ref{012005})-(\ref{012008}) are similar, we only show (\ref{012008}) is right below and assume (\ref{012005})-(\ref{012007}) is proven.
\color{black}
Without loss of generality, consider the homological equations (\ref{2405}),
which is
\begin{eqnarray*}
(\partial_{\omega}+\sqrt{-1}\Omega_i -\sqrt{-1}\Omega_j)F^{q_i\bar
q_j}+\delta_{ij}\widehat{N^{q_i\bar{q}_j}}
&=& R^{q_i\bar{q}_j}+W^{q_i\bar{q}_j},\qquad i,j\geq1.
\end{eqnarray*}By passing to Fourier coefficients,
\begin{eqnarray*}
\widehat{F^{q_i\bar{q}_j}}(k;\xi)&=&\frac{\widehat{R^{q_i\bar{q}_j}}(k;\xi)
+\widehat{W^{q_i\bar{q}_j}}(k;\xi)}{\sqrt{-1}(\langle
k,\omega(\xi)\rangle+\Omega_i(\xi)-\Omega_j(\xi))},\qquad
|k|+|i-j|\neq0,
\end{eqnarray*}and
\begin{equation}\label{081906}
\widehat{N^{q_i\bar{q}_j}}(\xi)=\widehat{R^{{q}_i\bar{q}_j}}(0;\xi)
+\widehat{W^{q_i\bar{q}_j}}(0;\xi).
\end{equation}
Then in view of the inequalities (\ref{081901}) and (\ref{081903}),
for each $\xi\in\tilde\Pi$ and $j'\geq 1$, we obtain
\begin{eqnarray*}
&&|\widehat{F^{q_i\bar{q}_j}}(k;\xi)|+|\widehat{\partial_{\xi_{j'}}F^{q_i\bar{q}_j}(k;\xi)}|\\
&\leq&\frac{3(|k|+1)^{2\tau+1}}{\eta^2}\left(|\widehat{R^{q_i\bar{q}_j}}
(k;\xi)+\widehat{
W^{q_i\bar{q}_j}}(k;\xi)|+|\partial_{\xi_{j'}}\widehat{R^{q_i\bar{q}_j}}
(k;\xi)+\partial_{\xi_{j'}}\widehat{W^{q_i\bar{q}_j}}(k;\xi)|\right).
\end{eqnarray*}
Similarly, we have
\begin{eqnarray*}
&&|\widehat{F^{q_i{q}_j}}(k;\xi)|+|\widehat{\partial_{\xi_{j'}}F^{q_i{q}_j}(k;\xi)}|\\
&\leq&\frac{3(|k|+1)^{2\tau+1}}{\eta^2}\left(|\widehat{R^{q_i{q}_j}}
(k;\xi)+\widehat{
W^{q_i{q}_j}}(k;\xi)|+|\partial_{\xi_{j'}}\widehat{R^{q_i{q}_j}}
(k;\xi)+\partial_{\xi_{j'}}\widehat{W^{q_i{q}_j}}(k;\xi)|\right),
\end{eqnarray*}and
\begin{eqnarray*}
&&|\widehat{F^{\bar q_i\bar{q}_j}}(k;\xi)|+|\widehat{\partial_{\xi_{j'}}F^{\bar q_i\bar{q}_j}(k;\xi)}|\\
&\leq&\frac{3(|k|+1)^{2\tau+1}}{\eta^2}\left(|\widehat{R^{\bar
q_i\bar{q}_j}} (k;\xi)+\widehat{ W^{\bar
q_i\bar{q}_j}}(k;\xi)|+|\partial_{\xi_{j'}}\widehat{R^{\bar
q_i\bar{q}_j}} (k;\xi)+\partial_{\xi_{j'}}\widehat{W^{\bar
q_i\bar{q}_j}}(k;\xi)|\right),
\end{eqnarray*}
while considering homological equations (\ref{2407}) and (\ref{2406}).
Hence, by the above three inequalities and Theorem
\ref{0005}, we get
\begin{eqnarray}
&&\nonumber|||X_{F^{2}}|||^T_{p,D(s-5\tilde\sigma,r-5\tilde\sigma,r-5\tilde\sigma)
\times\tilde\Pi}\\ \nonumber &\lessdot&\frac{
1}{\eta^2\sigma^{2\tau+1}}|||X_{R^{2}+
W^{2}}|||^T_{p,D(s-4\tilde\sigma,r-5\tilde\sigma,r-5\tilde\sigma)
\times\tilde\Pi}\\\nonumber&\leq&\frac{
1}{\eta^2\sigma^{2\tau+1}}\left(|||X_{R^{2}}|||^T_{p,D(s-4\tilde\sigma,r-5\tilde\sigma,r-5\tilde\sigma)
\times\tilde\Pi}+|||X_{W^2}|||^T_{p,D(s-4\tilde\sigma,r-5\tilde\sigma,r-5\tilde\sigma)
\times\tilde\Pi}\right)\\ \nonumber&\lessdot&\frac{
1}{\eta^2\sigma^{2\tau+1}}\left(|||X_{R^{low}}|||^T_{p,D(s,r,r)
\times\Pi}+|||X_{W}|||^T_{p,D(s-4\tilde\sigma,r-4\tilde\sigma,r-4\tilde\sigma)
\times\tilde\Pi}\right)\\&&\nonumber\mbox{(in view of the
inequalities (\ref{011708}) and (\ref{082401}))}\\ &\lessdot&
\frac{\varepsilon}{\eta^6\sigma^{6\tau+5}}, \label{071409}
\end{eqnarray}
where the last inequality is based on (\ref{081910}) and the following estimate
\begin{eqnarray}\label{070803}
|||X_{W}|||^T_{p,D(s-4\tilde\sigma,r-4\tilde\sigma,r-4\tilde\sigma)\times\tilde\Pi}
\lessdot\frac{\varepsilon}{\eta^4\sigma^{4\tau+4}}
\end{eqnarray} which is based on $W=\{R^{high},F^x+F^{1}\}$, and in view of
the inequalities (\ref{081911}), (\ref{012005}), (\ref{012006}) and
Theorem \ref{021102}.

The estimate of $
|||X_{N_+-N}|||^T_{p,D(s-\sigma/2,r-\sigma/2,r-\sigma/2)\times\tilde\Pi}
$ is trivial, so we omit the details.
\textbf{Step 2. Estimate
$||\Phi-id||_{p,D(s-\sigma/2,r-\sigma/2,r-\sigma/2)}$ and
$|||D\Phi-Id|||_{p,D(s-\sigma,r-\sigma,r-\sigma)}$.}

Note $\Phi=X_F^{t}|_{t=1}$. Based on Theorem \ref{081503} and the
inequality (\ref{012003}), we have
\begin{equation*}
||\Phi-id||_{p,D(s-\sigma/2,r-\sigma/2,r-\sigma/2)}
\lessdot\frac{\varepsilon}{\eta^{6}{\sigma^{6\tau+5}}}.\end{equation*}
Moreover, by Cauchy inequality in Lemma \ref{081907}, we
obtain
\begin{equation*}
|||D\Phi-Id|||_{p,D(s-\sigma,r-\sigma,r-\sigma)}\lessdot
\frac{\varepsilon}{\eta^6\sigma^{6\tau+6}}.
\end{equation*}
\textbf{Step 3. Estimate the norm of the frequency
shift\begin{equation*} ||\omega_{+}(\xi)-\omega(\xi)||+\sup_{j\geq
1}||\partial_{\xi_j}(\omega_{+}(\xi)-\omega(\xi))||,\end{equation*}and
\begin{equation*}
||\Omega_{+}(\xi)-\Omega(\xi)||+\sup_{j\geq
1}||\partial_{\xi_j}(\Omega_{+}(\xi)-\Omega(\xi))||.
\end{equation*}}
Based on Theorem \ref{012002} and the estimate (\ref{012004}), it is easy to verify
\begin{equation*}
||\omega_{+}(\xi)-\omega(\xi)||+\sup_{j\geq
1}||\partial_{\xi_j}(\omega_{+}(\xi)-\omega(\xi))||\lessdot
\frac{\varepsilon}{\eta^{4}{\sigma^{4\tau+4}}},\end{equation*}and
\begin{equation*}
||\Omega_{+}(\xi)-\Omega(\xi)||+\sup_{j\geq
1}||\partial_{\xi_j}(\Omega_{+}(\xi)-\Omega(\xi))||\lessdot
\frac{\varepsilon}{\eta^{4}{\sigma^{4\tau+4}}}.
\end{equation*}

\textbf{Step 4. Estimate
$|||X_{R_{+}^{low}}|||^T_{p,D(s-\sigma,r-\sigma,r-\sigma)\times\tilde\Pi}$
and
$|||X_{R_{+}^{high}}|||^T_{p,D(s-\sigma,r-\sigma,r-\sigma)\times\tilde\Pi}$.}

Firstly,
\begin{eqnarray}
&&\nonumber\max\left\{|||X_{\{R^{high},F\}^{low}}|||^T_{p,D(s-\sigma,r-\sigma,r-\sigma)\times\tilde\Pi},|||X_{\{R^{high},F\}^{high}}|||^T_{p,D(s-\sigma,r-\sigma,r-\sigma)\times\tilde\Pi}\right\}\\
&\leq&\nonumber|||X_{\{R^{high},F\}}|||^T_{p,D(s-\sigma,r-\sigma,r-\sigma)\times\tilde\Pi}\qquad \mbox{(in view of the inequalities (\ref{082403}) and (\ref{082901}))}\\
&\lessdot&\nonumber\max\left\{\frac{1}{\sigma},\frac{r}{\sigma}\right\}|||X_{R^{high}}|||^T_{p,D(s-\sigma/2,r-\sigma/2,r-\sigma/2)\times\tilde\Pi}
|||X_{F}|||^T_{p,D(s-\sigma/2,r-\sigma/2,r-\sigma/2)\times\tilde\Pi}\\
&&\nonumber(\mbox{based on Theorem
\ref{021102}})\\&\lessdot&\nonumber\frac{1}{\sigma}|||X_{R^{high}}|||^T_{p,D(s-\sigma/2,r-\sigma/2,r-\sigma/2)\times\tilde\Pi}
|||X_{F}|||^T_{p,D(s-\sigma/2,r-\sigma/2,r-\sigma/2)\times\tilde\Pi}\\&&
\nonumber\mbox{(in view of $0<\sigma<\frac1{10}\min\{s,{r}\}$ and
$0<s,r\leq
{1}$)}\\&\lessdot&\frac{\varepsilon}{\eta^6\sigma^{6\tau+6}}\label{082001}\qquad
(\mbox{in view of the inequalities (\ref{081911}) and
(\ref{012003})}).
\end{eqnarray}
Secondly, estimate
\begin{equation*}
|||X_{\int^1_0(1-t)\{\{N,F\},F\}\circ
X_F^tdt}|||^T_{p,D(s-\sigma,r-\sigma,r-\sigma)\times\tilde\Pi}.
\end{equation*}
Note \begin{equation*} \{N,F\}=N_+-N-R^{low}+\{R^{high},F\}^{low}.
\end{equation*}
Then in view of the inequalities (\ref{081910}), (\ref{012004}) and
(\ref{082001}), we obtain
\begin{eqnarray}\label{081913}
|||X_{\{N,F\}}|||^T_{p,D(s-\frac34\sigma,r-\frac34\sigma,r-\frac34\sigma)
\times\tilde\Pi}\lessdot\frac{\varepsilon}{\eta^6\sigma^{6\tau+6}},
.\end{eqnarray}In view of the inequalities (\ref{012003}),
(\ref{081913}) and  Theorem \ref{021102}, we have
\begin{eqnarray}
&&\nonumber|||X_{\{\{N,F\},F\}}|||^T_{p,D(s-\frac78\sigma,r-\frac78\sigma,r-\frac78\sigma)\times\tilde\Pi}\\
&\lessdot&\nonumber\frac{1}{\sigma}|||X_{\{N,F\}}|||^T_{p,D(s-\frac34\sigma,r-\frac34\sigma,r-\frac34\sigma)\times\tilde\Pi}|||X_{F}|||^T_{p,D(s-\frac34\sigma,r-\frac34\sigma,r-\frac34\sigma)\times\tilde\Pi}\\
&\lessdot&\label{081915}\left(\frac{\varepsilon}{\eta^6\sigma^{6\tau+6}}\right)^2.\end{eqnarray}
By Theorem \ref{081914}, for $|t|\leq 1$, we get
\begin{equation}\label{081916}
|||X_{\{\{N,F\},F\}\circ
X_{F}^t}|||^T_{p,D(s-\sigma,r-\sigma,r-\sigma)\times\tilde\Pi}\lessdot|||X_{\{\{N,F\},F\}}|||^T_{p,D(s-\frac78\sigma,r-\frac78\sigma,r-\frac78\sigma)\times\tilde\Pi}.
\end{equation}
Therefore, by (\ref{081915}) and
(\ref{081916}),
\begin{equation}\label{081917} |||X_{\int^1_0(1-t)\{\{N,F\},F\}\circ
X_F^tdt}|||^T_{p,D(s-\sigma,r-\sigma,r-\sigma)\times\tilde\Pi}
\lessdot\left(\frac{\varepsilon}{\eta^6\sigma^{6\tau+6}}\right)^2.\end{equation}
In view of the inequalities (\ref{081910}), (\ref{081911}) and
(\ref{012003}) and  following the proof of the inequality
(\ref{081917}), we obtain
\begin{equation}\label{081918} |||X_{\int^1_0\{R^{low},F\}\circ
X_F^tdt}|||^T_{p,D(s-\sigma,r-\sigma,r-\sigma)\times\tilde\Pi}
\lessdot\varepsilon\left(\frac{\varepsilon}{\eta^6\sigma^{6\tau+6}}\right)\end{equation}
and
\begin{equation}\label{081919} |||X_{\int^1_0(1-t)\{\{R^{high},F\},F\}\circ
X_F^tdt}|||^T_{p,D(s-\sigma,r-\sigma,r-\sigma)\times\tilde\Pi}
\lessdot\left(\frac{\varepsilon}{\eta^6\sigma^{6\tau+6}}\right)^2.\end{equation}
In view of the formula (\ref{082002}), then
\begin{eqnarray*}R_{+}^{low}\nonumber&=&\left(\int^1_0(1-t)\{\{N,F\},F\}\circ
X_F^tdt\right)^{low}
\\&&+\left(\int_0^1\{R^{low},F\}\circ X_F^tdt\right)^{low}+
\left(\int^1_0(1-t)\{\{R^{high},F\},F\}\circ X_F^tdt\right)^{low},
\end{eqnarray*}
and
\begin{eqnarray*}R_{+}^{high}\nonumber&=&R^{high}
+\{R^{high},F\}^{high}+\left(\int^1_0(1-t)\{\{N,F\},F\}\circ
X_F^tdt\right)^{high}
\\&&+\left(\int_0^1\{R^{low},F\}\circ X_F^tdt\right)^{high}+
\left(\int^1_0(1-t)\{\{R^{high},F\},F\}\circ X_F^tdt\right)^{high},
\end{eqnarray*}
Hence, in view of the inequalities (\ref{081911}), (\ref{082001}),
(\ref{081917}), (\ref{081918}) and (\ref{081919}), we obtain
\begin{equation*}
|||X_{R_{+}^{low}}|||^T_{p,D(s-\sigma,r-\sigma,r-\sigma)\times\tilde\Pi}\lessdot
\left(\frac{\varepsilon}{\eta^6\sigma^{6\tau+6}}\right)^2,
\end{equation*}and
\begin{eqnarray*}
|||X_{R_{+}^{high}}|||_{p,D(s-\sigma,r-\sigma,r-\sigma)\times\tilde\Pi}
\lessdot1+\frac{\varepsilon}{\eta^6\sigma^{6\tau+6}}+\left(\frac{\varepsilon}{\eta^6\sigma^{6\tau+6}}\right)^2.
\end{eqnarray*}

The proof of the inequality (\ref{091701}) (the estimate of the
measure of $\tilde\Pi$) is omitted, since the inequality
(\ref{091701}) can be proven by following the proof of the measure
of $\Pi_{\tilde\eta}$ in Section \ref{091406}.
\end{proof}
\subsection{Iterative lemma}
\subsubsection{Iterative constants}As usual, the KAM theorem is
proven by the Newton-type iteration procedure which involves an
infinite sequence of coordinate changes. In order to make our
iteration procedure run, we need the following
iteration constants:\\
1. $\eta_0=\eta$ is given, $\eta_{m}={\eta}2^{-m}$, $m=1,2,\dots$;\\
2.
$\varepsilon_0=\varepsilon=\eta^{12}\epsilon,\varepsilon_{m}=\eta^{12}\epsilon^{(4/3)^{m}}
,{m}=1,2,\dots;$\\
3.$\
\tau_0=0,\tau_{m}=(1^{-2}+\cdots+{m}^{-2})/2\sum_{j=1}^{\infty}j^{-2}$,
 ${m}=1,2,\dots$ (thus, $\tau_{m}<1/2$);\\
4. Given $0<s_0,r_0\leq1$. Let $0<\sigma\leq \min\{s_0,r_0\}$ and
$s_{m}=(1-\tau_{m}) \sigma$, $r_{m}=(1-\tau_{m})\sigma$,
$m=1,2,\dots$ (thus, $s_{m}>s_0/2, r_{m}>r_0/2$).

\subsubsection{Iterative lemma}\begin{lem}Consider a perturbation of the integrable Hamiltonian
$$ H_m(x,y,q,\bar q;\xi)=N_m(y,q,\bar q;\xi)+R_m(x,y,q,\bar
q;\xi),$$ where
\begin{equation*}
N_{m}(y,q,\bar q;\xi)=\sum_{j=1}^n \omega_{mj}(\xi) y_j+\sum_{j\geq
1} \Omega_{mj}(\xi) q_j\bar{q}_j\end{equation*} is a parameter
dependent integrable Hamiltonian and $$R_m(x,y,q,\bar
q;\xi)=R_m^{low}(x,y,q,\bar q;\xi)+R_m^{high}(x,y,q,\bar q;\xi)$$ is
the perturbation with the following form
\begin{eqnarray*}R_m^{low}(x,y,q,\bar q;\xi)=\sum_{\alpha\in\mathbb{N}^n,\beta,\gamma\in\mathbb{N}^{\mathbb{N}},2|\alpha|+|\beta|+|\gamma|\leq 2}R_m^{\alpha\beta\gamma}(x;\xi)y^{\alpha}q^{\beta}\bar
q^{\gamma}
\end{eqnarray*}and
\begin{eqnarray*}
R_m^{high}(x,y,q,\bar
q;\xi)=\sum_{\alpha\in\mathbb{N}^n,\beta,\gamma\in\mathbb{N}^{\mathbb{N}},2|\alpha|+|\beta|+|\gamma|\geq
3}R_m^{\alpha\beta\gamma}(x;\xi)y^{\alpha}q^{\beta}\bar q^{\gamma}.
\end{eqnarray*} Suppose the assumptions A and B are fulfilled for
$\omega_m(\xi)$ and $\Omega_m(\xi)$ with $m=0$ and
\begin{equation*}
||\omega_m(\xi)-\omega_0(\xi)||+\sup_{j\geq
1}||\partial_{\xi_j}(\omega_m(\xi)-\omega_0(\xi))||\leq
\sum_{i=1}^{m}\varepsilon^{2/3}_{i-1},\end{equation*}and
\begin{equation*}
||\Omega_m(\xi)-\Omega_0(\xi)||+\sup_{j\geq
1}||\partial_{\xi_j}(\Omega_m(\xi)-\Omega_0(\xi))||\leq
\sum_{i=1}^{m}\varepsilon^{2/3}_{i-1}.
\end{equation*}
Suppose
$R^{low}_m(x,y,q,\bar q;\xi)$ satisfies the smallness assumption
\begin{equation*}
 |||X_{R^{low}_m}|||^T_{p,D(s_m,r_m,r_m)\times\Pi_m}\leq \varepsilon_m,
\end{equation*}and $R_m^{high}(x,y,q,\bar q;\xi)$ satisfies
\begin{equation*}
 |||X_{R^{high}_m}|||^T_{p,D(s_m,r_m,r_m)\times\Pi_m}\leq \varepsilon+\sum_{i
 =1}^m \varepsilon_i^{2/3}.\end{equation*}
Let\begin{equation*}
\mathcal{R}_{kl}^m=\left\{\xi\in\Pi:\left|\langle
k,\omega_m(\xi)\rangle+\langle
l,\Omega_m(\xi)\rangle\right|\leq\frac{\eta_m}{(|
k|+1)^{\tau}}\right\},\qquad k\in
\mathbb{Z}^n,l\in\mathbb{Z}^{\mathbb{N}},
\end{equation*}and let
\begin{equation*}
{\Pi}_{m+1}=\Pi_m\backslash\bigcup_{k\in\mathbb{Z}^n,|l|\leq2,|k|+|l|\neq0}\mathcal{R}_{kl}^m.
\end{equation*}
Then for each $\xi\in{\Pi}_{m+1}$, the homological equation
$$\{N_m,F_m\}+R^{low}_m+\{R^{high}_m,F_m\}^{low}=N_{+m}-N_m$$ has a solution $F_m(x,y,q,\bar q;\xi)$ with the estimates
\begin{equation*}
 |||X_{F_m}|||^T_{p,D(s_{m+1},r_{m+1},r_{m+1})\times\Pi_{m+1}}\leq \varepsilon_m^{2/3},
 \end{equation*}and
 \begin{equation*}\label{071406}
|||X_{N_{+m}-N_m}|||^T_{p,D(s_{m+1},r_{m+1},r_{m+1})\times\Pi_{m+1}}\leq
\varepsilon^{2/3}_m.
\end{equation*}
Moreover,
\begin{equation*}H_{m+1}(x,y,q,\bar q;\xi):=H\circ X_{F_m}^t|_{t=1}=N_{m+1}(y,q,\bar q;\xi)+R_{m+1}(x,y,q,\bar q;\xi),
\end{equation*}where
\begin{equation*}
N_{m+1}(y,q,\bar q;\xi)=N_{+m}(y,q,\bar q;\xi)=\sum_{j=1}^n
\omega_{(m+1)j}(\xi)y_j+\sum_{j\geq1} \Omega_{(m+1)j}(\xi)
q_j\bar{q}_j,
\end{equation*}and
\begin{eqnarray*} &&R_{m+1}(x,y,q,\bar q;\xi)\\&=&R^{high}_m+\{R^{high}_m,F_m\}^{high}+ \int^1_0(1-t)\{\{N,F_{m}\},F_m\}\circ X_F^tdt
\\&&
+\int_0^1\{R^{low}_m,F_m\}\circ
X_{F_m}^tdt+\int^1_0(1-t)\{\{R_m^{high},F_m\},F_m\}\circ
X_{F_m}^tdt,
\end{eqnarray*}with the following estimates hold:\\
(1) for each $\xi\in\Pi_{m+1}$, the symplectic map
$\Psi_m=X_{F_m}^t|_{t=1}$ satisfies
\begin{equation*}
||\Psi_m-id||_{p,D(s_{m+1},r_{m+1},r_{m+1})}\leq \varepsilon^{2/3}_m
,\end{equation*} and
\begin{equation*}
|||D\Psi_m-Id|||_{p,D(s_{m+1},r_{m+1},r_{m+1})}\leq
\varepsilon^{2/3}_m.
\end{equation*}
\\
 (2) the frequencies $\omega_{m+1}(\xi)$ and $\Omega_{m+1}(\xi)$
satisfy
\begin{equation*}
||\omega_{m+1}(\xi)-\omega_0(\xi)||+\sup_{j\geq1}||\partial_{\xi_j}(\omega_{m+1}(\xi)-\omega_0(\xi))||\leq
\sum_{i=1}^{m+1}\varepsilon^{2/3}_{i-1},\end{equation*}and
\begin{equation*}
||\Omega_{m+1}(\xi)-\Omega_0(\xi)||+\sup_{j\geq1}||\partial_{\xi_j}(\Omega_{m+1}(\xi)-\Omega_0(\xi))||\leq
\sum_{i=1}^{m+1}\varepsilon^{2/3}_{i-1};
\end{equation*}(3) the perturbation $R_{m+1}(x,y,q,\bar q;\xi)$ satisfies
\begin{equation*}
 |||X_{R^{low}_{m+1}}|||^T_{p,D(s_{m+1},r_{m+1},r_{m+1})\times\Pi_{m+1}}\leq
 \varepsilon_{m+1},\end{equation*}and
 \begin{equation*}
|||X_{R_{m+1}^{high}}|||^T_{p,D(s_{m+1},r_{m+1},r_{m+1})\times\Pi_{m+1}}\leq
\varepsilon+\sum_{i=1}^{m+1} \varepsilon_i^{2/3} ;\end{equation*}
(4) the measure of the subset $\Pi_{m+1}$ of $\Pi_{m}$ satisfies
\begin{equation}
Meas\ \Pi_{m+1}\geq (Meas\ \Pi_{m})(1-O(\eta_m)).
\end{equation}
\end{lem}
\begin{proof}The proof is a standard KAM proof based on Theorem \ref{021102}, Theorem \ref{081914} and Theorem \ref{0005}. See details of KAM iteration in \cite{P1}.
\end{proof}
\subsection{Proof of Theorem \ref{T1}}
\begin{proof}Let $\Pi_{\eta}=\bigcap_{m=0}^{\infty}\Pi_m,
D(s_0/2,r_0/2,r_0/2)\subset\bigcap_{m=0}^{\infty}D(s_m,r_m,r_m)$ and
$\Psi=\prod_{m=0}^{\infty}\Psi_m$. By the standard argument, we
conclude that $\Psi,D\Psi,H_m,X_{H_{m}}$ converge uniformly on the
domain $D(s_0/2,r_0/2,r_0/2)\times\Pi_{\eta}$. Let
\begin{equation*}
\breve H(x,y,q,\bar q;\xi):=\lim_{m\rightarrow\infty}H_m=\breve
N(y,q,\bar q;\xi)+\breve R(x,y,q,\bar q;\xi),
\end{equation*}
where
$$\breve N(y,q,\bar q;\xi)=\sum_{j=1}^n\breve\omega_j(\xi)y_j+\sum_{j\geq1} \breve
\Omega_j(\xi) q_j\bar q_j$$ and $$\breve R(x,y,q,\bar
q;\xi)=\sum_{\alpha\in\mathbb{N}^n,\beta,\gamma\in\mathbb{N}^{\mathbb{N}},2|\alpha|+|\beta|+|\gamma|\geq
3}\breve R^{\alpha\beta\gamma}(x;\xi)y^{\alpha}q^{\beta}\bar
q^{\gamma}.$$ Moreover, by the standard KAM proof, we obtain the following estimates:\\
(1) for each $\xi\in\Pi_{\eta}$, the symplectic map $\Psi$ satisfies
\begin{equation*}
||\Psi-id||_{p,D(s_0/2,r_0/2,r_0/2)}\leq c\eta^6\epsilon,
\end{equation*}
and
\begin{equation*}
|||D{\Psi}-Id|||_{p,D(s_0/2,r_0/2,r_0/2)}\leq {c\eta^6\epsilon}.
\end{equation*}
 \\
(2) the frequencies $\breve \omega(\xi)$ and $\breve \Omega(\xi)$
satisfy
\begin{equation*}
||\breve\omega(\xi)-\omega(\xi)||+\sup_{j\geq1}||\partial_{\xi_j}(\breve\omega(\xi)-\omega(\xi))||\leq
{c\eta^8\epsilon},\end{equation*}and
\begin{equation*}
||\breve\Omega(\xi)-\Omega(\xi)||+\sup_{j\geq1}||\partial_{\xi_j}(\breve\Omega(\xi)-\Omega(\xi))||\leq
{c\eta^8\epsilon};
\end{equation*}
(3) the Hamiltonian vector field $X_{\breve R}$ satisfies
\begin{equation*}
 |||X_{\breve R}|||^T_{p,D(s_0/2,r_0/2,r_0/2)\times
 \Pi_{\eta}}\leq \varepsilon(1+{c\eta^6\epsilon}),\end{equation*}
 where $c>0$ is a constant depending on $s_0,r_0,n$ and $\tau$.
 Moreover, we fix $\tau> n+1$ and then the constant $c$ depends
 on $s_0,r_0$ and $n$;\\
 (4) the measure of
$\Pi_{\eta}$ satisfies \begin{equation*} \mbox{Meas}\
\Pi_{\eta}\geq(\mbox{Meas}\ \Pi)(1-O(\eta)).
\end{equation*}
\end{proof}
\subsection{Proof of Corollary \ref{L1}}\begin{proof}In view of (\ref{081601}),
\begin{equation*}\breve H(x,y,q,\bar q;\xi)=
\breve N(y,q,\bar q;\xi)+ \breve R(x,y,q,\bar q;\xi),
\end{equation*}where
\begin{equation*}
\breve N(y,q,\bar
q;\xi)=\sum_{j=1}^{n}\breve{\omega}_j(\xi)y_j+\sum_{j\geq1}\breve{\Omega}_j(\xi)q_j\bar
q_j
\end{equation*}
and
\begin{equation*}
\breve R(x,y,q,\bar
q;\xi)=\sum_{\alpha\in\mathbb{N}^n,\beta,\gamma\in\mathbb{N}^{\mathbb{N}},2|\alpha|+|\beta|+|\gamma|\geq
3}\breve
R^{\alpha\beta\gamma}(x;\xi)y^{\alpha}q^{\beta}\bar{q}^{\gamma}.
\end{equation*}It is easy to verify that
$$\mathcal{T}_{0}=\hat{\mathbb{T}}^n\times\{y=0\}\times\{q=0\}\times\{\bar
q=0\}$$ is an embedding torus with frequency
$\breve\omega(\xi)\in\breve\omega(\Pi_{\eta})$ of the Hamiltonian
${\breve H}(x,y,q,\bar q;\xi)$. Moreover,
 $\Psi^{-1}\mathcal{T}_{0}$ is an embedding torus of the original Hamiltonian $H(x,y,q,\bar q;\xi)$.
We finish the proof of the existence of KAM tori. The proof of time
$|t|\leq \delta^{-1}$ stability of KAM tori is omitted here, since it is
similar to the proof of $|t|\leq \delta^{-\mathcal{M}}$ stability of
KAM tori, which can be found in Section \ref{091702}.
\end{proof}
\section{Proof of Theorem \ref{T3}}
\subsection{Construct a partial normal form of order $\mathcal{M}+2$}
Basing on the normal form of order 2 obtained in Theorem \ref{T1}, we
will construct a partial normal form of order $\mathcal{M}+2$ in the
neighborhood of KAM tori by $\mathcal{M}$ times symplectic transformations.

Given a large $\mathcal{N}\in\mathbb{N}$, split the normal frequency
$\breve\Omega(\xi)$ and normal variable $(q,\bar q)$ into two parts
respectively, i.e.
$$ \breve\Omega(\xi)=(\tilde {\Omega}(\xi),\hat{\Omega}(\xi)),\qquad q=(\tilde
q,\hat q),\qquad \bar q=(\bar{\tilde q},\bar {\hat q}),$$ where
$$\tilde
{\Omega}(\xi)=(\breve\Omega_1(\xi),\dots,\breve\Omega_{\mathcal{N}}(\xi)),\quad
\tilde q=(q_1,\dots,q_{\mathcal{N}}),\quad \bar{\tilde q}=(\bar
q_{1},\dots,\bar q_{\mathcal{N}})$$ are the low frequencies  and
$$
\hat{\Omega}(\xi)=(\breve\Omega_{\mathcal{N}+1}(\xi),\breve\Omega_{\mathcal{N}+2}(\xi),\dots),\quad
\hat q=(q_{\mathcal{N}+1},q_{\mathcal{N}+2},\dots),\quad \bar{\hat
q}=(\bar q_{\mathcal{N}+1},\bar q_{\mathcal{N}+2},\dots)$$
 are the high frequencies. Given $0<\tilde \eta<1$,
and $\tau> n+3$, if the frequencies $\breve\omega(\xi)$ and
$\breve\Omega(\xi)$ satisfy the following inequalities
\begin{equation}\label{090302} \left|\langle
k,\breve\omega(\xi)\rangle+\langle \tilde
l,\tilde{\Omega}(\xi)\rangle+\langle
\hat{l},\hat{\Omega}(\xi)\rangle\right|\geq\frac{\tilde\eta}{4^{\mathcal{M}}(|
k|+1)^{\tau}C(\mathcal{N},\tilde l)},
\end{equation}for any $k\in\mathbb{Z}^n,\tilde
l\in\mathbb{Z}^{\mathcal{N}}, \hat l\in\mathbb{Z}^{\mathbb{N}}$
with $$|k|+|\tilde l|+|\hat l|\neq 0,\quad |\tilde l|+|\hat
l|\leq \mathcal{M}+2,\quad |\hat l|\leq 2,$$ where
\begin{equation}\label{121206}
C(\mathcal{N},\tilde l)=\mathcal{N}^{(|\tilde l|+4)^2},
\end{equation}
then we call that the
frequencies $\breve\omega(\xi)$ and $\breve\Omega(\xi)$ are
$(\tilde\eta,\mathcal{N},\mathcal{M})$-non-resonant. Define the
resonant sets $\mathcal{R}_{k\tilde l\hat l}$ by
\begin{equation}\label{091202}
\mathcal{R}_{k\tilde{l}\hat{l}}=\left\{\xi\in\Pi_{\eta}:\left|\langle
k,\breve\omega(\xi)\rangle+\langle \tilde
l,\tilde{\Omega}(\xi)\rangle+\langle
\hat{l},\hat{\Omega}(\xi)\rangle\right|\leq\frac{\tilde\eta}{4^{\mathcal{M}}(|
k|+1)^{\tau}C(\mathcal{N},\tilde l)}\right\}.
\end{equation}
Let
\begin{equation}\label{091208}
\mathcal{R}=\bigcup_{|k|+|\tilde l|+|\hat l|\neq 0, |\tilde l|+|\hat
l|\leq \mathcal{M}+2, |\hat l|\leq
2}\mathcal{R}_{k\tilde{l}\hat{l}},
\end{equation}and
\begin{equation}\label{091215}
\Pi_{\tilde\eta}=\Pi_{\eta}\setminus\mathcal{R},
\end{equation}
where $\Pi_{\eta}$ is defined in Theorem \ref{T1}. Then it is easy
to see that for each $\xi\in\Pi_{\tilde\eta}$, the frequencies
$\breve\omega(\xi)$ and $\breve\Omega(\xi)$ are
$(\tilde\eta,\mathcal{N},\mathcal{M})$-non-resonant.

In this
section, we always assume
\begin{equation*}
\alpha\in\mathbb{N}, \quad
\beta,\gamma\in\mathbb{N}^{\mathcal{N}},\quad
\mu,\nu\in\mathbb{N}^{\mathbb{N}}.
\end{equation*}
\begin{thm}\label{thm7.1} (Partial normal form of order $\mathcal{M}+2$)
Consider the normal form of order 2 $$\breve H(x,y,q,\bar
q;\xi)=\breve N(y,q,\bar q;\xi)+\breve R(x,y,q,\bar q;\xi)$$
obtained in Theorem \ref{T1}. Given any positive integer
$\mathcal{M}$ and $0<\tilde \eta<1$, there exist a small $\rho_0>0$
and a large positive integer $\mathcal{N}_0$ depending on
$s_0,r_0,n$ and $\mathcal{M}$. For each $0<\rho<\rho_0$ and any
integer $\mathcal{N}$ satisfying
\begin{equation}\label{091113}\mathcal{N}_0<\mathcal{N}<\left(\frac{\tilde\eta^2}{2\rho}\right)^{\frac{1}{2(\mathcal{M}+7)^2}}, \end{equation}
and for each $\xi\in\Pi_{\tilde\eta}$, then there is a symplectic
map
$$\Phi: D(s_0/4,4\rho,4\rho)\rightarrow D(s_0/2,5\rho,5\rho),$$ such
that
\begin{equation}\label{020302}\breve{\breve{H}}(x,y,q,\bar q;\xi):=
\breve H\circ\Phi={\breve N}(y,q,\bar q;\xi)+{Z}(y,q,\bar q;\xi)+{
P}(x,y,q,\bar q;\xi)+{ Q}(x,y,q,\bar q;\xi)
\end{equation}
is a partial normal form of order $\mathcal{M}+2$, where
\begin{eqnarray*}
{ Z}(y,q,\bar q;\xi)&=&\sum_{4\leq 2|\alpha|+2|\beta|+2|\mu|\leq
\mathcal{M}+2,|\mu|\leq1}{
Z}^{\alpha\beta\beta\mu\mu}(\xi)y^{\alpha}\tilde{q}^{\beta}\bar{\tilde
q}^{\beta}\hat{q}^{\mu}\bar{\hat{q}}^{\mu}
\end{eqnarray*}is the integrable term depending only on variables $y$ and $I_j=|q_j|^2,j\geq1$, and where
\begin{eqnarray*}  {
P}(x,y,q,\bar
q;\xi)&=&\sum_{2|\alpha|+|\beta|+|\gamma|+|\mu|+|\nu|\geq
\mathcal{M}+3,|\mu|+|\nu|\leq 2} { P}^{\alpha
\beta\gamma\mu\nu}(x;\xi)y^{\alpha}{\tilde q}^{\beta}{\bar{\tilde
q}}^{\gamma}{\hat q}^{\mu}{\bar{\hat q}^{\nu}},
\end{eqnarray*}and
\begin{equation*}
{ Q}(x,y,q,\bar q;\xi)=\sum_{|\mu|+|\nu|\geq 3}{ Q}^{\alpha
\beta\gamma\mu\nu}(x;\xi)y^{\alpha}{\tilde q}^{\beta}{\bar{\tilde
q}}^{\gamma}{\hat q}^{\mu}{\bar{\hat q}}^{\nu}.
\end{equation*}
Moreover, the following estimates hold:\\
(1) the symplectic map $\Phi$ satisfies
\begin{equation}\label{091812}
||\Phi-id||_{p,D(s_0/4,4\rho,4\rho)}\leq \frac{c\mathcal{N}^{98}\rho}{\tilde\eta^2},
\end{equation}
and
\begin{equation}
|||D\Phi-Id|||_{p,D(s_0/4,4\rho,4\rho)}\leq \frac{c\mathcal{N}^{98}}{\tilde\eta^2};
\end{equation}
\\
(2) the Hamiltonian vector fields $X_Z,X_P$ and $X_Q$ satisfy
\begin{equation*}
|||X_{{ Z}}|||^T_{p,D(s_0/4,4\rho,4\rho)\times\Pi_{\tilde\eta}}\leq
c\rho
\left(\frac1{\tilde\eta^2}\mathcal{N}^{2(\mathcal{M}+6)^2}\rho\right),
\end{equation*}
\begin{equation*}
|||X_{{ P}}|||^T_{p,D(s_0/4,4\rho,4\rho)\times\Pi_{\tilde\eta}}\leq
c
\rho\left(\frac1{\tilde\eta^2}\mathcal{N}^{2(\mathcal{M}+7)^2}\rho\right)^{\mathcal{M}},\qquad
 \end{equation*}and
\begin{equation*}
|||X_{{ Q}}|||^T_{p,D(s_0/4,4\rho,4\rho)\times\Pi_{\tilde\eta}}\leq
c \rho,
\end{equation*}
where $c>0$ is a constant depending on $s_0,r_0,n$ and
$\mathcal{M}$.
\end{thm}

To prove Theorem \ref{thm7.1}, we will give an iterative lemma first.
Take
\begin{equation}\label{090301} s'=\frac{s_0}{12\mathcal{M}}\qquad
\mbox{and}\qquad \rho'=\frac{\rho}{2\mathcal{M}}.
\end{equation}
Let $2\leq j_0\leq\mathcal{M}+2$ and denote
\begin{equation}\label{091201}
\mathcal{D}_{j_0}=D(s_0/2-3(j_0-2)s',5\rho-2(j_0-2)\rho',5\rho-2(j_0-2)\rho'),
\end{equation}
\begin{equation*}
\mathcal{D}_{j_0}'=D(s_0/2-(3(j_0-2)+1)s',5\rho-2(j_0-2)\rho',
5\rho-2(j_0-2)\rho')
\end{equation*}
and
\begin{equation*}
\mathcal{D}_{j_0}''=D(s_0/2-(3(j_0-2)+2)s',5\rho-(2(j_0-2)+1)\rho',
5\rho-(2(j_0-2)+1)\rho')
\end{equation*}
 Then it is easy to see
 \begin{equation*}
 \mathcal{D}_{j_0+1}\subset\mathcal{D}_{j_0}''\subset\mathcal{D}_{j_0}'\subset\mathcal{D}_{j_0}
 ,\end{equation*}and
\begin{equation*}
\mathcal{D}_2=D(s_0/2,5\rho,5\rho)\qquad \mbox{and} \qquad
\mathcal{D}_{\mathcal{M}+2}=D(s_0/4,4\rho,4\rho).
\end{equation*}
\begin{lem}\label{052003}
Consider the partial normal form of order $j_0$ ($2\leq
j_0\leq\mathcal{M}+1 $) \begin{equation} H_{j_0}(x,y,q,\bar
q;\xi)=\breve N(y,q,\bar q;\xi)+Z_{j_0}(y,q,\bar
q;\xi)+P_{j_0}(x,y,q,\bar q;\xi) +Q_{j_0}(x,y,q,\bar
q;\xi),\end{equation}where
\begin{eqnarray}
\nonumber\breve N(y,q,\bar q;\xi)&=&\sum_{j=1}^n
\breve\omega_j(\xi)y_j+\sum_{j\geq1} \breve\Omega_j(\xi)
q_j\bar{q}_j,\\
\label{090602}Z_{j_0}(y,q,\bar q;\xi)&=&\sum_{3\leq j\leq
j_0}Z_{j_0j}(y,q,\bar
q;\xi),\\
P_{j_0}(x,y,q,\bar q;\xi)&=&\sum_{j\geq j_0+1}P_{j_0j}(x,y,q,\bar
q;\xi),\\
Q_{j_0}(x,y,q,\bar q;\xi)&=&\sum_{j\geq 3}Q_{j_0j}(x,y,q,\bar
q;\xi),
\end{eqnarray}
with
\begin{eqnarray}Z_{j_0j}(y,q,\bar q;\xi)&=&\sum_{2|\alpha|+2|\beta|+2|\mu|=j, |\mu|\leq1}
Z_{j_0}^{\alpha\beta\beta\mu\mu}(\xi)y^{\alpha}\tilde{q}
^{\beta}\bar{\tilde{q}}^{\beta}\hat{q}^{\mu}\bar{\hat{q}}^{\mu},\\
\label{090604}P_{j_0j}(x,y,q,\bar
q;\xi)&=&\sum_{2|\alpha|+|\beta|+|\gamma|+|\mu|+|\nu|= j,
|\mu|+|\nu|\leq 2}
P_{j_0}^{\alpha\beta\gamma\mu\nu}(x;\xi)y^{\alpha}\tilde{q}^{\beta}
\bar{\tilde{q}}^{\gamma}\hat{q}^{\mu}\bar{\hat{q}}^{\nu},\\
Q_{j_0j}(x,y,q,\bar
q;\xi)&=&\sum_{2|\alpha|+|\beta|+|\gamma|+|\mu|+|\nu|=
j,|\mu|+|\nu|\geq3}
Q_{j_0}^{\alpha\beta\gamma\mu\nu}(x;\xi)y^{\alpha}\tilde{q}
^{\beta}\bar{\tilde{q}}^{\gamma}\hat{q}^{\mu}\bar{\hat{q}}^{\nu}.\end{eqnarray}
Suppose $Z_{j_0j}(y,q,\bar q;\xi),P_{j_0j}(x,y,q,\bar q;\xi)$ and
$Q_{j_0j}(x,y,q,\bar q;\xi)$ satisfy the following estimates,
\begin{eqnarray}\label{083106}
|||X_{Z_{j_0j}}|||^T_{p,\mathcal{D}_{j_0}\times\Pi_{\tilde\eta}}&\preceq&
\rho\left(\frac1{\tilde\eta^2}\mathcal{N}^{2(j_0+4)^2}\rho\right)^{j-3},\\
|||X_{P_{j_0j}}|||^T_{p,\mathcal{D}_{j_0}\times\Pi_{\tilde\eta}}&\preceq&\label{083107}
\rho\left(\frac1{\tilde\eta^2}\mathcal{N}^{2(j_0+5)^2}\rho\right)^{j-3},\\
\label{083108}
|||X_{Q_{j_0j}}|||^T_{p,\mathcal{D}_{j_0}\times\Pi_{\tilde\eta}}&\preceq&\rho\left(\frac1{\tilde\eta^2}\mathcal{N}^{2(j_0+5)^2}\rho\right)^{j-3},\end{eqnarray}
where $a\preceq b$ means there is a constant $c>0$ depending on
$s_0,n$ and $\mathcal{M}$ such that $a\leq cb$ (but independent of
$\rho$, $\tilde\eta$ and $\mathcal{N}$). Then there exists a
symplectic map
$\Phi_{j_0}:\mathcal{D}_{j_0+1}\rightarrow\mathcal{D}_{j_0}$, such
that
\begin{eqnarray*}
&&H_{j_0+1}(x,y,q,\bar q;\xi)\\&:=&H_{j_0}\circ \Phi_{j_0}(x,y,q,\bar q;\xi)\\
&=&\breve N(y,q,\bar q;\xi)+Z_{j_0+1}(y,q,\bar
q;\xi)+P_{j_0+1}(x,y,q,\bar q;\xi) +Q_{j_0+1}(x,y,q,\bar q;\xi),
\end{eqnarray*}
where
\begin{eqnarray*}
Z_{j_0+1}(y,q,\bar q;\xi)&=&\sum_{3\leq j\leq
j_0+1}Z_{(j_0+1)j}(y,q,\bar
q;\xi),\\
P_{j_0+1}(x,y,q,\bar q;\xi)&=&\sum_{j\geq
j_0+2}P_{(j_0+1)j}(x,y,q,\bar
q;\xi),\\
Q_{j_0+1}(x,y,q,\bar q;\xi)&=&\sum_{j\geq 3}Q_{(j_0+1)j}(x,y,q,\bar
q;\xi),
\end{eqnarray*}
with
\begin{eqnarray*}Z_{(j_0+1)j}(y,q,\bar q;\xi)&=&\sum_{2|\alpha|+2|\beta|+2|\mu|=j, |\mu|\leq1}
Z_{j_0+1}^{\alpha\beta\beta\mu\mu}(\xi)y^{\alpha}\tilde{q}
^{\beta}\bar{\tilde{q}}^{\beta}\hat{q}^{\mu}\bar{\hat{q}}^{\mu},\\
P_{(j_0+1)j}(x,y,q,\bar
q;\xi)&=&\sum_{2|\alpha|+|\beta|+|\gamma|+|\mu|+|\nu|= j,
|\mu|+|\nu|\leq 2}
P_{j_0+1}^{\alpha\beta\gamma\mu\nu}(x;\xi)y^{\alpha}\tilde{q}^{\beta}
\bar{\tilde{q}}^{\gamma}\hat{q}^{\mu}\bar{\hat{q}}^{\nu},\\
Q_{(j_0+1)j}(x,y,q,\bar
q;\xi)&=&\sum_{2|\alpha|+|\beta|+|\gamma|+|\mu|+|\nu|=
j,|\mu|+|\nu|\geq3}
Q_{j_0+1}^{\alpha\beta\gamma\mu\nu}(x;\xi)y^{\alpha}\tilde{q}
^{\beta}\bar{\tilde{q}}^{\gamma}\hat{q}^{\mu}\bar{\hat{q}}^{\nu}.\end{eqnarray*}
Moreover, the following estimates hold:\\
(1) the symplectic map $\Phi_{j_0}$ satisfies
\begin{eqnarray}\label{091011}
||\Phi_{j_0}-id||_{p,\mathcal{D}_{j_0}''}&\preceq&
\left(\frac1{\tilde\eta^2}\mathcal{N}^{2(j_0+5)^2}\rho\right)^{j_0-1}
\end{eqnarray}
and
\begin{equation}\label{090701}
|||D\Phi_{j_0}-Id|||_{p,\mathcal{D}_{j_0+1}}\preceq
\rho^{-1}\left(\frac1{\tilde\eta^2}\mathcal{N}^{2(j_0+5)^2}\rho\right)^{j_0-1}
;\end{equation}\\
(2) the Hamiltonian vector fields
$X_{Z_{(j_0+1)j}},X_{P_{(j_0+1)j}}$ and $X_{Q_{(j_0+1)j}}$ satisfy
\begin{eqnarray}
\label{091012}|||X_{Z_{(j_0+1)j}}|||^T_{p,\mathcal{D}_{j_0+1}\times\Pi_{\tilde\eta}}&\preceq&
\rho\left(\frac1{\tilde\eta^2}\mathcal{N}^{2(j_0+5)^2}\rho\right)^{j-3},\\
\label{091013}|||X_{P_{(j_0+1)j}}|||^T_{p,\mathcal{D}_{j_0+1}\times\Pi_{\tilde\eta}}&\preceq&
\rho\left(\frac1{\tilde\eta^2}\mathcal{N}^{2(j_0+6)^2}\rho\right)^{j-3},\\
\label{091014}|||X_{Q_{(j_0+1)j}}|||^T_{p,\mathcal{D}_{j_0+1}\times\Pi_{\tilde\eta}}
&\preceq&\rho\left(\frac1{\tilde\eta^2}\mathcal{N}^{2(j_0+6)^2}\rho\right)^{j-3}.
\end{eqnarray}
\end{lem}
\begin{proof}\textbf{Step 1. The derivative of homological
equation.}\\
Expand $P_{j_0}^{\alpha\beta\gamma\mu\nu}(x;\xi)$ into Fourier series
\begin{equation*}
P_{j_0}^{\alpha\beta\gamma\mu\nu}(x;\xi)=\sum_{k\in\mathbb{Z}^n}
\widehat{P_{j_0}^{\alpha\beta\gamma\mu\nu}}(k;\xi)e^{\sqrt{-1}\langle
k,x\rangle}.
\end{equation*}
To obtain the partial normal form of order $j_0+1$, we need to
eliminate all non-integrable terms in $ P_{j_0(j_0+1)}(x,y,q,\bar
q;\xi), $ which are
\begin{eqnarray*}
\sum_{2|\alpha|+|\beta|+|\gamma|+|\mu|+|\nu|= j_0+1, |\mu|+|\nu|\leq
2}\sum_{k\in\mathbb{Z}^n}\widehat{
P_{j_0}^{\alpha\beta\gamma\mu\nu}}(k;\xi)e^{\sqrt{-1}\langle k,x\rangle}y^{\alpha}\tilde{q}^{\beta}
\bar{\tilde{q}}^{\gamma}\hat{q}^{\mu}\bar{\hat{q}}^{\nu}
\end{eqnarray*}
with $|k|+|\beta-\gamma|+|\mu-\nu|\neq0$. To this end, let
\begin{equation}\label{090303}F_{j_0}(x,y,q,\bar
q;\xi)=\sum_{2|\alpha|+|\beta|+|\gamma|+|\mu|+|\nu|=j_0+1,|\mu|+|\nu|\leq
2}
F_{j_0}^{\alpha\beta\gamma\mu\nu}(x;\xi)y^{\alpha}\tilde{q}^{\beta}
\bar{\tilde{q}}^{\gamma}\hat{q}^{\mu}\bar{\hat{q}}^{\nu},
\end{equation}and let
$\Phi_{j_0}=X_{F_{j_0}}^t|_{t=1}$ be the time-$1$ map of the
Hamiltonian vector field $X_{F_{j_0}}$.

Using Taylor formula,
\begin{eqnarray}
H_{j_0+1}\nonumber&:=&H_{j_0}\circ
X^t_{F_{j_0}}|_{t=1}\\
\nonumber&=&(\breve N+Z_{j_0}+P_{j_0}+Q_{j_0})\circ
X^t_{F_{j_0}}|_{t=1}\\
&=&\label{090606}\breve N+\{\breve
N,F_{j_0}\}+\int^1_0(1-t)\{\{\breve N,F_{j_0}\},F_{j_0}\}
\circ X^t_{F_{j_0}}dt\\
&&\label{090607}+P_{j_0(j_0+1)}+\int_0^1\{P_{j_0(j_0+1)},F_{j_0}\}\circ X^t_{F_{j_0}}dt\\
&&\label{090608}+\left(Z_{j_0}+P_{j_0}-P_{j_0(j_0+1)}+Q_{j_0}\right)\circ
X_{F_{j_0}}^t|_{t=1}.
\end{eqnarray}
Then we obtain
the homological equation
\begin{equation}\label{008}
\{\breve N,F_{j_0}\}+P_{j_0(j_0+1)}=\widehat{Z_{j_0}},
\end{equation}
where
\begin{equation}\label{083109}
\widehat{Z_{j_0}}(y,q,\bar
q;\xi)=\sum_{2|\alpha|+2|\beta|+2|\mu|=j_0+1,|\mu|\leq1}
\widehat{P_{j_0}^{\alpha\beta\beta\mu\mu}}(0;\xi)
y^{\alpha}\tilde{q}^{\beta}\bar{\tilde{q}}^{\beta}
\hat{q}^{\mu}\bar{\hat{q}}^{\mu}. \end{equation}If the homological
equation (\ref{008}) is solvable, then in view of
(\ref{090606})-(\ref{090608}) we can define
\begin{equation}\label{090603}
Z_{j_0+1}(y,q,\bar q;\xi)=Z_{j_0}(y,q,\bar
q;\xi)+\widehat{Z_{j_0}}(y,q,\bar q;\xi),
\end{equation}
and
\begin{eqnarray} &&\nonumber
P_{j_0+1}+Q_{j_0+1}\\&=&\label{090610}\int^1_0(1-t)\{\{\breve
N,F_{j_0}\},F_{j_0}\}
\circ X^t_{F_{j_0}}dt\\
&&\label{090616}+\int_0^1\{P_{j_0(j_0+1)}+Z_{j_0},F_{j_0}\}\circ X^t_{F_{j_0}}dt\\
&&\label{090621}+\left(P_{j_0}-P_{j_0(j_0+1)}+Q_{j_0}\right)\circ
X_{F_{j_0}}^t|_{t=1}.
\end{eqnarray}
and $H_{j_0+1}(x,y,q,\bar q;\xi)$ has the following form
\begin{equation}
H_{j_0+1}=\breve N+Z_{j_0+1}+P_{j_0+1}+Q_{j_0+1}.
\end{equation}
\textbf{Step 2. The solution of homological equation (\ref{008}).}\\
By passing to Fourier coefficients, (\ref{008}) reads
\begin{eqnarray}\label{090707}
\widehat{F_{j_0}^{\alpha\beta\gamma\mu\nu}}(k;\xi)=
\frac{\widehat{P_{j_0}^{\alpha\beta\gamma\mu\nu}}(k;\xi)}{\sqrt{-1}(\langle
k,\breve\omega(\xi)\rangle+\langle
\beta-\gamma,\tilde\Omega(\xi)\rangle+\langle
\mu-\nu,\hat\Omega(\xi)\rangle)}\end{eqnarray} for
\begin{eqnarray*}
2|\alpha|+|\beta|+|\gamma|+|\mu|+|\nu|= j_0+1, |\mu|+|\nu|\leq
2,|k|+|\beta-\gamma|+|\mu-\nu|\neq0,
\end{eqnarray*}and otherwise
\begin{eqnarray*}
\widehat{F_{j_0}^{\alpha\beta\gamma\mu\nu}}(k;\xi)=
0.\end{eqnarray*}

Now we will estimate
$|||X_{F_{j_0}}|||^T_{p,\mathcal{D}_{j_0}'\times\Pi_{\tilde\eta}}.$

For each $\xi\in\Pi_{\tilde\eta}$, the frequencies
$\breve\omega(\xi)$ and $\breve\Omega(\xi)$ satisfy the
$(\tilde\eta,\mathcal{N},\mathcal{M})$-non-resonant conditions
(\ref{090302}), i.e.
\begin{equation*}
|\langle k,\breve\omega(\xi)\rangle+\langle
\beta-\gamma,\tilde\Omega(\xi)\rangle+\langle
\mu-\nu,\hat\Omega(\xi)\rangle|\geq \frac{\tilde
\eta}{4^{\mathcal{M}}(|k|+1)^{\tau}C(\mathcal{N},\beta-\gamma)}.
\end{equation*}
Then
\begin{eqnarray}
&&\left|\left(\langle k,\breve\omega(\xi)\rangle+\langle
\beta-\gamma,\tilde\Omega(\xi)\rangle+\langle
\mu-\nu,\hat\Omega(\xi)\rangle\right)^{-1}\right|\nonumber\\
\nonumber&\leq& \frac{4^{\mathcal{M}}}{\tilde
\eta}{(|k|+1)^{\tau}C(\mathcal{N},\beta-\gamma)} \nonumber\\&\leq&
\frac{4^{\mathcal{M}}}{\tilde
\eta}{(|k|+1)^{\tau}\mathcal{N}^{(|\beta-\gamma|+4)^2}}\nonumber\qquad\qquad\mbox{(in view of (\ref{121206}))}\\
\nonumber&\leq& \frac{4^{\mathcal{M}}}{\tilde
\eta}{(|k|+1)^{\tau}\mathcal{N}^{(|\beta|+|\gamma|+4)^2}}\\
&\leq&\frac{4^{\mathcal{M}}}{\tilde
\eta}{(|k|+1)^{\tau}\mathcal{N}^{(j_0+5)^2}}\label{083104}\qquad\qquad
\mbox{(in view of $|\beta|+|\gamma|\leq j_0+1$)}.
\end{eqnarray}
Moreover, in view of (\ref{081902}) in Assumption B (twist
conditions)
\begin{eqnarray*}
\partial_{\xi_j}\omega_i(\xi)=\delta_{ji}, \quad
\partial_{\xi_j}\Omega_{j'}(\xi)=\delta_{j(n+j')},\qquad 1\leq i\leq n,\ j,j'\ge 1,
\end{eqnarray*}
and the estimates (see (\ref{080103}) and (\ref{080104}))
\begin{equation*}
||\breve\omega(\xi)-\omega(\xi)||+\sup_{j\geq1}||\partial_{\xi_j}(\breve\omega(\xi)-\omega(\xi))||\leq
{c\eta^8\epsilon},\end{equation*}and
\begin{equation*}
||\breve\Omega(\xi)-\Omega(\xi)||+\sup_{j\geq1}||\partial_{\xi_j}(\breve\Omega(\xi)-\Omega(\xi))||\leq
{c\eta^8\epsilon},
\end{equation*}
we have
\begin{equation}\label{091203}
\partial_{\xi_i}\breve{\omega}_{i}(\xi)\geq1-{c\eta^8\epsilon}, \qquad
\partial_{\xi_j}\breve{\Omega}_{n+j}(\xi)\geq1-{c\eta^8\epsilon},
\end{equation}
\begin{equation}\label{091204}
|\partial_{\xi_j}\breve\omega_{i}(\xi)|\leq {c\eta^8\epsilon},\qquad
j\neq i,
\end{equation}and
\begin{equation}\label{091205}|\partial_{\xi_j}\breve\Omega_{n+j'}(\xi)|\leq
{c\eta^8\epsilon},\qquad j\neq j'.
\end{equation}Then for each $j\geq 1$, we
get
\begin{eqnarray}
\nonumber&&\left|\partial_{\xi_j}\left(\langle
k,\breve\omega(\xi)\rangle+\langle
\beta-\gamma,\tilde\Omega(\xi)\rangle+\langle
\mu-\nu,\hat\Omega(\xi)\rangle\right)\right|\\
\nonumber&\leq&
(|k|+j_0+1)(c\eta^8\epsilon+1)\qquad \mbox{(in view of the estimates (\ref{091203})-(\ref{091205}))}\\
\nonumber&\leq& 2(|k|+j_0+1)\\
&\leq& 2(|k|+\mathcal{M}+2)\qquad \mbox{(in view of $2\leq j_0\leq
\mathcal{M}+1$)}\label{090601}.
\end{eqnarray} Hence,
\begin{eqnarray}
&&\nonumber\left|\partial_{\xi_j}\left(\langle
k,\breve\omega(\xi)\rangle+\langle
\beta-\gamma,\tilde\Omega(\xi)\rangle+\langle
\mu-\nu,\hat\Omega(\xi)\rangle\right)^{-1}\right|\\
&=&\left|\frac{\partial_{\xi_j}\left(\langle
k,\breve\omega(\xi)\rangle+\langle
\beta-\gamma,\tilde\Omega(\xi)\rangle+\langle
\mu-\nu,\hat\Omega(\xi)\rangle\right)}{\left(\langle
k,\breve\omega(\xi)\rangle+\langle
\beta-\gamma,\tilde\Omega(\xi)\rangle+\langle
\mu-\nu,\hat\Omega(\xi)\rangle\right)^2}\right|\nonumber\\
&\leq&\frac{2\cdot 4^{2\mathcal{M}}}{\tilde \eta^2}{(|k|+1)^{2\tau}
\mathcal{N}^{2(j_0+5)^2}(|k|+\mathcal{M}+2)}\nonumber\\&&
\nonumber\mbox{(based on the
inequalities (\ref{083104}) and (\ref{090601}))}\\
&\leq& \frac{2\cdot 4^{2\mathcal{M}}}{\tilde
\eta^2}{(\mathcal{M}+2)(|k|+1)^{2\tau+1} \mathcal{N}^{2(j_0+5)^2}}
,\label{083105}
\end{eqnarray}where the last inequality is based on $|k|+\mathcal{M}+2\leq (|k|+1)(\mathcal{M}+2).$ Then in view of the formula (\ref{090707}) and the inequalities
(\ref{083104}) and (\ref{083105}), for each $j\geq 1$, we obtain
\begin{eqnarray}
&&\nonumber|\widehat{F_{j_0}^{\alpha\beta\gamma\mu\nu}}(k;\xi)|+|\partial_{\xi_j}\widehat{F_{j_0}^{\alpha\beta\gamma\mu\nu}}(k;\xi)|
\\&\preceq&\label{091008}\frac{(|k|+1)^{2\tau+1}\mathcal{N}^{2(j_0+5)^2}}{\tilde\eta^2}
\left(|\widehat{P_{j_0}
^{\alpha\beta\gamma\mu\nu}}(k;\xi)|+|\partial_{\xi_j}\widehat{P_{j_0}
^{\alpha\beta\gamma\mu\nu}}(k;\xi)|\right).\end{eqnarray}  Noting
$\tau> n+3$ will be fixed and recalling $s'=s_0/(12\mathcal{M})$
(see (\ref{090301})), we obtain
\begin{eqnarray}
&&\label{1191}\nonumber
|||X_{F_{j_0}}|||^T_{p,\mathcal{D}_{j_0}'\times\Pi_{\tilde\eta}}\\
\nonumber&\preceq& \frac{{\mathcal{N}^{2(j_0+5)^2}}}{\tilde\eta^2}
|||X_{P_{j_0(j_0+1)}}|||^T_{p,\mathcal{D}_{j_0}\times\Pi_{\tilde\eta}}\nonumber\\&&\nonumber\mbox{(based
on the inequality (\ref{091008}) and
Theorem \ref{0005})}\\
&\preceq&\nonumber\frac{{\mathcal{N}^{2(j_0+5)^2}}}{\tilde\eta^2}\cdot
\rho\left(\frac1{\tilde\eta^2}\mathcal{N}^{2(j_0+5)^2}\rho\right)^{j_0-2}\quad \mbox{(based on the inequality (\ref{083107}) for $j=j_0+1$)}\\
&=&\label{090401}\left(\frac1{\tilde\eta^2}\mathcal{N}^{2(j_0+5)^2}\rho\right)^{j_0-1}\label{083112}.\end{eqnarray}
Noting $\Phi_{j_0}=X_{F_{j_0}}^t|_{t=1}$ and basing on the inequality
(\ref{090401}) and
 Theorem \ref{081503}, we obtain
 \begin{equation*}\label{080301}
||\Phi_{j_0}-id||_{p,\mathcal{D}_{j_0}''}\preceq
\left(\frac1{\tilde\eta^2}\mathcal{N}^{2(j_0+5)^2}\rho\right)^{j_0-1}
.\end{equation*}Moreover,
\begin{equation*}
|||D\Phi_{j_0}-Id|||_{p,\mathcal{D}_{j_0+1}}\preceq\rho^{-1}
\left(\frac1{\tilde\eta^2}\mathcal{N}^{2(j_0+5)^2}\rho\right)^{j_0-1}
,\end{equation*}which follows the generalized Cauchy estimate in
Lemma \ref{081907}. We finish the proof of the inequalities
(\ref{091011}) and (\ref{090701}).
\\
\textbf{Step 3. Estimate $|||X_{Z_{j_0+1}}|||^T_{p,\mathcal{D}_{j_0+1}
\times\Pi_{\tilde\eta}},|||X_{P_{j_0+1}}|||^T_{p,\mathcal{D}_{j_0+1}
\times\Pi_{\tilde\eta}}$ and $|||X_{Q_{j_0+1}}|||^T_{p,\mathcal{D}_{j_0+1}\times\Pi_{\tilde\eta}}$.} \\
Following the notation as (\ref{090602}), rewrite $Z_{j_0+1}(y,q,\bar q;\xi)$ as
\begin{equation*}
Z_{j_0+1}(y,q,\bar q;\xi)=\sum_{3\leq j\leq
j_0+1}Z_{(j_0+1)j}(y,q,\bar q;\xi),
\end{equation*}where
\begin{equation*}Z_{(j_0+1)j}(y,q,\bar q;\xi)=\sum_{2|\alpha|+2|\beta|+2|\mu|=j, |\mu|\leq1}
Z_{j_0+1}^{\alpha\beta\beta\mu\mu}(\xi)y^{\alpha}\tilde{q}
^{\beta}\bar{\tilde{q}}^{\beta}\hat{q}^{\mu}\bar{\hat{q}}^{\mu}.
\end{equation*}
In view of (\ref{090602}), (\ref{083109}) and (\ref{090603}), we
have
\begin{equation}\label{090605}
Z_{(j_0+1)j}(y,q,\bar q;\xi)={Z_{j_0j}}(y,q,\bar q;\xi),\qquad 3\leq
j\leq j_0
\end{equation}
and
\begin{equation}\label{090611} Z_{(j_0+1)(j_0+1)}(y,q,\bar
q;\xi)=\widehat{Z_{j_0}}(y,q,\bar q;\xi).
\end{equation}
For $3\leq j\leq j_0$,
\begin{eqnarray}
&&\nonumber|||X_{{Z_{(j_0+1)j}}}|||^T_{p,\mathcal{D}_{j_0+1}\times\Pi_{\tilde\eta}}
\\
&\preceq&\nonumber
\rho\left(\frac1{\tilde\eta^2}\mathcal{N}^{2(j_0+4)^2}\rho\right)^{j-3}\qquad
\mbox{(in view of the inequality (\ref{083106}))}\\
&\leq&
\rho\left(\frac1{\tilde\eta^2}\mathcal{N}^{2(j_0+5)^2}\rho\right)^{j-3}\label{090201}.
\end{eqnarray}
When $j=j_0+1$, we have\begin{eqnarray}
&&\nonumber|||X_{{Z_{(j_0+1)(j_0+1)}}}|||^T_{p,\mathcal{D}_{j_0+1}\times
\Pi_{\tilde\eta}}\\&=&\nonumber|||X_{\widehat{Z_{j_0}}}|||^T_{p,\mathcal{D}_{j_0+1}\times\Pi_{\tilde\eta}}\qquad
\qquad\mbox{(in view of (\ref{090611}))}
\\&\leq&
|||X_{P_{j_0(j_0+1)}}|||^T_{p,\mathcal{D}_{j_0+1}\times\Pi_{\tilde\eta}}\nonumber\qquad\quad\mbox{(in
view of (\ref{083109}))}\\&\preceq&
\rho\left(\frac1{\tilde\eta^2}\mathcal{N}^{2(j_0+5)^2}\rho\right)^{j_0-2}\label{091010}\
\mbox{(in view of the inequality (\ref{083107}) for $j=j_0+1$)}.
\end{eqnarray}
Based on the inequalities (\ref{090201}) and (\ref{091010}), we
finish the proof of the inequality (\ref{091012}).

To obtain the estimates (\ref{091013}) and (\ref{091014}),  we will
estimate the $p$-tame norm of the terms
(\ref{090610})-(\ref{090621}) respectively.

Firstly we consider the term (\ref{090610}). Let
\begin{equation}\label{101001}
\breve N^{(0)}=\breve N\circ X_{F_{j_0}}^t|_{t=0}=\breve N,
\end{equation}
and \begin{equation}\label{101002} \breve N^{(j)}=\{\breve N^{(j-1)},F_{j_0}\},
\qquad j\geq1.
\end{equation}
Using Taylor formula again, we have
\begin{equation}\label{091020}
\int^1_0(1-t)\{\{\breve N,F_{j_0}\},F_{j_0}\} \circ
X^t_{F_{j_0}}dt=\sum_{j\geq 2}\frac1{j!}\breve N^{(j)}.
\end{equation}
Moreover, note that $\breve N^{(j)}(x,y,q,\bar q;\xi)$ has the
following form
\begin{equation}\label{090613}
\breve N^{(j)}(x,y,q,\bar
q;\xi)=\sum_{2|\alpha|+|\beta|+|\gamma|+|\mu|+|\nu|=j(j_0+1)+2-2j}\breve
N^{(j)\alpha\beta\gamma\mu\nu}(x,\xi)y^{\alpha}\tilde{q}^{\beta}
\bar{\tilde{q}}^{\gamma}\hat{q}^{\mu}\bar{\hat{q}}^{\nu}.
\end{equation}In view of the homological equation (\ref{008})
$$\{\breve N,F_{j_0}\}+P_{j_0(j_0+1)}=\widehat{Z_{j_0}},$$and the formulas (\ref{101001})
and (\ref{101002}), we have
\begin{equation}
\breve N^{(1)}=-P_{j_0(j_0+1)}+\widehat{Z_{j_0}}.
\end{equation}
In view of the inequality (\ref{083107}) for $j=j_0+1$ and the formula (\ref{083109}),
\begin{equation}\label{090612}
|||X_{\breve N^{(1)}}|||^T_{p,\mathcal{D}_{j_0}\times
\Pi_{\tilde\eta}}\leq
|||X_{P_{j_0(j_0+1)}}|||^T_{p,\mathcal{D}_{j_0}\times
\Pi_{\tilde\eta}}\preceq
\rho\left(\frac1{\tilde\eta^2}\mathcal{N}^{2(j_0+5)^2}\rho\right)^{j_0-2}.
\end{equation}
For $j\geq 2$, let $$s_j=\frac{s_0}{4j\mathcal{M}}\qquad \mbox{ and}
\qquad \rho_j=\frac{\rho}{2j\mathcal{M}}.$$Then we have
\begin{eqnarray}\nonumber
&&\frac1{j!}|||X_{\breve
N^{(j)}}|||_{p,\mathcal{D}_{j_0+1}\times\Pi_{\tilde\eta}}^{T}\\&\leq&\nonumber
\frac1{j!}\left(\left(4C\max\left\{\frac{4j\mathcal{M}}{s_0},10j\mathcal{M}\right\}\right)^{j-1}
\left(|||X_{\breve
N^{(1)}}|||_{p,\mathcal{D}_{j_0}\times\Pi_{\tilde\eta}}^{T}\right)\left(|||X_{F_{j_0}}|||_{p,\mathcal{D}_{j_0'}\times\Pi_{\tilde\eta}}^{T}\right)^{j-1}\right)\\&&\mbox{(based
on Theorem
{\ref{021102}})}\nonumber\\
&\preceq&\nonumber \frac{j^{j-1}}{j!}
\left(\rho\left(\frac1{\tilde\eta^2}\mathcal{N}^{2(j_0+5)^2}\rho\right)^{j_0-2}\right)
\left( C_1\left(\frac1{\tilde\eta^2}\mathcal{N}^{2(j_0+5)^2}\rho\right)^{j_0-1}\right)^{j-1}\\
&&\nonumber\mbox{(based on the inequalities (\ref{083112}) and (\ref{090612}),}\\&&\nonumber \mbox{and $C_1>0$ is a constant depending on $s_0$ and $\mathcal{M}$)}\\
&\leq&\nonumber \frac{e}{j}
\left(\rho\left(\frac1{\tilde\eta^2}\mathcal{N}^{2(j_0+5)^2}\rho\right)^{j_0-2}\right)
\left(\tilde c
e\left(\frac1{\tilde\eta^2}\mathcal{N}^{2(j_0+5)^2}\rho\right)^{j_0-1}\right)^{j-1}\qquad
\\&&\nonumber\mbox{(using the inequality $j^j<j!e^j$)}\nonumber\\&\leq&\nonumber
\frac{e}{j}
\left(\rho\left(\frac1{\tilde\eta^2}\mathcal{N}^{2(j_0+5)^2}\rho\right)^{j_0-2}\right)
\left(\left(\frac1{\tilde\eta^2}\mathcal{N}^{2(j_0+6)^2}\rho\right)^{j_0-1}\right)^{j-1}\\&&
\nonumber\mbox{(noting $\mathcal{N}$ is large depending on $s_0,r_0,n$ and $\mathcal{M}$)}\\
&\preceq&\label{090615}\rho\left(\frac1{\tilde\eta^2}\mathcal{N}^{2(j_0+6)^2}
\rho\right)^{j(j_0+1)-1-2j}.
\end{eqnarray}
In view of (\ref{090613}), note that the index
$\alpha,\beta,\gamma,\mu,\nu$ of $\breve N^{(j)} $ satisfy
$$2|\alpha|+|\beta|+|\gamma|+|\mu|+|\nu|=j(j_0+1)+2-2j\geq j_0+2 \ (\mbox{since $j,j_0\geq 2$})$$
and $$j(j_0+1)-1-2j=(j(j_0+1)+2-2j)-3.$$Then in view of
(\ref{091020}) and the inequality (\ref{090615}), we finish the
estimate of $p$-tame norm of the term (\ref{090610}).

Next, we will give the estimate of $p$-tame norm of the term (\ref{090616}). Let $$W_i=Z_{j_0i},\qquad
\mbox{for}\ 3\leq i\leq j_0,$$ and $$W_{j_0+1}=P_{j_0(j_0+1)}.$$Then
\begin{equation}\label{091021}\sum_{3\leq i\leq j_0+1}W_i=P_{j_0(j_0+1)}+Z_{j_0}.
\end{equation} Let
\begin{equation*} W_i^{(0)}=W_i\circ X_{F_{j_0}}^t|_{t=0}=W_i,
\end{equation*}
and \begin{equation*} W^{(j)}_i=\{W_i^{(j-1)},F_{j_0}\}, \qquad
j\geq1.
\end{equation*}
In view of (\ref{091021}) and using Taylor formula, we have
\begin{equation}\label{091023}
\int_0^1\{P_{j_0(j_0+1)}+Z_{j_0},F_{j_0}\}\circ
X^t_{F_{j_0}}dt=\sum_{3\leq i\leq j_0+1}\left(\sum_{j\geq
1}\frac1{j!}W_i^{(j)}\right).
\end{equation}
Note that $W_i^{(j)}(x,y,q,\bar q;\xi)$ has the following form
\begin{equation}\label{091022}
W_i^{(j)}(x,y,q,\bar
q;\xi)=\sum_{2|\alpha|+|\beta|+|\gamma|+|\mu|+|\nu|=j(j_0+1)+i-2j}W_i^{(j)\alpha\beta\gamma\mu\nu}(x,\xi)y^{\alpha}\tilde{q}^{\beta}
\bar{\tilde{q}}^{\gamma}\hat{q}^{\mu}\bar{\hat{q}}^{\nu}.
\end{equation}
 Using the proof of  (\ref{090615}) and in view of (\ref{083106}) and (\ref{083107}), we have
\begin{eqnarray}
\frac1{j!}|||X_{W_i^{(j)}}|||_{p,\mathcal{D}_{j_0+1}\times\Pi_{\tilde\eta}}^{T}
\preceq\label{091024}\rho\left(\frac1{\tilde\eta^2}\mathcal{N}^{2(j_0+6)^2}
\rho\right)^{(j(j_0+1)-i-2j)-3}.
\end{eqnarray}
In view of (\ref{091022}), note that the indices
$\alpha,\beta,\gamma,\mu,\nu$ of $W_i^{(j)}$ satisfy
$$2|\alpha|+|\beta|+|\gamma|+|\mu|+|\nu|=j(j_0+1)+i-2j\geq j_0+2,$$
since $j\geq 1,j_0\geq 2,i\geq 3.$ Then in view of (\ref{091023})
and the inequality (\ref{091024}), we finish the estimate of
$p$-tame norm of the term (\ref{090616}).

Finally, we will give the estimate of $p$-tame norm of the term (\ref{090621}). Let $$U_i=P_{j_0i},\qquad
\mbox{for}\ i\geq j_0+2,$$ and $$V_{i}=Q_{j_0i},\qquad \mbox{for}\
i\geq 3.$$Then
\begin{equation}\label{090617}\sum_{i\geq j_0+2}U_i+\sum_{i\geq 3}V_i=P_{j_0}-P_{j_0(j_0+1)}+Q_{j_0}.
\end{equation}For simplicity, denote $T_i=U_i$ or $V_i$. Let
\begin{equation*} T_i^{(0)}=T_i\circ X_{F_{j_0}}^t|_{t=0}=T_i,
\end{equation*}
and \begin{equation*} T^{(j)}_i=\{T_i^{(j-1)},F_{j_0}\}, \qquad
j\geq1.
\end{equation*}
Using Taylor formula again, we have
\begin{eqnarray}&&\left(P_{j_0}-P_{j_0(j_0+1)}+Q_{j_0}\right)\circ
X_{F_{j_0}}^t|_{t=1}\nonumber\\
&=&P_{j_0}-P_{j_0(j_0+1)}+Q_{j_0}+\sum_{i\geq
j_0+2}\left(\sum_{j\geq 1}\frac1{j!}U_i^{(j)}\right)+\sum_{i\geq
3}\left(\sum_{j\geq 1}\frac1{j!}V_i^{(j)}\right)\label{090622}.
\end{eqnarray}
For $j\geq 1$, note that $T_i^{(j)}(x,y,q,\bar q;\xi)$ has the
following form
\begin{equation}\label{090620}
T_i^{(j)}(x,y,q,\bar
q;\xi)=\sum_{2|\alpha|+|\beta|+|\gamma|+|\mu|+|\nu|=j(j_0+1)+i-2j}T_i^{(j)\alpha\beta\gamma\mu\nu}(x,\xi)y^{\alpha}\tilde{q}^{\beta}
\bar{\tilde{q}}^{\gamma}\hat{q}^{\mu}\bar{\hat{q}}^{\nu}.
\end{equation}
Then basing on the inequalities (\ref{083107}) and (\ref{083108}), and following the proof of the inequalities (\ref{090615}) again, we have
\begin{eqnarray}
\frac1{j!}|||X_{T_i^{(j)}}|||_{p,\mathcal{D}_{j_0+1}\times\Pi_{\tilde\eta}}^{T}
\preceq
\label{091025}\rho\left(\frac1{\tilde\eta^2}\mathcal{N}^{2(j_0+6)^2}
\rho\right)^{(j(j_0+1)-i-2j)-3}.
\end{eqnarray}
In view of (\ref{090620}), note that the index
$\alpha,\beta,\gamma,\mu,\nu$ of $T_i^{(j)}$ satisfy
$$2|\alpha|+|\beta|+|\gamma|+|\mu|+|\nu|=j(j_0+1)+i-2j\geq j_0+2,$$
since $j\geq 1,j_0\geq 2,i\geq3.$ Then in view of (\ref{090622}) and
the inequality (\ref{091025}), we finish the estimate of $p$-tame
norm of the term (\ref{090621}). Hence in view of the inequalities
(\ref{090615}), (\ref{091024}) and (\ref{091025}), we obtain
\begin{eqnarray*}
|||X_{P_{(j_0+1)j}}|||^T_{p,\mathcal{D}_{j_0}\times\Pi_{\tilde\eta}}&\preceq&
\rho\left(\frac1{\tilde\eta^2}\mathcal{N}^{2(j_0+6)^2}\rho\right)^{j-3},\\
|||X_{Q_{(j_0+1)j}}|||^T_{p,\mathcal{D}_{j_0}\times\Pi_{\tilde\eta}}
&\preceq&\rho\left(\frac1{\tilde\eta^2}\mathcal{N}^{2(j_0+6)^2}\rho\right)^{j-3}.
\end{eqnarray*}Then we finish the proof of Lemma \ref{052003}.
\end{proof}

\noindent {\bf Proof of Theorem \ref{thm7.1}.}
\begin{proof} By Lemma \ref{052003}, we will finish the proof of Theorem \ref{thm7.1}. By Theorem \ref{T1}, we obtain a normal form of order 2
around KAM tori, which is\begin{equation*}\breve H(x,y,q,\bar
q;\xi)= \breve N(y,q,\bar q;\xi)+ \breve R(x,y,q,\bar q;\xi),
\end{equation*}where
\begin{equation*}
\breve N(y,q,\bar
q;\xi)=\sum_{j=1}^{n}\breve{\omega}_j(\xi)y_j+\sum_{j\geq1}\breve{\Omega}_j(\xi)q_j\bar
q_j
\end{equation*}
and
\begin{equation*}
\breve R(x,y,q,\bar
q;\xi)=\sum_{\alpha\in\mathbb{N}^n,\beta,\gamma\in\mathbb{N}^{\mathbb{N}},2|\alpha|+|\beta|+|\gamma|\geq
3}\breve
R^{\alpha\beta\gamma}(x;\xi)y^{\alpha}q^{\beta}\bar{q}^{\gamma}.
\end{equation*}
Then following the notations in Lemma \ref{052003}, denote
$$H_2(x,y,q,\bar q;\xi)=\breve H(x,y,q,\bar q;\xi)$$ and
$\breve R(x,y,q,\bar q;\xi)$ can be rewritten as
\begin{equation}\label{090702}
\breve R(x,y,q,\bar q;\xi)=Z_2(y,q,\bar q;\xi)+P_2(x,y,q,\bar
q;\xi)+Q_2(x,y,q,\bar q;\xi),
\end{equation}
where
\begin{eqnarray}
\label{101101}
Z_{2}(y,q,\bar q;\xi)&=&\sum_{3\leq j\leq 2}Z_{2j}(y,q,\bar
q;\xi)=0,\\
\nonumber P_{2}(x,y,q,\bar q;\xi)&=&\sum_{j\geq 3}P_{2j}(x,y,q,\bar
q;\xi),\\
\nonumber Q_{2}(x,y,q,\bar q;\xi)&=&\sum_{j\geq 3}Q_{2j}(x,y,q,\bar q;\xi),
\end{eqnarray}
and
\begin{eqnarray}
\nonumber Z_{2j}(y,q,\bar
q;\xi)&=&\sum_{2|\alpha|+2|\beta|+2|\mu|=j, |\mu|\leq1}
Z_{2}^{\alpha\beta\beta\mu\mu}(\xi)y^{\alpha}\tilde{q}
^{\beta}\bar{\tilde{q}}^{\beta}\hat{q}^{\mu}\bar{\hat{q}}^{\mu},\\
\label{090703}P_{2j}(x,y,q,\bar
q;\xi)&=&\sum_{2|\alpha|+|\beta|+|\gamma|+|\mu|+|\nu|= j,
|\mu|+|\nu|\leq 2}
P_{2}^{\alpha\beta\gamma\mu\nu}(x;\xi)y^{\alpha}\tilde{q}^{\beta}
\bar{\tilde{q}}^{\gamma}\hat{q}^{\mu}\bar{\hat{q}}^{\nu},\\
\label{090704}Q_{2j}(x,y,q,\bar
q;\xi)&=&\sum_{2|\alpha|+|\beta|+|\gamma|+|\mu|+|\nu|=
j,|\mu|+|\nu|\geq3}
Q_{2}^{\alpha\beta\gamma\mu\nu}(x;\xi)y^{\alpha}\tilde{q}
^{\beta}\bar{\tilde{q}}^{\gamma}\hat{q}^{\mu}\bar{\hat{q}}^{\nu}.\end{eqnarray}
Denote $$W_j(x,y,q,\bar q;\xi)=P_{2j}(x,y,q,\bar q;\xi)\quad \mbox
{or} \quad Q_{2j}(x,y,q,\bar q;\xi), \qquad j\geq 3.$$Then
\begin{eqnarray*}
&&|||X_{W_j}|||^T_{p,D(s_0/2,r_0/2,r_0/2)\times\Pi_{\eta}}\\
&\leq&|||X_{\breve
R}|||^T_{p,D(s_0/2,r_0/2,r_0/2)\times\Pi_{\eta}}\\&\leq&
 \varepsilon(1+{c\eta^6\epsilon})\qquad \qquad\qquad\qquad\qquad\mbox{(by (\ref{080105}))}.
\end{eqnarray*}
Hence, by replacing $r_0/2$ with $5\rho$, we obtain
\begin{eqnarray}\label{090705}|||X_{W_j}|||^T_{p,D(s_0/2,5\rho,5\rho)\times\Pi_{\eta}}
\leq \varepsilon(1+{c\eta^6\epsilon})\left(\frac{10\rho}{r_0}\right)^{j-2} .
\end{eqnarray}
Moreover, in view of the inequality (\ref{091113}) (i.e. noting
$\mathcal{N}$ is large depending on $s_0,r_0,n$ and $\mathcal{M}$)
and $0<\tilde\eta<1$,
\begin{eqnarray}\label{090706}|||X_{W_j}|||^T_{p,D(s_0/2,5\rho,5\rho)\times\Pi_{\eta}}
\leq\rho\left(\frac1{\tilde\eta^2}\mathcal{N}\rho\right)^{j-3}.
\end{eqnarray}In view of the formula (\ref{101101}) and the inequality (\ref{090706}), the assumptions (\ref{083106})-(\ref{083108}) in Lemma \ref{052003}
hold for $j_0=2$.

Let
\begin{eqnarray}\label{091108}
\Phi=\Phi_{2}\circ\Phi_3\circ\cdots\circ\Phi_{\mathcal{M}+1}.
\end{eqnarray}
Then, based on Lemma \ref{052003}, $\breve{\breve H}:=\breve
H\circ\Phi$ has the form of
\begin{equation*}
\breve{\breve H}(x,y,q,\bar q;\xi)={\breve N}(y,q,\bar q;\xi)+{
Z}_{\mathcal{M}+2}(y,q,\bar q;\xi)+P_{\mathcal{M}+2}(x,y,q,\bar
q;\xi)+Q_{\mathcal{M}+2}(x,y,q,\bar q;\xi),
\end{equation*}where
\begin{eqnarray*}
Z_{\mathcal{M}+2}(y,q,\bar q;\xi)&=&\sum_{3\leq j\leq
\mathcal{M}+2}Z_{(\mathcal{M}+2)j}(y,q,\bar
q;\xi),\\
P_{\mathcal{M}+2}(x,y,q,\bar q;\xi)&=&\sum_{j\geq
\mathcal{M}+3}P_{(\mathcal{M}+2)j}(x,y,q,\bar
q;\xi),\\
Q_{\mathcal{M}+2}(x,y,q,\bar q;\xi)&=&\sum_{j\geq
3}Q_{(\mathcal{M}+2)j}(x,y,q,\bar q;\xi),
\end{eqnarray*}
with
\begin{eqnarray*}Z_{(\mathcal{M}+2)j}(y,q,\bar q;\xi)&=&\sum_{2|\alpha|+2|\beta|+2|\mu|=j, |\mu|\leq1}
Z_{\mathcal{M}+2}^{\alpha\beta\beta\mu\mu}(\xi)y^{\alpha}\tilde{q}
^{\beta}\bar{\tilde{q}}^{\beta}\hat{q}^{\mu}\bar{\hat{q}}^{\mu},\\
P_{(\mathcal{M}+2)j}(x,y,q,\bar
q;\xi)&=&\sum_{2|\alpha|+|\beta|+|\gamma|+|\mu|+|\nu|= j,
|\mu|+|\nu|\leq 2}
P_{\mathcal{M}+2}^{\alpha\beta\gamma\mu\nu}(x;\xi)y^{\alpha}\tilde{q}^{\beta}
\bar{\tilde{q}}^{\gamma}\hat{q}^{\mu}\bar{\hat{q}}^{\nu},\\
Q_{(\mathcal{M}+2)j}(x,y,q,\bar
q;\xi)&=&\sum_{2|\alpha|+|\beta|+|\gamma|+|\mu|+|\nu|=
j,|\mu|+|\nu|\geq3}
Q_{\mathcal{M}+2}^{\alpha\beta\gamma\mu\nu}(x;\xi)y^{\alpha}\tilde{q}
^{\beta}\bar{\tilde{q}}^{\gamma}\hat{q}^{\mu}\bar{\hat{q}}^{\nu}.\end{eqnarray*}
Moreover, the following estimates hold:
\begin{eqnarray}
\label{091104}|||X_{Z_{(\mathcal{M}+2)j}}|||^T_{p,{D}(s_0/4,4\rho,4\rho)
\times\Pi_{\tilde\eta}}&\preceq&
\rho\left(\frac1{\tilde\eta^2}\mathcal{N}^{2(\mathcal{M}+6)^2}\rho\right)^{j-3},\\
\label{091105}|||X_{P_{(\mathcal{M}+2)j}}|||^T_{p,{D}(s_0/4,4\rho,4\rho)
\times\Pi_{\tilde\eta}}&\preceq&
\rho\left(\frac1{\tilde\eta^2}\mathcal{N}^{2(\mathcal{M}+7)^2}\rho\right)^{j-3},\\
\label{091106}|||X_{Q_{(\mathcal{M}+2)j}}|||^T_{p,{D}(s_0/4,4\rho,4\rho)\times\Pi_{\tilde\eta}}
&\preceq&\rho\left(\frac1{\tilde\eta^2}\mathcal{N}^{2(\mathcal{M}+7)^2}\rho\right)^{j-3}.
\end{eqnarray}
In view of the inequality
(\ref{091113}), we obtain
\begin{equation}\label{091107}
\frac1{\tilde\eta^2}\mathcal{N}^{2(\mathcal{M}+7)^2}\rho<\frac12.
\end{equation}Therefore, recalling the definition of $\Phi$ (see
(\ref{091108})) and by the inequalities (\ref{091011})
for $j_0\geq 3$, we have
\begin{equation}\label{080502}
||\Phi-id||_{p,D(s_0/4,4\rho,4\rho)} \preceq
\frac{\mathcal{N}^{98}\rho}{{\tilde\eta}^2}.
\end{equation}
Moreover, in view of the inequalities
(\ref{090701}) for $j_0\geq 3$, we obtain
\begin{equation}\label{080503}
|||D\Phi-Id|||_{p,D(s_0/4,4\rho,4\rho)}\preceq\frac
{\mathcal{N}^{98}}{{\tilde\eta}^2}.
\end{equation}Finally, let \begin{eqnarray*}{Z}=Z_{\mathcal{M}+2} , \quad
{P}=P_{\mathcal{M}+2}, \quad {Q}=Q_{\mathcal{M}+2}.
\end{eqnarray*}Basing on the estimates (\ref{091104})-(\ref{091106}), we have
\begin{eqnarray*}
|||X_{{ Z}}|||^T_{p,D(s_0/4,4\rho,4\rho)\times\Pi_{\tilde\eta}}&\leq
&\sum_{3\leq j\leq \mathcal{M}+2}|||X_{
Z_{(\mathcal{M}+2)j}}|||^T_{p,\mathcal{D}_{\mathcal{M}+2}\times\Pi_{\tilde\eta}}\\
&= &\sum_{4\leq j\leq \mathcal{M}+2}|||X_{
Z_{(\mathcal{M}+2)j}}|||^T_{p,\mathcal{D}_{\mathcal{M}+2}\times\Pi_{\tilde\eta}}\qquad \mbox{(noting $Z_{(\mathcal{M}+2)3}=0$)}\\
&\preceq&
\sum_{j=4}^{\mathcal{M}+2}\rho\left(\frac1{\tilde\eta^2}\mathcal{N}^{2(\mathcal{M}+6)^2}\rho\right)^{j-3}
\\
&\preceq&\rho\left(\frac1{\tilde\eta^2}\mathcal{N}^{2(\mathcal{M}+6)^2}\rho\right)\qquad \mbox{(by (\ref{091107}))},\\
|||X_{{ P}}|||^T_{p,D(s_0/4,4\rho,4\rho)\times\Pi_{\tilde\eta}}
&\leq& \sum_{j\geq
\mathcal{M}+3}|||X_{P_{(\mathcal{M}+2)j}}|||^T_{p,\mathcal{D}_{\mathcal{M}+2}\times\Pi_{\tilde\eta}}\\
&\preceq&\sum_{j\geq\mathcal{M}+3}
\rho\left(\frac1{\tilde\eta^2}\mathcal{N}^{2(\mathcal{M}+7)^2}\rho\right)^{j-3}\\
&\preceq&
\rho\left(\frac1{\tilde\eta^2}\mathcal{N}^{2(\mathcal{M}+7)^2}\rho\right)^{\mathcal{M}} \qquad \mbox{(by (\ref{091107}))},\\
|||X_{{
Q}}|||^T_{p,D(s_0/4,4\rho,4\rho)\times\Pi_{\tilde\eta}}&\leq&
\sum_{j\geq3}|||X_{Q_{(\mathcal{M}+2)j}}|||^T_{p,\mathcal{D}_{\mathcal{M}+2}
\times\Pi_{\tilde\eta}}\\
&\preceq& \sum_{j\geq
3}\rho\left(\frac1{\tilde\eta^2}\mathcal{N}^{2(\mathcal{M}+7)^2}\rho\right)^{j-3}\\
&\preceq&\rho\qquad\qquad\qquad\qquad\qquad\qquad \mbox{(by
(\ref{091107}))}.
\end{eqnarray*}
\end{proof}
\subsection{The measure of the set $\Pi_{\tilde\eta}$ satisfying $(\tilde\eta,\mathcal{N},\mathcal{M})$-non-resonant
conditions}\label{091406}
In this subsection, we will show
\begin{equation}\label{091301} Meas\ \Pi_{\tilde\eta} \geq(Meas\ \Pi_{\eta})(1-c\tilde\eta),
\end{equation}
where $c>0$ is a constant depending on $n$.

Firstly, we will estimate the measure of the resonant sets
$\mathcal{R}_{k\tilde l\hat l}$.
\\$\textbf{Case 1.}$\\
For $|k|\neq0$, without loss of generality, we assume
\begin{equation}\label{091206}
|k_1|=\max_{1\leq i\leq n}\{|k_1|,\dots,|k_n|\}. \end{equation} Then
\begin{eqnarray*}
&&|\partial_{\xi_1}(\langle k,\breve\omega(\xi)\rangle+\langle
\tilde
l,\tilde{\Omega}(\xi)\rangle+\langle\hat l,\hat{\Omega}(\xi)\rangle)|\\
&\geq&
|k_1||\partial_{\xi_1}\breve\omega_1(\xi)|-|\partial_{\xi_1}(\sum_{i=2}^nk_i\breve\omega_i(\xi)+\langle
\tilde l,\tilde{\Omega}(\xi)\rangle+\langle\hat l,\hat{\Omega}(\xi)\rangle)|\\
&\geq&|k_1|(1-c\eta^8\epsilon)-(\sum_{i=2}^n|k_i|+|\tilde l|+|\hat l|)c\eta^{8}\epsilon\qquad \mbox{(in view of the inequalities (\ref{091203})-(\ref{091205}))}\\
&\geq&|k_1|-(|k|+\mathcal{M}+2)c\eta^{8}\epsilon\qquad \mbox{(in
view of $|\tilde l|+|\hat l|\leq \mathcal{M}+2$)}\\&\geq&
\frac14|k_1|\qquad\qquad\qquad \mbox{(by
(\ref{091206}) and
$\mathcal{M}\leq(2c\eta^{8}\epsilon)^{-1})$}\\
&\geq&\frac14.
\end{eqnarray*}
Hence,
\begin{equation}\label{101102}
Meas\ \mathcal{R}_{k\tilde l\hat
l}\leq\frac{4\tilde\eta}{4^\mathcal{M}(|
k|+1)^{\tau}C(\mathcal{N},\tilde l)}\cdot Meas\ \Pi_{\eta}.
\end{equation}
$\textbf{Case 2.}$ \\
If $|k|=0$ and $|\tilde l|\neq0$, without loss of generality, we
assume \begin{equation}\label{091207}|\tilde {l}_1|=\max_{1\leq
j\leq \mathcal{N}}\{|\tilde l_1|,\dots,|\tilde{l}_{\mathcal{N}}|\}.
\end{equation}
Then
\begin{eqnarray*}
&&|\partial_{\xi_{n+1}}(\langle k,\breve\omega(\xi)\rangle+\langle
\tilde
l,\tilde{\Omega}(\xi)\rangle+\langle\hat l,\hat{\Omega}(\xi)\rangle)|\\
&\geq& |\tilde
l_1||\partial_{\xi_{n+1}}\tilde{\Omega}_1(\xi)|-|\partial_{\xi_{n+1}}(\langle
\tilde l,\tilde{\Omega}(\xi)\rangle+\langle\hat l,\hat{\Omega}(\xi)\rangle-\tilde l_1\tilde{\Omega}_1(\xi))|\\
&\geq&|\tilde l_1|(1-c\eta^8\epsilon)-(\sum_{i=2}^{\mathcal{N}}|\tilde l_i|+|\hat l|)c\eta^8\epsilon\qquad \mbox{(in view of the inequalities (\ref{091203})-(\ref{091205}))}\\
&\geq&|\tilde l_1|-(\mathcal{M}+2)c\eta^8\epsilon\qquad\qquad\qquad \mbox{(in view of $|\tilde l|+|\hat l|\leq\mathcal{M}+2$)}\\
&\geq& \frac14|\tilde l_1|\qquad\qquad\qquad\qquad\qquad \mbox{(in
view of $\mathcal{M}\leq(2c\eta^{8}\epsilon)^{-1}$)}\\
&\geq&\frac14.
\end{eqnarray*}
Hence,
\begin{equation}\label{101103}
Meas\ \mathcal{R}_{0\tilde l\hat
l}\leq\frac{4\tilde\eta}{4^\mathcal{M}C(\mathcal{N},\tilde l)}\cdot Meas \
\Pi_{\eta}.
\end{equation}
$\textbf{Case 3.}$ \\If $|k|=0,|\tilde l|=0$ and
$1\leq|\hat{l}|\leq2$, then it is easy to see that $ |\langle \hat
l, \hat \Omega(\xi)\rangle|$ is not small, i.e. \begin{equation}\label{101104}\mbox{the sets
$\mathcal{R}_{k\tilde l\hat l}$ are empty for $|k|=0,|\tilde l|=0$
and $1\leq|\hat{l}|\leq2$.}\end{equation}

Recall $$\mathcal{R}=\bigcup_{|k|+|\tilde l|+|\hat l|\neq 0, |\tilde
l|+|\hat l|\leq \mathcal{M}+2, |\hat l|\leq
2}\mathcal{R}_{k\tilde{l}\hat{l}}\qquad \mbox{(see
(\ref{091208}))}.$$ In view of the estimates (\ref{101102}),
(\ref{101103}) and (\ref{101104}), to estimate the measure of the
resonant set $\mathcal{R}$, we just need to count the number of the
non-resonant sets $\mathcal{R}_{k\tilde l\hat l}$. More precisely,
for fixed $k\in\mathbb{Z}$, we will show that the number of the sets
$\mathcal{R}_{k\tilde l\hat l}$ is finite.

In view of (\ref{006}) in Assumption A
\begin{equation*} |\Omega_i-\Omega_j|\geq  {
c_1}|i-j|(i+j),
\end{equation*}and the estimates
(\ref{091203})-(\ref{091205}),  then for $i\neq j,$ we have
\begin{equation}\label{091209} |\breve\Omega_i-\breve\Omega_j|\geq c_1|i-j|(i+j)-2c\eta^8\epsilon\geq \frac {c_1}2|i-j||i+j|.
\end{equation}
For $|\hat l|=1$,
\begin{equation*}
|\langle
\hat{l},\hat{\Omega}(\xi)\rangle|=|\breve\Omega_{i}(\xi)|,\qquad
\mbox{for some $i\geq \mathcal{N}+1$}.
\end{equation*}Then
in view of the inequality (\ref{091209}),
\begin{equation}\label{091210}
|\langle \hat{l},\hat{\Omega}(\xi)\rangle|\geq \frac{
c_1}{2}|i|^2\geq \frac{ c_1}{2}|i|.
\end{equation} For $|\hat l|=2$,
\begin{equation*}
|\langle
\hat{l},\hat{\Omega}(\xi)\rangle|=|\pm\breve\Omega_{i}(\xi)\pm\breve\Omega_{j}(\xi)|,\qquad
\mbox{for some $i,j\geq \mathcal{N}+1$},
\end{equation*}
then in view of the inequality (\ref{091209}) again,
\begin{equation}\label{091211}
|\langle \hat{l},\hat{\Omega}(\xi)\rangle|\geq \frac {
c_1}2|i-j|(i+j)\geq \frac {c_1}2\max\{i,j\}.
\end{equation}
If
$$\max\{i,j\}\geq \frac 4{c_1}
\left((|k|+1)(||\omega(\xi)||+1)+c_2j_0\mathcal{N}^2+1\right),$$
where the constant $c_2$ is given in (\ref{091212}), then in view of
the inequalities (\ref{091210}) and (\ref{091211}), we have
\begin{eqnarray}
\nonumber|\langle \hat{l},\hat{\Omega}(\xi)\rangle|&\geq&\frac{
c_1}{2}\left(\frac 4{ c_1}
\left((|k|+1)(||\omega(\xi)||+1)+c_2j_0\mathcal{N}^2+1\right)\right)\\
&=&2(|k|+1)(||\omega(\xi)||+1)+2c_2j_0\mathcal{N}^2+2\label{091213}.
\end{eqnarray}
Then for $|\tilde l|+|\hat l|=j_0$ and $1\leq |\hat l|\leq 2$, we
have
\begin{eqnarray}\nonumber&&\left|\langle k,\breve\omega(\xi)\rangle+\langle \tilde
l,\tilde{\Omega}(\xi)\rangle+\langle
\hat{l},\hat{\Omega}(\xi)\rangle\right|\\
\nonumber &\geq&|\langle \hat{l},\hat{\Omega}(\xi)\rangle|-|\langle
k,\breve\omega(\xi)\rangle+\langle \tilde
l,\tilde{\Omega}(\xi)\rangle|\\
\nonumber&\geq&2(|k|+1)(||\omega(\xi)||+1)+2c_2j_0\mathcal{N}^2+2-
(1+c\eta^8\epsilon)(|k|||\omega(\xi)||+c_2|\tilde l|\mathcal{N}^2)\\
\nonumber&&\mbox{(in view of the estimates (\ref{091213}), (\ref{091212}), (\ref{080103}) and (\ref{080104}))}\\
\label{091214}\nonumber&\geq&1,\end{eqnarray} which is not small.

 Then the number of the non-empty non-resonant sets
$\mathcal{R}_{k\tilde l\hat l}$ is less than
\begin{equation}\label{101201}A:=\left( 2\mathcal{N}+1\right)^{|\tilde l|}\left(
\frac 4{c_1}
\left((|k|+1)(||\omega(\xi)||+1)+c_2j_0\mathcal{N}^2+1\right)+1\right)^2.
\end{equation}
Hence,
\begin{eqnarray*}\label{11173} Meas\ \mathcal{R}&\leq& \sum_{|k|+|\tilde l|+|\hat l|\neq 0, |\tilde l|+|\hat
l|\leq \mathcal{M}+2, |\hat l|\leq 2} Meas \  \mathcal{R}_{k\tilde
l\hat l}\\
&\leq& \sum_{k\in\mathbb{Z}^n,|\tilde l|\leq
\mathcal{M}+2}\frac{4\tilde\eta A}{4^\mathcal{M}(|
k|+1)^{\tau}C(\mathcal{N},\tilde l)}\cdot Meas\ \Pi_{\eta}\\
&\leq & c\tilde\eta \cdot Meas\ \Pi_{\eta},
\end{eqnarray*}
where $ c>0$ is a constant depending on $c_1,c_2$ and $n$. This finishes the proof.

\subsection{Proof of Theorem \ref{T3}}\label{091702}
\begin{proof}Based on Theorem \ref{thm7.1}, for given positive integer
$0\leq \mathcal{M}\leq (2\eta^8\epsilon)^{-1}$ and $0<\tilde
\eta<1$, there exsists a small $\delta_0>0$ depending on $s_0,r_0,n,\tilde\eta$ and
$\mathcal{M}$, such that, for each $0<\delta<\delta_0$, $\xi\in\Pi_{\tilde\eta}$ and
the positive
integer $\mathcal{N}$ satisfying
\begin{equation}\label{091302}
\delta^{-\frac{\mathcal{M}+1}{p-1}}\leq\mathcal{N}+1<\delta^{-\frac{\mathcal{M}+1}{p-1}}+1,
\end{equation}
 there is a symplectic map
$$\Phi: D(s_0/4,4\delta,4\delta)\rightarrow D(s_0/2,5\delta,5\delta),$$where
\begin{equation*}\breve{\breve{H}}(x,y,q,\bar q;\xi):=
\breve H\circ\Phi={\breve N}(y,q,\bar q;\xi)+{Z}(y,q,\bar q;\xi)+{
P}(x,y,q,\bar q;\xi)+{ Q}(x,y,q,\bar q;\xi)
\end{equation*}
is a partial normal form of order $\mathcal{M}+2$, where
\begin{eqnarray*}
{ Z}(y,q,\bar q;\xi)&=&\sum_{4\leq 2|\alpha|+2|\beta|+2|\mu|\leq
\mathcal{M}+2,|\mu|\leq1}{
Z}^{\alpha\beta\beta\mu\mu}(\xi)y^{\alpha}\tilde{q}^{\beta}\bar{\tilde
q}^{\beta}\hat{q}^{\mu}\bar{\hat{q}}^{\mu}
\end{eqnarray*}is the integrable term depending only on variables $y$ and $I_j=|q_j|^2,j\geq1$,
\begin{eqnarray*} {
P}(x,y,q,\bar
q;\xi)&=&\sum_{2|\alpha|+|\beta|+|\gamma|+|\mu|+|\nu|\geq
\mathcal{M}+3,|\mu|+|\nu|\leq 2} { P}^{\alpha
\beta\gamma\mu\nu}(x;\xi)y^{\alpha}{\tilde q}^{\beta}{\bar{\tilde
q}}^{\gamma}{\hat q}^{\mu}\bar{\hat q}^{\nu},
\end{eqnarray*}and
\begin{equation*}
{ Q}(x,y,q,\bar q;\xi)=\sum_{|\mu|+|\nu|\geq 3}{ Q}^{\alpha
\beta\gamma\mu\nu}(x;\xi)y^{\alpha}{\tilde q}^{\beta}{\bar{\tilde
q}}^{\gamma}{\hat q}^{\mu}\bar{\hat q}^{\nu}.
\end{equation*}

Moreover, the following estimates hold:\\
(1) the symplectic map $\Phi$ satisfies
\begin{equation*}
||\Phi-id||_{p,D(s_0/4,4\delta,4\delta)}\leq
\frac{c\mathcal{N}^{98}\delta}{\tilde\eta^2},
\end{equation*}
and
\begin{equation*}
|||D\Phi-Id|||_{p,D(s_0/4,4\delta,4\delta)}\leq
\frac{c\mathcal{N}^{98}}{\tilde\eta^2};
\end{equation*}
\\
(2) the Hamiltonian vector fields $X_Z, X_P$ and $X_Q$ satisfy
\begin{equation*}
|||X_{{
Z}}|||^T_{p,D(s_0/4,4\delta,4\delta)\times\Pi_{\tilde\eta}}\leq
c\delta
\left(\frac1{\tilde\eta^2}\mathcal{N}^{2(\mathcal{M}+6)^2}\delta\right),
\end{equation*}
\begin{equation}\label{091304}
|||X_{{
P}}|||^T_{p,D(s_0/4,4\delta,4\delta)\times\Pi_{\tilde\eta}}\leq c
\delta\left(\frac1{\tilde\eta^2}\mathcal{N}^{2(\mathcal{M}+7)^2}\delta\right)^{\mathcal{M}},\qquad
 \end{equation}and
\begin{equation}\label{091305}
|||X_{{
Q}}|||^T_{p,D(s_0/4,4\delta,4\delta)\times\Pi_{\tilde\eta}}\leq c
\delta,
\end{equation}
where $c>0$ is a constant depending on $s_0,r_0,n$ and
$\mathcal{M}$.
Moreover, , the
frequency $\breve\omega(\xi)$ and $\breve\Omega(\xi)$ satisfy
$(\tilde\eta,\mathcal{N},\mathcal{M})$-non-resonant conditions.

By a direct calculation,
\begin{eqnarray}
&&\nonumber c\delta\left(\frac1{\tilde\eta^2}\mathcal{N}^{2(\mathcal{M}+7)^2}\delta\right)^\mathcal{M}\\
&=&\nonumber\delta^{\mathcal{M}+1}\left(\frac{c\mathcal{N}^{2(\mathcal{M}+7)^2\mathcal{M}}}{\tilde\eta^{2\mathcal{M}}}\right)
\\\nonumber&\leq&\delta^{\mathcal{M}+1}\left(\frac c{\tilde\eta^{2\mathcal{M}}}\delta^{-\frac{2(\mathcal{M}+7)^2\mathcal{M}(\mathcal{M}+1)}{p-1}}\right)
\qquad\mbox{(in view of the inequality (\ref{091302}))}\\
&\leq&\nonumber\delta^{\mathcal{M}+1}\left(\frac
c{\tilde\eta^{2\mathcal{M}}}\delta^{-\frac{1}{4}}\right)\qquad \mbox{(in view of $p\geq8(\mathcal{M}+7)^4+1$)}\\
&\leq&\label{091303} \delta^{\mathcal{M}+\frac58}\qquad \mbox{(by
assuming $\delta$ is very small)}.
\end{eqnarray}
In view of the inequalities (\ref{091304}) and (\ref{091303}), we
have
\begin{equation}\label{091307}
|||X_{{
P}}|||^T_{p,D(s_0/4,4\delta,4\delta)\times\Pi_{\tilde\eta}}\leq
\delta^{\mathcal{M}+\frac12}.
 \end{equation}
 Moreover, by the inequalities (\ref{081105}) and
 (\ref{091307}), we obtain
\begin{equation}\label{111901}
|||X_{{
P}}|||_{\mathcal{P}^p,D(s_0/4,4\delta,4\delta)\times\Pi_{\tilde\eta}}\leq
\delta^{\mathcal{M}+\frac12}.
 \end{equation}
On the other hand,
\begin{eqnarray}
\nonumber||\hat{z}||_1&=&\sqrt{\sum_{|j|\geq\mathcal{N}+1}|z_j|^2j^2}\\
&=&\nonumber\sqrt{\sum_{|j|\geq\mathcal{N}+1}|z_j|^2j^{2p}/j^{2(p-1)}}\\
&\leq&\nonumber\frac{||\hat{z}||_p}{(\mathcal{N}+1)^{p-1}}\\
&\leq&\label{091306}\delta^{\mathcal{M}+1}||\hat{z}||_p\qquad\qquad\qquad
\mbox{(in view of the inequality (\ref{091302}))}.
\end{eqnarray}
Then we have the following estimate
\begin{equation}\label{091308}
|||X_{{
Q}}|||_{\mathcal{P}^p,D(s_0/4,4\delta,4\delta)\times\Pi_{\tilde\eta}}\leq
\delta^{\mathcal{M}+\frac32},
 \end{equation}where we use that $Q(x,y,q,\bar
q;\xi)$ has the form of
\begin{equation*}
{ Q}(x,y,q,\bar q;\xi)=\sum_{|\mu|+|\nu|\geq 3}{ Q}^{\alpha
\beta\gamma\mu\nu}(x;\xi)y^{\alpha}{\tilde q}^{\beta}{\bar{\tilde
q}}^{\gamma}{\hat q}^{\mu}{\bar{\hat q}}^{\nu}
\end{equation*}and has $p$-tame property on the domain
$D(s_0/4,4\delta,4\delta)\times\Pi_{\tilde\eta}$, and use the
inequalities (\ref{091305}) and (\ref{091306}).

 Assume $\breve{\breve w}(t)$ is a solution
of Hamiltonian vector field $X_{\breve{\breve H}}$ with the initial
datum
$$\breve{\breve w}(0)=(\breve{\breve w}_x(0),\breve {\breve w}_y(0),\breve {\breve w}_q(0),
\breve {\breve w}_{\bar q}(0))\in D(s_0/8,4\delta,4\delta)$$
satisfying
\begin{equation*}
d_{p}(\breve{\breve w}(0),{\mathcal{T}}_{0})\leq\delta.
\end{equation*}
Noting
$\mathcal{T}_{0}=\hat{\mathbb{T}}^n\times\{y=0\}\times\{q=0\}\times\{\bar
q=0\}$ and in view of (\ref{111904}) and (\ref{111905}), we have
\begin{equation}\label{0011}
\frac1{4\delta}||\breve{\breve{w}}_y(0)||+||\breve{\breve{w}}_q(0)||_p+||\breve{\breve{w}}_{\bar
q}(0)||_p\leq\delta.
\end{equation}Firstly, we will estimate $||\breve{\breve{w}}_q(t)||_p$ and $||\breve{\breve{w}}_{\bar q}(t)||_p$. For
each $(x,y,q,\bar q)\in\mathcal{P}^p$, denote
$$\widetilde {{N}}(x,y,q,\bar q):=||q||_p^2=\sum_{j\geq1} |q_j|^2|j|^{2p}.$$
In view of the inequality (\ref{0011}), we have $ \widetilde
{{N}}(\breve{\breve{w}}(0))\leq \delta^2.$ Define
\begin{equation*}
T:=\inf\left\{|t|:\widetilde {{N}}(\breve{\breve{w}}(t))>
4\delta^2\right\},
\end{equation*}
and then $\widetilde {{N}}(\breve{\breve{w}}(t))\leq 4\delta^2$ for
all $|t|\leq T$. Specially, we have
\begin{equation}
\widetilde {{N}}(\breve{\breve{w}}(t))= 4\delta^2,\qquad \mbox{for
$t=T$ or $t=-T$}.
\end{equation}Without loss of generality, assume
\begin{equation}\label{091309}
\widetilde {{N}}(\breve{\breve{w}}(T))= 4\delta^2.
\end{equation}Now we are position to show that $T>\delta^{-\mathcal{M}}$. For each $(x,y,q,\bar
q)\in D(s_0/8,4\delta,4\delta)$,
\begin{eqnarray}
&&\nonumber\left|\left\{P,\widetilde N\right\}(x,y,q,\bar
q)\right|\\&=&\nonumber\left|\sqrt{-1}\sum_{j\geq 1}(P_{
q_j}\widetilde N_{\bar q_j}-P_{\bar q_j}\widetilde
N_{q_j})\right|\\
&\leq&\nonumber\left(||P_q||_p+||P_{\bar
q}||_p\right)||q||_p\\
&\leq&\nonumber4\delta|||X_{P}|||_{\mathcal{P}^p,D(s_0/8,4\delta,4\delta)\times\Pi_{\eta}}||q||_p\\
&\leq&\label{111902}16\delta^{\mathcal{M}+\frac52}\qquad \qquad
\mbox{(in view of the inequality (\ref{111901}) and $||q||_p\leq
4\delta$)},
\end{eqnarray}and
\begin{eqnarray}
&&\nonumber\left|\left\{Q,\widetilde N\right\}(x,y,q,\bar
q)\right|\\&=&\nonumber\left|\sqrt{-1}\sum_{j\geq 1}(Q_{
q_j}\widetilde N_{\bar q_j}-Q_{\bar q_j}\widetilde
N_{q_j})\right|\\
&\leq&\nonumber\left(||Q_q||_p+||Q_{\bar
q}||_p\right)||q||_p\\
&\leq&\nonumber4\delta|||X_{Q}|||_{\mathcal{P}^p,D(s_0/8,4\delta,4\delta)\times\Pi_{\eta}}||q||_p\\
&\leq&\label{111903}16\delta^{\mathcal{M}+\frac52}\qquad
\qquad\mbox{(in view of the inequality (\ref{091308}) and
$||q||_p\leq 4\delta$)}.
\end{eqnarray}
Since $\widetilde {{N}}(\breve{\breve{w}}(t))\leq 4\delta^2$ for all
$|t|\leq T$, then
\begin{eqnarray}
\nonumber&&\left|\widetilde {{N}}(\breve{\breve{w}}(t))-\widetilde {{N}}(\breve{\breve{w}}(0))\right|\\
&=&\nonumber\left|\int^{t}_0\left\{\breve{\breve H},\widetilde {{N}}\right\}(\breve{\breve{w}}(s))ds\right|\\
&=&\nonumber\left|\int^{t}_0\left\{P+Q, \widetilde {{N}}\right\}(\breve{\breve{w}}(s))ds\right|\qquad\qquad \mbox{(noting $\{\breve N+Z,\widetilde N\}=0$)}\\
&\leq&\nonumber\int^{t}_0\left|\left\{P+Q, \widetilde {{N}}\right\}(\breve{\breve{w}}(s))\right|ds\\
&\leq&\nonumber32\delta^{\mathcal{M}+\frac52}|t|\qquad\qquad \mbox{(based on the inequalities (\ref{111902}) and (\ref{111903}))}\\
&\leq& 32\delta^{\mathcal{M}+\frac52}T.\label{032703}
\end{eqnarray}
Assume by contradiction that $T\leq\delta^{-\mathcal{M}}$, and then
\begin{eqnarray*}
4\delta^2&=&\widetilde {{N}}(\breve{\breve{w}}(T))\\&\leq&\widetilde
{{N}}(\breve{\breve{w}}(0))+ |\widetilde
{{N}}(\breve{\breve{w}}(T))-\widetilde
{{N}}(\breve{\breve{w}}(0))|\\&\leq&\delta^2+32\delta^{\mathcal{M}+\frac52}\delta^{-\mathcal{M}}\\
&<&2\delta^2,
\end{eqnarray*}
which is impossible.

Secondly, we will estimate $||\breve{\breve{w}}_y(t)||$. For $1\leq
j\leq n$, let
$$\widetilde {{Y}}_j(x,y,q,\bar q):=y_j.$$
Then in view of the inequality (\ref{0011}), we have $ |\widetilde
{{Y}}_j(\breve{\breve{w}}(0))|\leq 4\delta^2.$ Define
\begin{equation*}
T_{j}:=\inf\left\{|t|:|\widetilde {{Y}}_j(\breve{\breve{w}}(t))|>
8\delta^2\right\},
\end{equation*}
and then $|\widetilde {{Y}}_j(\breve{\breve{w}}(t))|\leq 8\delta^2$
for all $|t|\leq T_j$. Specially, we have
\begin{equation}
|\widetilde {{Y}}_j(\breve{\breve{w}}(t))|=8\delta^2,\qquad
\mbox{for $t=T_j$ or $t=-T_j$}.
\end{equation}Without loss of generality, assume
\begin{equation}
|\widetilde {{Y}}_j(\breve{\breve{w}}(T_j))|=8\delta^2.
\end{equation}
Now, we will show that $T_j>\delta^{-\mathcal{M}}$. For each
$(x,y,q,\bar q)\in D(s_0/8,4\delta,4\delta)$,
\begin{eqnarray}
&&\nonumber\left|\left\{P,\widetilde Y_j\right\}(x,y,q,\bar
q)\right|\\&=&\nonumber\left|\sum_{1\leq i\leq n}P_{ x_i}\widetilde Y_{jy_i}\right|\\
&\leq&\nonumber||P_x||\\
&\leq&\nonumber
16\delta^2|||X_{P}|||_{\mathcal{P}^p,D(s_0/8,4\delta,4\delta)\times\Pi_{\tilde\eta}}\\&\leq&\label{111906}
16\delta^{\mathcal{M}+\frac52}\qquad\mbox{(in view of the inequality
(\ref{111901}))},
\end{eqnarray}
abd\begin{eqnarray} &&\nonumber\left|\left\{Q,\widetilde
Y_j\right\}(x,y,q,\bar
q)\right|\\&=&\nonumber\left|\sum_{1\leq i\leq n}Q_{ x_i}\widetilde Y_{jy_i}\right|\\
&\leq&\nonumber||Q_x||\\
&\leq&\nonumber
16\delta^2|||X_{Q}|||_{\mathcal{P}^p,D(s_0/8,4\delta,4\delta)\times\Pi_{\tilde\eta}}\\&\leq&\label{111907}
16\delta^{\mathcal{M}+\frac52}\qquad\mbox{(in view of the inequality
(\ref{091308}))}.
\end{eqnarray} In view of $\widetilde {{Y}}_j(\breve{\breve{w}}(t))\leq 4\delta^2$ for
all $|t|\leq T_j$, then we have
\begin{eqnarray}
\nonumber&&\left|\widetilde {{Y}}_j(\breve{\breve w}(t))-\widetilde {{Y}}_j(\breve{\breve{w}}(0))\right|\\
&=&\nonumber\left|\int^{t}_0\left\{\breve{\breve H},\widetilde {{Y}}_j\right\}(\breve{\breve{w}}(s))ds\right|\\
&=&\nonumber\left|\int^{t}_0\left\{P+Q, \widetilde {{Y}}_j\right\}(\breve{\breve{w}}(s))ds\right|\\
&\leq&\nonumber\int^{t}_0\left|\left\{P+Q, \widetilde {{Y}}_j\right\}(\breve{\breve{w}}(s))\right|ds\\
&\leq&\nonumber32\delta^{\mathcal{M}+\frac52}|t|\qquad\qquad \mbox{(by (\ref{111906}) and (\ref{111907}))}\\
&\leq& 32\delta^{\mathcal{M}+\frac52}T.\label{032702}
\end{eqnarray}
Assume by contradiction that $T_j\leq\delta^{-\mathcal{M}}$, and
then
\begin{eqnarray*}
8\delta^2&=&\widetilde {{Y}}(\breve{\breve
w}(T_j))\\&\leq&\widetilde {{Y}}_j(\breve{\breve{w}}(0))+
|\widetilde {{Y}}_j(\breve{\breve{w}}(T_j))-\widetilde
{{Y}}_j(\breve{\breve{w}}(0))|\\&\leq&4\delta^2+32\delta^{\mathcal{M}+\frac52}T\\
&<&5\delta^2,
\end{eqnarray*}
which is impossible.

Hence, for all $|t|\leq \delta^{-\mathcal{M}}$, we obtain
\begin{equation*}
\frac1{4\delta}||\breve{\breve
{w}}_y(t)||+||\breve{\breve{w}}_q(t)||_p+||\breve{\breve{w}}_{\bar
q}(t)||_p\leq2\delta,
\end{equation*}
i.e.
\begin{equation}\label{091813}
d_{p}(\breve{\breve {w}}(t),\mathcal{T}_{0})\leq2\delta.
\end{equation}
\begin{rem}\label{032901}In view of (\ref{032703}) and (\ref{032702}), we can obtain a better estimates about the change of $\tilde N(\breve{\breve w}(t))$ and $\tilde Y(\breve{\breve w}(t))$, that is, for all $|t|\leq \delta^{-\mathcal{M}}$,
\begin{equation}
|\tilde N(\breve{\breve w}(t))-\tilde N(\breve{\breve w}(0))|,|\tilde Y(\breve{\breve w}(t))-\tilde Y(\breve{\breve w}(0))|\leq 32\delta^{5/2}.
\end{equation}
\end{rem}
Based on the partial normal form (\ref{020302}) constructed in
Theorem \ref{T3}, for each $\xi\in\Pi_{\tilde\eta}$ the KAM tori
$\mathcal{T}$ of the original Hamiltonian $H(x,y,q,\bar q;\xi)$ can
be defined by $\mathcal{T}=(\Psi\circ \Phi)^{-1}{\mathcal{T}}_{0}$.
Assume $w(t)$ is a solution of original Hamiltonian vector field
$X_{{H}}$ with the initial datum $w(0)=( w_x(0), w_y(0), w_q(0),
w_{\bar q}(0))$ satisfying
\begin{equation*}
d_{p}(w(0),\mathcal{T})\leq \delta.
\end{equation*}
Then there exists $w^{*}\in\mathcal{T}$ such that
\begin{equation}\label{032701}
d_{p}(w(0),w^{*})\leq \frac98\delta.
\end{equation}
Hence,
\begin{eqnarray}\nonumber&&d_{p}(\Psi\circ\Phi \circ w(0),\mathcal{T}_{0})\\
&\leq&\nonumber d_{p}(\Psi\circ\Phi \circ w(0),\Psi\circ\Phi \circ w^{*})\qquad \nonumber\mbox{(in view of $\mathcal T=(\Psi\circ\Phi)^{-1}\mathcal{T}_0$)}\\
&\leq&\nonumber d_{p}(\Psi\circ\Phi \circ w(0),w(0))+d_{p}(w(0), w^{*})+d_{p}(w^{*},\Psi\circ\Phi \circ w^{*})\\
&\leq&\nonumber4\delta||\Psi\circ\Phi \circ w(0)-w(0)||_{\mathcal{P}^p,D(s_0/4,4\delta,4\delta)}+
\frac98\delta+4\delta||\Psi\circ\Phi \circ w^*-w^*||_{\mathcal{P}^p,D(s_0/4,4\delta,4\delta)}\\&&\nonumber \mbox{(in view of (\ref{111904}) and (\ref{032701}))}\\
&\leq&\nonumber
4\delta|||\Phi\circ\Psi-id|||_{p,D(s_0/4,4\delta,4\delta)}\left(||w(0)||_{\mathcal{P}^p,D(s_0/4,4\delta,4\delta)}+||w^*||_{\mathcal{P}^p,D(s_0/4,4\delta,4\delta)}\right)+\frac98\delta\quad
\\&\leq&
\frac43\delta,\label{032902}
\end{eqnarray}
where the last inequality is based on the inequality
(\ref{080101}) and (\ref{091812}) which imply $\Phi\circ\Psi$ is close to an identity map. Hence, based on the estimate
(\ref{091813}), Remark \ref{032901} and noting $\breve{\breve w}(t)=\Psi\circ\Phi \circ
w(t)$, we obtain \begin{equation} d_p(\Psi\circ\Phi \circ
w(t),\mathcal{T}_0)\leq \frac53\delta,\qquad \mbox{for all $|t|\leq
\delta^{-\mathcal{M}}$}.
 \end{equation} Moreover, we have
\begin{equation*}
d_{p}(w(t),(\Psi\circ\Phi)^{-1}\mathcal{T}_{0})\leq
2\delta, \qquad\mbox{for all}\ |t|\leq \delta^
{-\mathcal{M}},
\end{equation*}
which follows form the proof of (\ref{032902}).
\end{proof}\section{Proof of Theorem \ref{T4}}
\begin{proof}Consider the nonlinear
Schr$\ddot{\mbox{o}}$dinger equation (\ref{26})
\begin{equation*}
{\bf i}u_{t}=u_{xx}-M_{\xi}u+\varepsilon|u|^2u
\end{equation*}
subject to Dirichlet boundary conditions $u(t,0)=u(t,\pi)=0.$

 \textbf{Step 1.  Rewrite equation (\ref{26}) as a Hamiltonian.}\\
The eigenvalues and eigenfunctions of
$L_{M_{\xi}}=-\partial_{xx}+M_{\xi}$ with Dirichlet boundary
conditions are $\lambda_j=j^2+\xi_j$ and $\phi_j(x)=\sqrt{2/\pi}\sin
jx$, respectively. Write
\begin{equation*}
u(t,x)=\sum_{j\geq 1} w_j(t)\phi_j(x),
\end{equation*}
and then the Hamiltonian takes the form $H(w ,\bar {w})=H_0(w,\bar
w)+\varepsilon P(w,\bar w)$, where
\begin{equation*}
H_0(w,\bar w)=\sum_{j\geq 1} \lambda_jw_j\bar{w}_j
\end{equation*}
and \begin{equation*} P(w,\bar w)=\frac14\sum_{i,j,k,l\geq 1}
P_{ijkl}w_{i}\bar{w}_{j}w_{k}\bar{w}_{l},\end{equation*} with
$$P_{ijkl}=\int_0^{\pi}
\phi_{i}\phi_{j}\phi_{k}\phi_{l}dx.$$ \begin{rem} It is easy to
verify that $P_{ijkl}=0$ unless $i\pm j\pm k\pm l=0,$ for some
combination of plus and minus signs. Thus, only a codimension-one
set of coefficients is actually different from zero, and the sum
extends only over $i\pm j\pm k\pm l=0.$ In particular,
\begin{equation*}
P_{ijij}=\frac1{2\pi}(2+\delta_{ij}).
\end{equation*}
\end{rem}
In view of example 3.2 in \cite{BG}, it is proven that there exists
a constant $c_p>0$ such that
\begin{equation}\label{060701}
||\widetilde{X_{P}}(\underline{z}^{(1)},\underline{z}^{(2)},\underline{z}^{(3)})||_p\leq
c_p||(\underline{z}^3)||_{p,1},\qquad \underline{z}=(w,\bar{w}).
\end{equation}
In particular, when $p=1$, the inequality (\ref{060701}) reads
\begin{equation}\label{091901}
||\widetilde{X_{P}}(\underline{z}^{(1)},\underline{z}^{(2)},\underline{z}^{(3)})||_1\leq
c_1||(\underline{z}^3)||_{1,1}.
\end{equation}
The inequalities (\ref{060701}) and (\ref{091901}) shows that the
Hamiltonian vector filed $X_P(\underline z)$ has $p$-tame property.

\textbf{Step 2. Introduce the action-angle variables.}\\
Now we choose any finite $n$ of modes
$\phi_{j_1},\phi_{j_2},\dots,\phi_{j_n}$ as tangent direction and
the other as normal direction. Let
\begin{equation}\label{091902}
\tilde w=(w_{j_1},\dots,w_{j_n}) \qquad\mbox{and}\qquad q
=(w_j)_{j\notin\{{j_1},\dots,{j_n}\}}
\end{equation}be the tangent variable and normal variable, respectively. Then rewrite $P(w,\bar w)$ in the multiple-index as
\begin{eqnarray}\label{060103}
P(w,\bar{w})=\sum_{|\mu|+|\nu|+|\beta|=4}P^{\mu\nu\beta}\tilde
w^{\mu}\bar{\tilde w}^{\nu}z^{\beta},\qquad
\mu,\nu\in\mathbb{N}^n,\beta\in\mathbb{N}^{\mathbb{N}},
\end{eqnarray} where $z=(q,\bar{q})$ and $P^{\mu\nu\beta}=P_{ijkl}$ for some corresponding $i,j,k,l.$
In the tangent direction, introduce the action-angle variables
\begin{equation}\label{091904}
w_{j_i}=\sqrt{2(\zeta_i+y_i)}e^{ \sqrt{-1}x_i},\qquad \bar{
w}_{j_i}=\sqrt{2(\zeta_i+y_i)}e^{-\sqrt{-1}x_i},\qquad\qquad
j_i\in\{{j_1},\dots,{j_n}\},
\end{equation}where $\zeta_i\in[1,2]$ is the initial datum and will be considered as a constant. The symplectic
structure is $dy\wedge dx+\sqrt{-1}dq\wedge d\bar{q}$, where
$y=(y_1,\dots,y_n)$ and $x=(x_1,\dots,x_n)$.
 Hence, (\ref{060103})
is turned into
\begin{eqnarray}
P(w,\bar w)&=&\nonumber P(x,y,
z)\\&=&\nonumber\sum_{|\mu|+|\nu|+|\beta|=4}2^{\frac12(|\mu|+|\nu|)}P^{\mu\nu\beta}\sqrt{(\zeta+y)^{\mu}}e^{\sqrt{-1}\langle\mu,x\rangle}
\sqrt{(\zeta+y)^{\nu}}e^{-\sqrt{-1}\langle\nu,x\rangle}z^{\beta}\\
\nonumber
&=&\sum_{|\mu|+|\nu|+|\beta|=4}(2^{\frac12(|\mu|+|\nu|)}P^{\mu\nu\beta}
e^{\sqrt{-1}\langle\mu-\nu,x\rangle})\sqrt{(\zeta+y)^{\mu+\nu}}
z^{\beta}\nonumber\\
&=&\sum_{|\beta|\leq4}P^{\beta}(x,y) z^{\beta}\label{121001},
\end{eqnarray}
where \begin{eqnarray*}\zeta=(\zeta_1,\dots,\zeta_n)\qquad
\mbox{and}\qquad P^{\beta}(x,y)=\sum_{|\mu|+|\nu|=4-|\beta|}
2^{\frac12(|\mu|+|\nu|)}P^{\mu\nu\beta}e^{\sqrt{-1}\langle\mu-\nu,x\rangle}\sqrt{(\zeta+y)^{\mu+\nu}}.
\end{eqnarray*}

\textbf{Step 3. Show that the Hamiltonian vector field $X_P$ has $p$-tame property.}\\
In this step, we will show that the Hamiltonian vector field
$X_P(x,y,z)$ has $p$-tame property after introducing the
action-angle variables. Assume
$P(x,y,z)=\sum_{|\beta|\leq4}P^{\beta}(x,y) z^{\beta}$ is defined on
the domain $(x,y,z)\in D(s_0,r_0,r_0)$ for some
$0<s_0,r_0\leq 1$. \\
\textbf{Step 3.1. Estimate $|||P_x|||_{D(s_0,r_0,r_0)\times\Pi}.$}\\
For $1\leq i\leq n$,
\begin{eqnarray}
&&\nonumber P_{x_i}(x,y,z)\\&=&\nonumber
P_{w_{j_i}}(x,y,z)\cdot(w_{j_i})_{x_i}+P_{\bar
w_{j_i}}(x,y,z)\cdot(\bar
w_{j_i})_{x_i}\\
&=&\label{091905}P_{w_{j_i}}(x,y,z)(\sqrt{-1}\sqrt{\zeta_i+y_i}e^{\sqrt{-1}x_i})+P_{\bar
w_{j_i}}(x,y,z)(-\sqrt{-1}\sqrt{\zeta_i+y_i}e^{-\sqrt{-1}x_i}).
\end{eqnarray}
In view of (\ref{060103}) and (\ref{121001}), we have
\begin{eqnarray}\label{091903}
P_{w_{j_i}}(x,y,z)=P_{w_{j_i}}(w,\bar{w})=\sum_{|\mu|+|\nu|+|\beta|=4}\mu_{i}P^{\mu\nu\beta}\tilde
w^{\mu-1_i}\bar{\tilde w}^{\nu}z^{\beta},
\end{eqnarray}where
\begin{equation*}
\mu-1_i=(\mu_1,\dots,\mu_i-1,\dots,\mu_n).
\end{equation*}Denote
\begin{equation*}
P_{w_{j_i}}^\beta(\tilde w,\bar{\tilde
w})=\sum_{|\mu|+|\nu|=4-|\beta|}\mu_{i}P^{\mu\nu\beta}\tilde
w^{\mu-1_i}\bar{\tilde w}^{\nu}.
\end{equation*}
Then in view of (\ref{091903}), we have
\begin{equation}\label{121002}
P_{w_{j_i}}(x,y,z)=P_{w_{j_i}}(w,\bar{w})=P_{w_{j_i}}^\beta(\tilde
w,\bar{\tilde w})z^{\beta}.
\end{equation}
Furthermore, let
\begin{equation*}P_{\tilde{w}_{j_i}}^\beta(\tilde w,\bar{\tilde w})=\sum_{|\mu|+|\nu|=4-|\beta|}\nu_{i}P^{\mu\nu\beta}\tilde
w^{\mu}\bar{\tilde w}^{\nu-1_i},
\end{equation*}where $\nu-1_i=(\nu_1,\dots,\nu_i-1,\dots,\nu_n)$, and then
\begin{equation}\label{121003}
P_{\tilde{w}_{j_i}}(x,y,z)=P_{\tilde{w}_{j_i}}(w,\bar{w})=P_{\tilde{w}_{j_i}}^\beta(\tilde
w,\bar{\tilde w})z^{\beta}.
\end{equation}
Denote
\begin{equation}\label{091908}
P_{x_i}^{\beta}(x,y)=P_{w_{j_i}}^\beta(\tilde w,\bar{\tilde
w})(\sqrt{-1}\sqrt{\zeta_i+y_i}e^{\sqrt{-1}x_i})+P_{\tilde{w}_{j_i}}^\beta(\tilde
w,\bar{\tilde w})(-\sqrt{-1}\sqrt{\zeta_i+y_i}e^{-\sqrt{-1}x_i}).
\end{equation}In view of the formulas (\ref{091905}) and (\ref{121002})-(\ref{091908}), we have
\begin{equation}\label{121004}P_{x_i}(x,y,z)=\sum_{|\beta|\leq
4}P_{x_i}^{\beta}(x,y)z^{\beta},\qquad
\beta\in\mathbb{N}^{\mathbb{N}}.
\end{equation}In view of the inequality (\ref{091901}) and for
$1\leq i\leq n$, we have
\begin{eqnarray}\label{0008}
|\widetilde{P_{w_{j_i}}}(\underline{z}^{(1)},\underline{z}^{(2)},\underline{z}^{(3)})|,|\widetilde{P_{\bar
w_{j_i}}}(\underline{z}^{(1)},\underline{z}^{(2)},\underline{z}^{(3)})|\leq
c_1(||\underline{z}^{(1)}||_1||\underline{z}^{(2)}||_1||\underline{z}^{(3)}||_1).
\end{eqnarray}In view of (\ref{121004}), based on the inequality (\ref{0008}) and the definition of ${\lfloor \cdot\rceil}_{D(s_0,r_0)\times\Pi}$ (see (\ref{092703}) for the details), we obtain
 \begin{eqnarray}\label{121005}
&&\nonumber\widetilde{\lfloor{P_{x_i}}\rceil}_{D(s_0,r_0)\times\Pi}({z}^{(1)},{z}^{(2)},{z}^{(3)})\\
&\leq& c(\sum_{1\leq i\leq 3}||z^{(i)}||_1+\sum_{1\leq i\neq j\leq
3}||z^{(i)}||_1||z^{(j)}||_1+||z^{(1)}||_1||z^{(2)}||_1||z^{(3)}||_1),\label{121005}
 \end{eqnarray}where $c>0$ is a constant depending on $s_0 ,r_0$ and $n$.
Based on the inequality (\ref{121005}) for $1\leq i\leq n$ and
Definition \ref{071802},
 \begin{equation}\label{091909}
 |||P_{x}|||_{D(s_0,r_0,r_0)\times\Pi}<\infty.
 \end{equation}
Following the proof of the inequality (\ref{091909}), we obtain
\begin{eqnarray}\label{060202} |||P_y|||_{D(s_0,r_0,r_0)\times\Pi}<\infty.
\end{eqnarray}Now, we would like to show
\begin{equation*}
|||P_{ z}|||^{T}_{p,D(s_0,r_0,r_0)\times\Pi}<\infty.
\end{equation*}
In view of (\ref{060103}),
\begin{eqnarray}
P_{z}(x,y,z)&=&P_z(w,\bar{w})\nonumber\\
&=&\sum_{|\mu|+|\nu|+|\beta|=4}\beta P^{\mu\nu\beta}\tilde
w^{\mu}\bar{\tilde w}^{\nu}z^{\beta-1}.
\end{eqnarray}Let
\begin{equation}
P_{z}^{\beta-1}(x,y)=P_{z}^{\beta-1}(\tilde w,\bar{\tilde
w})=\sum_{|\mu|+|\nu|=4-|\beta|}\beta P^{\mu\nu\beta}\tilde
w^{\mu}\bar{\tilde w}^{\nu},
\end{equation}
and then
\begin{equation}
P_{z}(x,y,z)=\sum_{|\beta|\leq 4}P_{z}^{\beta-1}(x,y)z^{\beta-1}.
\end{equation}
Hence, we obtain
\begin{eqnarray}\label{091910}
||\widetilde{\lfloor P_z\rceil}_{D(s_0,r_0)\times\Pi}||_p\leq
c(\sum_{1\leq i\leq 3}||{z}^{(i)}||_p+\sum_{1\leq i\neq j\leq
3}||{z}^{(i)}||_1 ||{z}^{(j)}||_p+||({z}^{3})||_{p,1}),
\end{eqnarray}
where $c>0$ is a constant depending on $s_0,r_0,n$ and $p$, and the
above inequality is based on the inequality (\ref{060701}) and the
definition of ${\lfloor \cdot\rceil}_{D(s_0,r_0)\times\Pi}$ (see
(\ref{092703}) for the details), and noting that
$\underline{z}=(\tilde w,q,\bar{\tilde w},\bar q)$ and $z=(q,\bar
q)$. In particular, when $p=1$, the inequality (\ref{091910}) reads
\begin{eqnarray}\label{091911}
||\widetilde{\lfloor P_z\rceil}_{D(s_0,r_0)\times\Pi}||_1\leq \breve
c(\sum_{1\leq i\leq 3}||{z}^{(i)}||_1+\sum_{1\leq i\neq j\leq
3}||{z}^{(i)}||_1 ||{z}^{(j)}||_1+||({z}^{3})||_{1,1}),
\end{eqnarray}where $\breve c>0$ is a constant depending on
$s_0,r_0$ and $n$. Based on the inequalities (\ref{091910}) and
(\ref{091911}), we obtain
\begin{equation}\label{060203}
|||P_{z}|||^{T}_{p,D(s_0,r_0,r_0)\times \Pi}<\infty.
\end{equation}
In view of the inequalities (\ref{091909}), (\ref{060202}) and
(\ref{060203}), we get
\begin{equation}\label{091912}
|||X_P|||^{T}_{p,D(s_0,r_0,r_0)\times\Pi}<\infty.
\end{equation}
Finally, we obtain a Hamiltonian $H(x,y,q,\bar q;\xi)$ having the
following form
\begin{equation}H(x,y,q,\bar q;\xi)=N(x,y,q,\bar q;\xi)+R(x,y,q,\bar
q;\xi),\end{equation} where
\begin{equation}
N(x,y,q,\bar q;\xi)=H_{0}(w,\bar w)=\sum_{1\leq j\leq n}
\omega_j(\xi)y_j+\sum_{j\geq 1}\Omega_j(\xi)q_j\bar q_j,
\end{equation}with the tangent frequency \begin{equation}\label{001}
\omega({\xi})=(\omega_1(\xi),\dots,\omega_n(\xi))=(j_1^2+\xi_{j_1},\dots,j_{n}^2+\xi_{j_n}),
\end{equation}and the normal frequency
\begin{equation}\label{002}
\Omega(\xi)=(\Omega_1(\xi),\Omega_2(\xi),\dots,),\qquad
\Omega_j(\xi)=j^2+\xi_j,\quad j\notin\{j_1,\dots,j_n\},
\end{equation}
and the perturbation
\begin{equation}\label{003}
R(x,y,q,\bar q;\xi)=\varepsilon P(x,y,q,\bar q)
\end{equation}
is independent of parameters $\xi$.

In view of the formulas (\ref{001}) and (\ref{002}), it is easy to
show that Assumption A and Assumption B in Theorem \ref{T1} hold.
Basing on the inequality (\ref{091912}) and noting $R(x,y,q,\bar
q;\xi)=\varepsilon P(x,y,q,\bar q)$ (see (\ref{003})), we obtain
$|||X_R|||^T_{p,D(s_0,r_0,r_0)\times\Pi}$ satisfies the small
assumption. Hence, all assumptions in Theorem \ref{T1} hold. According to Theorem \ref{T1}, there exists a subset
$\Pi_{\eta}\subset\Pi$ with the estimate
\begin{equation*}
\mbox{Meas}\ \Pi_{\eta}\geq(\mbox{Meas}\ \Pi)(1-O(\eta)).
\end{equation*}For each $\xi\in\Pi_{\eta}$, there is a symplectic map
$$\Psi: D(s_0/2,r_0/2,r_0/2)\rightarrow D(s_0,r_0,r_0),$$ such that $H(x,y,q,\bar q;\xi)$ can be transformed into a normal form of order 2 with the following form
\begin{equation*}\breve H(x,y,q,\bar q;\xi):=H\circ\Psi=
\breve N(y,q,\bar q;\xi)+ \breve R(x,y,q,\bar q;\xi),
\end{equation*}where
\begin{equation*}
\breve N(y,q,\bar
q;\xi)=\sum_{j=1}^{n}\breve{\omega}_j(\xi)y_j+\sum_{j\geq1}\breve{\Omega}_j(\xi)q_j\bar
q_j
\end{equation*}
and
\begin{equation*}
\breve R(x,y,q,\bar
q;\xi)=\sum_{\alpha\in\mathbb{N}^n,\beta,\gamma\in\mathbb{N}^{\mathbb{N}},2|\alpha|+|\beta|+|\gamma|\geq
3}\breve
R^{\alpha\beta\gamma}(x;\xi)y^{\alpha}q^{\beta}\bar{q}^{\gamma}.
\end{equation*}

Let $\delta$ be given in the statement of Theorem \ref{T3} and $\mathcal{N}$ be given in (\ref{091302}). Take $\tilde\eta$ satisfying $\delta^{1/6}<\tilde
\eta<{\mathcal{N}}^{-2}$ and
$\tilde\Pi=\Pi_{\tilde\eta}$. Based on Theorem \ref{T3}, for each $\xi\in\tilde\Pi$, equation (\ref{26}) possesses a KAM
torus $\mathcal{T}_{\xi}$ in Sobolev space $H^p_0([0,\pi])$ which is
 defined by
\begin{equation}
H_0^p[0,\pi]=\left\{v(x)=\sum_{j\geq
1}\widehat{v}(j)\phi_j(x)\mid\sum_{j\geq1}|\widehat{v}(j)|^2j^{2p}<\infty\right\},
\end{equation}where $\phi_j(x)=\sqrt{2/\pi}\sin
jx$ and $$ \widehat{v}(j)=\sqrt{\frac2\pi}\int_{0}^{\pi}v(x)\sin
jxdx$$ is the $j$-th Fourier coefficient of $v(x)$.

 Furthermore,
for any solution
$u(t,x)$ of equation (\ref{26}) with the
initial datum satisfying \color{black}
$${d}_{H^p_0[0,\pi]}(u(0,x),\mathcal{T}_{\xi})\leq \delta,$$
then
\begin{equation*}{d}_{H^p_0[0,\pi]}(u(t,x),\mathcal{T}_{\xi})\leq
2\delta,\qquad \mbox{for all} \ |t|\leq
{\delta}^{-\mathcal{M}}.
\end{equation*}
\end{proof}
\section{Appendix: technical lemmas}

\begin{lem}
Consider two functions $U(x;\xi)$ and $V(x;\xi)$ defined on the
domain $D(s)\times\Pi$, which are analytic about the variable $x\in
D(s)$ and $C^1$-smooth about the parameter $\xi\in\Pi$, then the
following inequality holds
\begin{equation*}
||UV||_{D(s)\times\Pi}\leq ||U||_{D(s)\times\Pi}
||V||_{D(s)\times\Pi}.
\end{equation*}
\end{lem}
\begin{proof}
Let $W(x;\xi)=U(x;\xi)V(x;\xi)$ with its Fourier series
$W(x;\xi)=\sum_{k\in\mathbb{Z}^n}\widehat{W}(k;\xi)e^{\sqrt{-1}\langle
k,x\rangle}.$ By a direct calculation,
\begin{equation*}
\widehat{W}(k;\xi)=\sum_{k_1+k_2=k}\widehat{U}(k_1;\xi)\widehat{V}(k_2;\xi),
\end{equation*}
and
\begin{equation*}
\partial_{\xi_j}\widehat{W}(k;\xi)=\sum_{k_1+k_2=k}
\left(\partial_{\xi_j}\widehat{U}(k_1;\xi)\widehat{V}(k_2;\xi)
+\widehat{U}(k_1;\xi)\partial_{\xi_j}\widehat{V}(k_2;\xi)\right),
\end{equation*}
where $\widehat{U}(k;\xi)$ and $\widehat{V}(k;\xi)$ are $k$-th
Fourier coefficients of $U(x;\xi)$ and $V(x;\xi)$, respectively.
Hence,
\begin{eqnarray*}
&&\sum_{k\in\mathbb{Z}^n}|\widehat{W}(k;\xi)|e^{|k|s}\\
&=&\sum_{k\in\mathbb{Z}^n}|\sum_{k_1+k_2=k}\widehat{U}(k_1;\xi)
\widehat{V}(k_2;\xi)|e^{|k_1+k_2|s}\\
&\leq&\sum_{k\in\mathbb{Z}^n}\sum_{k_1+k_2=k}|\widehat{U}(k_1;\xi)|
|\widehat{V}(k_2;\xi)|e^{(|k_1|+|k_2|)s}\\
&=&\left(\sum_{k_1\in\mathbb{Z}^n}|\widehat{U}(k_1;\xi)|e^{|k_1|s}\right)\left
(\sum_{k_2\in\mathbb{Z}^n}|\widehat{V}(k_2;\xi)|e^{|k_2|s}\right)
\end{eqnarray*}and
\begin{eqnarray*}
&&\sum_{k\in\mathbb{Z}^n}|\partial_{\xi_j}\widehat{W}(k;\xi)|e^{|k|s}\\
&=&\sum_{k\in\mathbb{Z}^n}|\sum_{k_1+k_2=k}
\partial_{\xi_j}\widehat{U}(k_1;\xi)\widehat{V}(k_2;\xi)+\widehat{U}(k_1;\xi)
\partial_{\xi_j}\widehat{V}(k_2;\xi))|e^{|k_1+k_2|s}\\
&\leq&\sum_{k\in\mathbb{Z}^n}\sum_{k_1+k_2=k}\left(
|\partial_{\xi_j}\widehat{U}(k_1;\xi)||
\widehat{V}(k_2;\xi)|+|\widehat{U}(k_1;\xi)||\partial_{\xi_j}
\widehat{V}(k_2;\xi)|\right)e^{(|k_1|+|k_2|)s}\\
&=&\left(\sum_{k_1\in\mathbb{Z}^n}
|\partial_{\xi_j}\widehat{U}(k_1;\xi)|e^{|k_1|s}\right)\left(\sum_{k_2\in\mathbb{Z}^n}|
\widehat{V}(k_2;\xi)|e^{|k_2|s}\right)\\&&+\left(\sum_{k_1\in\mathbb{Z}^n}
|\widehat{U}(k_1;\xi)|e^{|k_1|s}\right)\left(\sum_{k_2\in\mathbb{Z}^n}|
\partial_{\xi_j}\widehat{V}(k_2;\xi)|e^{|k_2|s}\right).
\end{eqnarray*}
Therefore,
\begin{eqnarray*}
&&\sum_{k\in\mathbb{Z}^n}\left(
|\widehat{W}(k;\xi)|+|\partial_{\xi_j}\widehat{W}(k;\xi)|\right)e^{|k|s}
\\
&\leq&\left(\sum_{k_1\in\mathbb{Z}^n}|\widehat{U}(k_1;\xi)|e^{|k_1|s}\right)
\left(\sum_{k_2\in\mathbb{Z}^n}|\widehat{V}(k_2;\xi)|e^{|k_2|s}\right)\\
&&+\left(\sum_{k_1\in\mathbb{Z}^n}
|\partial_{\xi_j}\widehat{U}(k_1;\xi)|e^{|k_1|s}\right)\left(\sum_{k_2\in\mathbb{Z}^n}|
\widehat{V}(k_2;\xi)|e^{|k_2|s}\right)\\&&+\left(\sum_{k_1\in\mathbb{Z}^n}
|\widehat{U}(k_1;\xi)|e^{|k_1|s}\right)\left(\sum_{k_2\in\mathbb{Z}^n}|
\partial_{\xi_j}\widehat{V}(k_2;\xi)|e^{|k_2|s}\right)\\
&\leq&\left(\sum_{k\in\mathbb{Z}^n}(|\widehat{U}(k;\xi)|+|\partial_{\xi_j}\widehat{U}
(k;\xi)|)e^{|k|s}\right)\left(\sum_{k\in\mathbb{Z}^n}(|\widehat{V}(k;\xi)|+
|\partial_{\xi_j}\widehat{V}(k;\xi)|)e^{|k|s}\right).\end{eqnarray*}
Thus,
\begin{eqnarray*}
&&||W||_{D(s)\times\Pi}\\&=&\sup_{\xi\in\Pi,j\geq1}\sum_{k\in\mathbb{Z}^n}\left(
|\widehat{W}(k;\xi)|+|\partial_{\xi_j}\widehat{W}(k;\xi)|\right)e^{|k|s}\\
&\leq&\sup_{\xi\in\Pi,j\geq1}\left(\left(\sum_{k\in\mathbb{Z}^n}(|\widehat{U}(k;\xi)|+|\partial_{\xi_j}\widehat{U}
(k;\xi)|)e^{|k|s}\right)\left(\sum_{k\in\mathbb{Z}^n}(|\widehat{V}(k;\xi)|+
|\partial_{\xi_j}\widehat{V}(k;\xi)|)e^{|k|s}\right)\right)\\
&\leq&\sup_{\xi\in\Pi,j\geq1}\left(\sum_{k\in\mathbb{Z}^n}(|\widehat{U}(k;\xi)|+|\partial_{\xi_j}\widehat{U}
(k;\xi)|)e^{|k|s}\right)\sup_{\xi\in\Pi,j\geq1}\left(\sum_{k\in\mathbb{Z}^n}(|\widehat{V}(k;\xi)|+
|\partial_{\xi_j}\widehat{V}(k;\xi)|)e^{|k|s}\right)\\
&=&||U||_{D(s)\times\Pi}||V||_{D(s)\times\Pi}.
\end{eqnarray*}
\end{proof}\begin{lem}\label{051701}
Consider two functions
$U(x,y;\xi)=\sum_{\alpha\in\mathbb{N}^n}U^{\alpha}(x;\xi)y^{\alpha}$
and
$V(x,y;\xi)=\sum_{\alpha\in\mathbb{N}^n}V^{\alpha}(x;\xi)y^{\alpha}$
defined on the domain $D(s,r)\times\Pi$, which are analytic about
the variable $(x,y)\in D(s,r)$ and $C^1$-smooth about the parameter
$\xi\in\Pi$, then the following inequality holds
\begin{equation*}
||UV||_{D(s,r)\times\Pi}\leq ||U||_{D(s,r)\times\Pi}
||V||_{D(s,r)\times\Pi}.
\end{equation*}
\end{lem}
\begin{proof}
Let $W(x,y;\xi)=U(x,y;\xi)V(x,y;\xi)$ with its Taylor series about
the variable $y$
$$W(x,y;\xi)=\sum_{\beta\in\mathbb{N}^n}W^{\beta}(x;\xi)y^{\beta}.$$
By a direct calculation,
\begin{equation*}
W^{\beta}(x;\xi)=\sum_{\alpha+\alpha'=\beta}U^{\alpha}(x;\xi)V^{\alpha'}(x;\xi).
\end{equation*}
Hence, by the definition of $||\cdot||_{D(s,r)\times\Pi}$,
\begin{eqnarray*}
||W||_{D(s,r)\times\Pi}&=&\sum_{\beta\in\mathbb{N}^n}||W^{\beta}||_{D(s)\times\Pi}
r^{2|\beta|}\\
&=&\sum_{\beta\in\mathbb{N}^n}||\sum_{\alpha+\alpha'=\beta}
U^{\alpha}V^{\alpha'}||_{D(s)\times\Pi}r^{2|\alpha+\alpha'|}
\\
&\leq&\sum_{\beta\in\mathbb{N}^n}\sum_{\alpha+\alpha'=\beta}
||U^{\alpha}V^{\alpha'}||_{D(s)\times\Pi}r^{2|\alpha+\alpha'|}\\
&\leq&\sum_{\beta\in\mathbb{N}^n}\sum_{\alpha+\alpha'=\beta}
||U^{\alpha}||_{D(s)\times\Pi}||V^{\alpha'}||_{D(s)\times\Pi}r^{2|\alpha+\alpha'|}\\
&=&\left(\sum_{\alpha\in\mathbb{N}^n}||U^{\alpha}||_{D(s)\times\Pi}r^{2|\alpha|}\right)
\left(\sum_{\alpha'\in\mathbb{N}^n}||V^{\alpha'}||_{D(s)\times\Pi}r^{2|\alpha'|}\right)\\
&=&||U||_{D(s,r)\times\Pi}||V||_{D(s,r)\times\Pi}.
\end{eqnarray*}
\end{proof}
Now we will give two generalized Cauchy estimates:
\begin{lem}\label{050402} Consider a function
$W(x;\xi)=\sum_{k\in\mathbb{Z}^n}\widehat{W}(k;\xi)e^{\sqrt{-1}\langle
k,x\rangle}$ defined on the domain $D(s)\times\Pi$, which is
analytic about the variable $x\in D(s)$ and $C^1$-smooth about the
parameter $\xi\in\Pi$, then the following generalized Cauchy
estimate holds
\begin{equation*}
||W_x||_{D(s-\sigma)\times\Pi}\leq \frac
1{e\sigma}||W||_{D(s)\times\Pi},
\end{equation*}
where $0<\sigma<s$.
\end{lem}
\begin{proof}
By a direct calculation,
\begin{equation*}
W_{x_j}(x;\xi)=\sum_{k\in\mathbb{Z}^n}\sqrt{-1}k_j\widehat{W}(k;\xi)e^{\sqrt{-1}\langle
k,x\rangle}.
\end{equation*}
Hence,
\begin{eqnarray*}
&&\sum_{k\in\mathbb{Z}^n}|k_j|(|\widehat{W}(k;\xi)|+|\partial_{\xi_j}\widehat{W}(k;\xi)|)e^{|k|(s-\sigma)}\\
&\leq&\sum_{k\in\mathbb{Z}^n}(|\widehat{W}(k;\xi)|+|\partial_{\xi_j}\widehat{W}(k;\xi)|)e^{|k|s}|k|e^{-|k|\sigma}\\
&\leq&\left(\sup_{k\in\mathbb{Z}^n}|k|e^{-|k|\sigma}\right)\left(\sum_{k\in\mathbb{Z}^n}(|\widehat{W}(k;\xi)|+|\partial_{\xi_j}\widehat{W}(k;\xi)|)e^{|k|s}\right)
\\
&\leq&\frac
1{e\sigma}\left(\sum_{k\in\mathbb{Z}^n}(|\widehat{W}(k;\xi)|+|\partial_{\xi_j}\widehat{W}(k;\xi)|)e^{|k|s}\right).
\end{eqnarray*}
Hence,
\begin{eqnarray*}
&&||W_{x_j}||_{D(s-\sigma)\times\Pi}\\
&=&\sup_{\xi\in\Pi,j\geq1}\left(
\sum_{k\in\mathbb{Z}^n}|k_j|(|\widehat{W}(k;\xi)|
+|\partial_{\xi_j}\widehat{W}(k;\xi)|)e^{|k|(s-\sigma)}\right)\\
&\leq&\sup_{\xi\in\Pi,j\geq1}\left(\frac
1{e\sigma}\left(\sum_{k\in\mathbb{Z}^n}
(|\widehat{W}(k;\xi)|+|\partial_{\xi_j}\widehat{W}(k;\xi)|)e^{|k|s}\right)\right)\\
&=&\frac 1{e\sigma}||W||_{D(s)\times\Pi}.
\end{eqnarray*}
Moreover,
\begin{equation*}
||W_x||_{D(s-\sigma)\times\Pi}= \sup_{1\leq j\leq
n}||W_{x_j}||_{D(s-\sigma)\times\Pi}\leq \frac
{1}{e\sigma}||W||_{D(s)\times\Pi}.
\end{equation*}
\end{proof}\begin{lem}\label{050403}
Consider a function
$W(x,y;\xi)=\sum_{\alpha\in\mathbb{N}^n}W^{\alpha}(x;\xi)y^{\alpha}$
defined on the domain $D(s,r)\times\Pi$, which is analytic about the
variable $(x,y)\in D(s,r)$ and $C^1$-smooth about the parameter
$\xi\in\Pi$, then the following generalized Cauchy estimate holds
\begin{equation}\label{071901}
||W_x||_{D(s-\sigma,r)\times\Pi}\leq \frac
1{e\sigma}||W||_{D(s,r)\times\Pi},
\end{equation}and
\begin{equation}\label{071902}
||W_y||_{D(s,r-\sigma')\times\Pi}\leq \frac
1{r\sigma'}||W||_{D(s,r)\times\Pi},
\end{equation}
where $0<\sigma<s$ and $0<\sigma'<r/2$.
\end{lem}
\begin{proof}The inequality (\ref{071901}) can be obtained directly by Lemma
\ref{050402}. Now, we will prove the inequality (\ref{071902}). By a
direct calculation,
\begin{equation*}
W_{y_j}(x,y;\xi)=\sum_{\alpha\in\mathbb{N}^n}\alpha_j
W^{\alpha}(x;\xi)y^{\alpha-1_j},
\end{equation*}
where
\begin{equation*}
\alpha-1_j=(\alpha_1,\dots,\alpha_{j-1},\alpha_j-1,\alpha_{j+1},\dots,\alpha_n).
\end{equation*}
 Hence,
\begin{eqnarray*}
&&||W_{y_j}||_{D(s,r-\sigma')\times\Pi}\\
&=&\sum_{\alpha\in\mathbb{N}^n}||\alpha_j W^{\alpha}||_{D(s)\times\Pi}(r-\sigma')^{2(|\alpha|-1)}\\
&\leq&\sum_{\alpha\in\mathbb{N}^n}|\alpha|||W^{\alpha}||_{D(s)\times\Pi}(r-\sigma')^{2(|\alpha|-1)}\\
&=&\frac1{2(r-\sigma')}\sum_{\alpha\in\mathbb{N}^n}||W^{\alpha}||_{D(s)\times\Pi}2|\alpha|(r-\sigma')^{2|\alpha|-1}\\
&\leq&\frac1{2(r-\sigma')}\sum_{\alpha\in\mathbb{N}^n}||W^{\alpha}||_{D(s)\times\Pi}r^{2|\alpha|}/\sigma'\qquad
 \mbox{(using the inequality }\ k(r-\sigma')^{k-1}\leq\frac{r^{k}}{\sigma'})\\
&\leq&\frac1{r\sigma'}\sum_{\alpha\in\mathbb{N}^n}||W^{\alpha}||_{D(s)\times\Pi}r^{2|\alpha|}\qquad
(\mbox{based on $0<\sigma'<r/2$})
 \\&=&\frac1{r\sigma'}||W||_{D(s,r)\times\Pi}.
\end{eqnarray*}
Therefore,
\begin{equation*}||W_y||_{D(s,r-\sigma')\times\Pi}=
\sup_{1\leq j\leq n}||W_{y_j}||_{D(s,r-\sigma')\times\Pi} \leq
\frac1{r\sigma'}||W||_{D(s,r)\times\Pi}.
\end{equation*}
\end{proof}
Assume $F(z)$ is a homogeneous polynomial with the corresponding
multi-linear form $\widetilde F(z^{(1)},\dots,z^{(h)})$. Then the
following lemma holds:
\begin{lem}\label{042201}Consider two homogeneous polynomials $F(z)=
\sum_{|\beta|=h_1}F^{\beta}z^{\beta}$ and
$G(z)=\sum_{|\beta'|=h_2}G^{\beta'}z^{\beta'}$. Then
\begin{equation*}
\widetilde {F
G}(z^{(1)},\dots,z^{(h)})=\frac1{h!}\sum_{\tau_h}\widetilde
F(z^{(\tau_h(1))},\dots,z^{(\tau_h(h_1))})\widetilde
G(z^{(\tau_h(h_1+1))},\dots,z^{(\tau_h(h))}),
\end{equation*}
where $h=h_1+h_2$ and $\tau_h$ is an $h$-permutation.
\end{lem}
\begin{proof}It suffices to consider $F(z)$ and $G(z)$ as two
monomials with the following forms
$$F(z)=F^{\beta}z^{\beta}=F^{\beta}z_{j_1}\cdots z_{j_{h_1}}\quad \mbox{and}\quad
G(z)=G^{\beta'}z^{\beta'}=G^{\beta'}z_{i_1}\cdots z_{i_{h_2}}.$$
Then
\begin{equation*}
W(z)=F(z)G(z)=F^{\beta}G^{\beta'}z_{j_1}\cdots
z_{j_{h_1}}z_{i_1}\cdots
z_{i_{h_2}}=F^{\beta}G^{\beta'}z_{k_1}\cdots z_{k_{h}},
\end{equation*}
where $(k_1,\dots,k_h)=(j_1,\dots,j_{h_1},i_{1},\dots,i_{h_2})$.
 Hence
\begin{eqnarray*}
&&\frac1{h!}\sum_{\tau_h}\widetilde
F(z^{(\tau_h(1))},\dots,z^{(\tau_h(h_1))})\widetilde
G(z^{(\tau_h(h_1+1))},\dots,z^{(\tau_h(h))})\\
&=&\frac1{h!}\sum_{\tau_h}\left(\frac{F^{\beta}}{h_1!}
\sum_{\tau_{h_1}}z_{j_1}^{(\tau_{h_1}\circ\tau_{h}(1))}\cdots
z_{j_{h_1}}^{(\tau_{h_1}\circ\tau_{h}(h_1))}\right)\left
(\sum_{\tau_{h_2}}\frac{G^{\beta'}}{h_2!}z_{i_1}^{(\tau_{h_2}\circ\tau_{h}(h_1+1))}\cdots
z_{i_{h_2}}^{(\tau_{h_2}\circ\tau_{h}(h))}\right)\\
&=&\frac{F^{\beta}G^{\beta'}}{h!h_1!h_2!}
\sum_{\tau_{h},\tau_{h_1},\tau_{h_2}}z_{k_1}^{(\tau_{h_1}\circ\tau_{h}(1))}\cdots
z_{k_{h}}^{(\tau_{h_2}\circ\tau_{h}(h))}\\
&=&\frac1{h_1!h_2!}\sum_{\tau_{h_1},\tau_{h_2}}\frac{F^{\beta}G^{\beta'}}{h!
} \sum_{\tau_{h}} z_{k_1}^{(\tau_{h_1}\circ\tau_{h}(1))}\cdots
z_{k_{h}}^{(\tau_{h_2}\circ\tau_{h}(h))}\\
&=&\frac1{h_1!h_2!}\sum_{\tau_{h_1},\tau_{h_2}}\frac{F^{\beta}G^{\beta'}}{h!
} \sum_{\tau_{h}'} z_{k_1}^{(\tau'_h(1))}\cdots
z_{k_h}^{(\tau'_h(h))}\qquad \qquad\qquad\mbox{where}\quad\tau'_h=(\tau_{h_1},\tau_{h_2})\circ\tau\\
&=&\frac1{h_1!h_2!}\sum_{\tau_{h_1},\tau_{h_2}}\widetilde
W(z^{(1)},\dots, z^{(h)})\\
&=&\widetilde W(z^{(1)},\dots, z^{(h)}).
\end{eqnarray*}
\end{proof}
\begin{rem}
Usually, the permutation $\tau_{h_1}$ is defined on the set
$\{1,\dots,h_1\}$. However, $\tau_{h}(1)$ may be larger than $h_1$,
since $\tau_h$ is an $h=h_1+h_2$ permutation. In fact, we define
$\tau_{h_1}$ on the set $(\tau_h(1),\dots,\tau_{h}(h_1))$ as
follows. Assume
$\tau_{h}(l_1)<\dots<\tau_{h}(l_i)<\tau_h(l_{i+1})<\dots<\tau_{h}(l_{h_1})$
for $\{l_1,\dots,l_{h_1}\}=\{1,\dots,h_1\}$. Then
 define the map
$\bigwedge:\{\tau_h(1),\dots,\tau_{h}(h_1)\}\rightarrow\{1,\dots,{h_1}\}$
by $\bigwedge(\tau_h(i))=j$, if $i=l_j$, and the inverse of
$\bigwedge$ by $\bigwedge^{-1} (j)=\tau_{h}(i)$. Hence, we have $
\tau_{h_1}\circ\tau_h(i):=\bigwedge^{-1}\circ \tau_{h_1} \circ
\bigwedge (\tau_h(i)). $
\end{rem}
\begin{lem}\label{012101}
\begin{eqnarray*}
f(z)&:=&\frac1{\tilde m!}\frac1{h-1}\sum_{\tau_{\tilde
m}}\sum_{j=1}^{h-1}||z^{\tau_{\tilde
m}(1)}||_1\cdots||z^{\tau_{\tilde m}(j-1)}||_1 ||z^{\tau_{\tilde
m}(j)}||_p||z^{\tau_{\tilde m}(j+1)}||_1\cdots||z^{\tau_{\tilde{m}}
({\tilde{m}})}||_1\\
&=&||(z^{\tilde m})||_{p,1}.
\end{eqnarray*}
\end{lem}
\begin{proof}
 By a direct
calculation,
\begin{eqnarray*}
&&{\tilde m}!(h-1)f(z)\\
&=& \sum_{\tau_{\tilde m}}\sum_{j=1}^{h-1} ||z^{\tau_{\tilde
m}(1)}||_1\cdots||z^{\tau_{\tilde m}(j-1)}||_1
||z^{\tau_{\tilde m}(j)}||_p||z^{\tau_{\tilde m}(j+1)}||_1\cdots||z^{\tau_{\tilde m}({\tilde m})}||_1\\
&=&\sum_{j=1}^{h-1}\left(\sum_{\tau_{\tilde m}} ||z^{\tau_{\tilde
m}(1)}||_1\cdots||z^{\tau_{\tilde m}(j-1)}||_1
||z^{\tau_{\tilde m}(j)}||_p||z^{\tau_{\tilde m}(j+1)}||_1\cdots||z^{\tau_{\tilde m}({\tilde m})}||_1\right)\\
&=&\sum_{j=1}^{h-1}\left(\sum_{j_0=1}^{\tilde m}\sum_{\tau_{\tilde
m}(j)= j_0} ||z^{\tau_{\tilde m}(1)}||_1\cdots||z^{\tau_{\tilde
m}(j-1)}||_1
||z^{(j_0)}||_p||z^{\tau_{\tilde m}(j+1)}||_1\cdots||z^{\tau_{\tilde m}({\tilde m})}||_1\right)\\
&=&\sum_{j=1}^{h-1}\left(\sum_{j_0=1}^{\tilde
m}||z^{(j_0)}||_p\sum_{\tau_{\tilde m}(j)= j_0} ||z^{\tau_{\tilde
m}(1)}||_1\cdots||z^{\tau_{\tilde m}(j-1)}||_1
||z^{\tau_{\tilde m}(j+1)}||_1\cdots||z^{\tau_{\tilde m}({\tilde m})}||_1\right)\\
&=&\sum_{j=1}^{h-1}\left(\sum_{j_0=1}^{\tilde
m}||z^{(j_0)}||_p\left(({\tilde m}-1)!||z^{(1)}||_1
\cdots||z^{(j_0-1)}||_1 ||z^{(j_0+1)}||_1\cdots||z^{({\tilde m})}||_1\right)\right)\\
&=&({\tilde m}-1)!\sum_{j=1}^{h-1}\left(\sum_{j_0=1}^{\tilde
m}\left(||z^{(1)}||_1 \cdots||z^{(j_0-1)}||_1
||z^{(j_0)}||_p||z^{(j_0+1)}||_1\cdots||z^{({\tilde m})}||_1\right)
\right)\\
&=&({\tilde m}-1)!\sum_{j=1}^{h-1}\left({\tilde m}||(z^{{\tilde m}})||_{p,1}\right)\\
 &=&{\tilde m}!(h-1)||(z^{{\tilde m}})||_{p,1}.
\end{eqnarray*}Hence
\begin{eqnarray*}
f(z)=||(z^{{\tilde m}})||_{p,1}.
\end{eqnarray*}
\end{proof}
Let $E$ and $F$ be two complex Banach spaces with norm $||\cdot||_E$
and $||\cdot||_F$, and let $G$ be an analytic map from an open
subset of $E$ into $F$. The first derivative $d_vG$ of $G$ at $v$ is
a linear map from $E$ into $F$, whose induced operator norm is
\begin{equation*}
||d_vG||_{F,E}=\max_{u\neq0}\frac{||d_vG(u)||_F}{||u||_E}.
\end{equation*}
The Cauchy inequality can be stated as follows.
\begin{lem}\label{081907}
Let $G$ be an analytic map from the open ball of radius $r$ around
$v$ in $E$ into $F$ such that $||G||_F\leq M$ on this ball. Then
\begin{equation}
||d_vG||_{F,E}\leq \frac{M}r.
\end{equation}
\end{lem}
\begin{proof}
See details in Lemma A.3 of \cite{P1}.
\end{proof}


\begin{thebibliography}{99}







\bibitem{Bam1} D. Bambusi, Birkhoff normal form for some nonlinear
PDEs., Comm. Math. Phys., $\mathbf{234}$, 253--285 (2003)

\bibitem{Bam3} D. Bambusi, A Birkhoff normal form theorem for some semilinear
PDEs., Hamiltonian dynamical systems and applications, 213--247
(2008)


\bibitem{BDGS}D. Bambusi, J. M. Delort,
B. Gr$\acute{\mbox{e}}$bert and J. Szeftel, Almost global existence
for Hamiltonian semilinear Klein-Gordon equations with small Cauchy
data on Zoll manifolds, Comm. Pure Appl. Math., $\mathbf{60}$,
1665--1690 (2007)

\bibitem{BG} D. Bambusi and B. Gr$\acute{\mbox{e}}$bert, Birkhoff normal form for
partial differential equations with tame modulus, Duke Math. J.,
$\mathbf{135}$, 507--567 (2006)


\bibitem{BN} D. Bambusi and N. N. Nekhoroshev, A property of exponential
stability in nonlinear wave equations near the fundamental linear
mode, Phys. D, $\mathbf{122}$, 73--104 (1998)



\bibitem{BB1} M. Berti and L. Biasco, Branching of Cantor manifolds of elliptic tori and application to PDEs, Comm. Math. Phys. {\bf 305}, 741--796 (2011)










\bibitem{Bour1}
J. Bourgain, Construction of quasi-periodic
solutions for Hamiltonian
    perturbations of linear
     equations and application to nonlinear pde,
     Int. Math. Research Notices, {\bf 11} 475--497  (1994)
     \bibitem{B96} J. Bourgain, Construction of approximative and almost
periodic solutions of perturbed linear Schr$\ddot{\mbox{o}}$dinger
and wave equations, Geom. Funct. Anal., $\mathbf{6}$, 201-230 (1996)

\bibitem{Bour2}
J. Bourgain, On Melnikov's persistence problem,
Math. Research Letters, {\bf 4}  445-458 (1997)





\bibitem{B3} J. Bourgain, Quasi-periodic solutions of Hamiltonian perturbations of 2D linear
Schr$\ddot{\mbox{o}}$dinger equation, Ann. of
 Math., $\mathbf{148}$,
363-439 (1998)
\bibitem{B6}J. Bourgain, On diffusion in high-dimensional Hamiltonian
systems and PDE, J. Anal. Math., {\bf 80}, 1--35 (2000)


\bibitem{Bour5} J. Bourgain,  Recent progress on quasi-periodic lattice
Schr$\ddot{\mbox{o}}$dinger operators and Hamiltonian PDEs, Russian Math. Surveys {\bf 59} :2   231-246 (2004)





\bibitem{Bour6}
 J. Bourgain, Green function estimates for lattice
Schr\"odinger operators  and applications, Annals of Math. Studies, Princeton University Press (2005)








\bibitem{CKSTT} J. Colliander, M. Keel, G. Staffilani, H. Takaota
and T. Tao, Transfer of energy to high frequencies in the cubic
defocusing nonlinear Schr$\ddot{\mbox{o}}$dinger equation, Invent.
Math., {\bf 181}, 39--113 (2010)

\bibitem{C-W} Craig,
W. and  Wayne, C.,  Newton's method and periodic
   solutions of nonlinear wave equation,  Comm. Pure. Appl.
   Math., $\mathbf  {46}$, 1409--1501 (1993)

\bibitem{DS} J. M. Delort and J. Szeftel, Long-time existence for small data
nonlinear Klein-Gordon equations on tori and spheres, Int. Math.
Res. Not., $\mathbf{37}$, 1897--1966 (2004)



\bibitem{Elia2} L. H. Eliasson, A talk in Fudan University, (2007)

\bibitem{E-K} L. H. Eliasson and S. B. Kuksin, KAM for the nonlinear
Schr${\ddot{\mbox {o}}}$dinger equation, Ann. of Math.,
$\mathbf{172}$, 371--435 (2010)

\bibitem{GT}B. Gr$\acute{\mbox{e}}$bert and L. Thomann, KAM for the quantum harmonic oscillator, Comm. Math. Phys., $\mathbf{307}$, 383--427 (2011)

\bibitem{K3}S. B. Kuksin, Hamiltonian perturbations of infinite-dimensional linear
systems with an imaginary spectrum, Funct. Anal. Appl.,
$\mathbf{21}$, 192--205 (1987)

\bibitem{K1}S. B. Kuksin, Nearly integrable infinite-dimensional Hamiltonian systems,
Springer-Verlag, Berlin (1993)

\bibitem{Kuk3}S. B. Kuksin, Fifteen Years of KAM for PDE,  {\it geometry, topology, and
 mathematical physics.}  Amer. Math. Soc. Transl. Ser. 2, 212:  237--258 (2004)




\bibitem{K-P-2} S. B. Kuksin and J. P$\ddot{\mbox{o}}$schel, Invariant Cantor
manifolds of quasi-periodic oscillations for a nonlinear
Schr$\ddot{\mbox{o}}$dinger equation, Ann. of Math., $\mathbf{143}$,
149--179 (1996)






\bibitem{P2} J. P$\ddot{\mbox{o}}$schel, On elliptic lower-dimensional tori in Hamiltonian systems, Math. Z., $\mathbf{202}$, 559--608 (1989)


\bibitem{P1} J. P$\ddot{\mbox{o}}$schel, A KAM-theorem
for some nonlinear partial differential equations, Ann. Scuola Norm.
Sup. Pisa Cl. Sci., $\mathbf{23}$, 119--148 (1996)

\bibitem{W} C. E. Wayne, Periodic and quasi-periodic solutions of nonlinear wave equations via KAM theory, Comm. Math. Phys., $\mathbf{127}$, 479--528 (1990)






\end{thebibliography}
\end{document}